\RequirePackage{fix-cm}

\documentclass[numbook]{svjour3}       
\usepackage{graphicx}
\usepackage{amsmath}
\usepackage{amsfonts}
\usepackage{enumitem} 							
\usepackage[labelfont=bf,tableposition=top]{caption} 	
\usepackage{subfig}									
\usepackage{epstopdf}								
\usepackage{floatrow}
\newfloatcommand{capbtabbox}{table}[][\FBwidth]
\newsavebox{\bigleftbox}
\floatstyle{plaintop}
\restylefloat{table}
\usepackage{color}
\usepackage{mathrsfs}  
\usepackage{multirow}
\usepackage[capitalise]{cleveref}
\usepackage{titlesec}
\usepackage{etoolbox}
\usepackage{lipsum}
\usepackage{chngpage}
\usepackage{wrapfig}
\usepackage{url}

\newcommand{\QQ}{\mathcal{Q}}

\def\B#1{\mbox{\boldmath{$#1$}}}

\newcommand{\norm}[1]{|\!| #1 |\!|}

\newcommand{\dO}{\,\text{d}\Omega}
\newcommand{\dS}{\,\text{d}S}

\newcommand*\closedset[1]{%
   \,
   \vbox{%
     \hrule height 0.5pt
     \kern0.25ex
     \hbox{%
       \kern -0.3em
       \ifmmode#1\else\ensuremath{#1}\fi
       \kern -0.3em
     }
   }
   \,
}

\makeatletter
\newcommand*\rel@kern[1]{\kern#1\dimexpr\macc@kerna}
\newcommand*\widebar[1]{%
  \begingroup
  \def\mathaccent##1##2{%
    \rel@kern{1.2}%
    \overline{\rel@kern{-1.2}\macc@nucleus\rel@kern{-0.2}}%
    \rel@kern{-0.2}%
  }%
  \macc@depth\@ne
  \let\math@bgroup\@empty \let\math@egroup\macc@set@skewchar
  \mathsurround\z@ \frozen@everymath{\mathgroup\macc@group\relax}%
  \macc@set@skewchar\relax
  \let\mathaccentV\macc@nested@a
  \macc@nested@a\relax111{#1}%
  \endgroup
}
\makeatother

\makeatletter
\newsavebox\myboxA
\newsavebox\myboxB
\newlength\mylenA
\newcommand*\widebarQ[2][0.6]{%
    \sbox{\myboxA}{$\m@th#2$}%
    \setbox\myboxB\null
    \ht\myboxB=\ht\myboxA%
    \dp\myboxB=\dp\myboxA%
    \wd\myboxB=#1\wd\myboxA
    \sbox\myboxB{$\m@th\overline{\copy\myboxB}$}
    \setlength\mylenA{\the\wd\myboxA}
    \addtolength\mylenA{-\the\wd\myboxB}%
    \ifdim\wd\myboxB<\wd\myboxA%
       \rlap{\hskip 0.5\mylenA\usebox\myboxB}{\usebox\myboxA}%
    \else
        \hskip -0.5\mylenA\rlap{\usebox\myboxA}{\hskip 0.5\mylenA\usebox\myboxB}%
    \fi}
\makeatother
\usepackage[a4paper,
            left=1in,
            right=1.1in,
            top=1.25in,
            bottom=1.76in]{geometry}

\smartqed

\begin{document}


\title{A DEIM driven reduced basis method for the diffuse Stokes/Darcy model coupled at parametric phase-field interfaces}


\titlerunning{A DEIM driven RB method for Stokes/Darcy coupled at parametric phase-field interfaces}        

\author{Stein K.F. Stoter   \and
        Etienne Jessen      \and
        Viktor Niedens      \and
        Dominik Schillinger}


\institute{S.K.F. Stoter \at
              Institute of Mechanics and Computational Mechanics\\
              Leibniz Univerit\"at Hannover\\
              Appelstrasse 9A, 30167 Hannover \\
              Tel.: +49 511 762 14093\\
              Fax: +49 511 762 19053\\
              \email{Stein.Stoter@ibnm.uni-hannover.de}           
           \and
           E. Jessen \at
              Institute of Mechanics and Computational Mechanics\\
              Leibniz Univerit\"at Hannover\\
           \and
           V. Niedens \at
              Institute of Mechanics and Computational Mechanics\\
              Leibniz Univerit\"at Hannover\\
           \and
           D. Schillinger \at
              Institute of Mechanics and Computational Mechanics\\
              Leibniz Univerit\"at Hannover\\
}


\maketitle

\begin{abstract}
In this article, we develop a reduced basis method for efficiently solving the coupled Stokes/Darcy equations with parametric internal geometry. To accommodate possible changes in topology, we define the Stokes and Darcy domains implicitly via a phase-field indicator function. In our reduced order model, we approximate the parameter-dependent phase-field function with a discrete empirical interpolation method (DEIM) that enables affine decomposition of the associated linear and bilinear forms. In addition, we introduce a modification of DEIM that leads to non-negativity preserving approximations, thus guaranteeing positive-semidefiniteness of the system matrix. We also present a strategy for determining the required number of DEIM modes for a given number of reduced basis functions. We couple reduced basis functions on neighboring patches to enable the efficient simulation of large-scale problems that consist of repetitive subdomains. We apply our reduced basis framework to efficiently solve the  inverse problem of characterizing the subsurface damage state of a complete in-situ leach mining site.

\end{abstract}

\newpage
\section{Introduction}
\label{intro}


Fluid flow through porous media is typically modeled with the Darcy equations. When there are large cracks and voids in the porous medium, then a homogenization of the material into a single permeability tensor is no longer appropriate. The creeping flow in those domains can more accurately be modeled with the Stokes equations. 
The coupling relations between the Darcy and Stokes flow regimes were first studied by Beavers and Joseph \cite{Beavers:67.1}, and later supplemented by Saffman \cite{Saffman:71.1}. 
The resulting Stokes/Darcy equations coupled via Beavers-Joseph-Saffman interface conditions play an important role in many disciplines, for instance modeling groundwater flow \cite{Discacciati:02.1}, petroleum and karst reservoirs \cite{Popov2009,Wu2006}, in-situ leach mining \cite{Forbes1994}), perfusion of blood through tissue \cite{Stoter2017,Mosharaf2019}, filtration devices \cite{Nassehi1998,Hanspal2006}, and chemical reactors \cite{Schneider2005}. 
The wide range of applications has lead to significant interest in methods that computationally approximate the coupled Stokes/Darcy problem, including their numerical analysis and algorithmic treatment \cite{Layton:02.1,Discacciati:07.1,Chen:11.1,Baber:12.1,Discacciati:09.1,Burman:07.1,Badia:09.1,Chidyagwai:09.1,Cao2010}.
Many of these applications are characterized by incomplete or uncertain data in terms of the geometry and topology of the two flow domains (e.g., uncertain subsurface soil characteristics, limited resolution of CT or MRI scans). At the same time, they generally involve a wide range of length scales (e.g., from small rock cavities to a complete mining site, from small capillaries to a full organ). 

In this article, our objective is to enable the efficient simulation of such multiscale systems with parametrically defined free-flow Stokes domains at the lowest scale.
To achieve this objective, we make use of model order reduction by means of reduced basis methods \cite{Chinesta2017,quarteroni2015,Hesthaven2016} on small repetitive subdomains. 
A reduced basis method replaces a computationally expensive high-dimensional finite element discretization of a parametrized partial differential equation (PDE) by 
a small low-dimensional set of basis functions that have high approximation power with respect to the solution manifold of the parametrized problem \cite{Willcox2002,Binev2011}.
We develop a reduced basis method for the coupled Stokes/Darcy model on a variable internal geometry, and model the large-scale system consisting of many repetitive subdomains by coupling together many such reduced basis functions. 
Similar approaches to coupling reduced basis functions on repetitive subdomains have been used in \cite{Antil2010,Antil2011,Baiges2013,Rademacher2014,Ferrero2018,Xiao2019}.

In a reduced basis context, the reduced basis functions themselves are defined on a high-dimensional finite element approximation space. Consequently, the high-resolution finite element mesh must remain fixed. This requirement conflicts with our aim of flexibly handling the geometry of the internal Stokes and Darcy domains. Usually, parametric geometry is handled 
by mapping back to a reference domain \cite{quarteroni2015,Hesthaven2016}, 
but such an approach does not permit changes of topology as this would degenerate the Jacobian of the mapping.
We therefore require a method that is able to model a topologically flexible Stokes domain that may merge and disperse within the Darcy domain. 
Diffuse interface methods \cite{Maury:01.1,Raetz:06.1,Ramiere:07.1,Li:09.1,Lervaag:14.1}, also known as diffuse domain or phase-field methods, offer such a flexible framework for solving coupled boundary value problems on non-boundary fitted meshes. The internal geometry is implicitly represented by a phase-field function, which smoothly transitions from zero to one. 
The phase-field indicator function and its gradient can be leveraged to replace integrals on subdomains or interfaces by weighted volumetric integrals on the complete domain. The resulting phase-field formulation is equivalent to the sharp-boundary interface problem when the width of the diffuse interface (controlled by a characteristic length-scale parameter) limits to zero. Phase-field geometry representations have been widely applied, for instance in growth modeling of tissues and crystals \cite{Oden2010,Golovin1998}, for tracking the evolution of crack patterns \cite{Miehe:10.2,Borden:12.1,Schillinger:15.2,Mikelic:15.1}, for enforcing boundary conditions in imaging data based analysis \cite{Nguyen:17.1,Nguyen:18.1}, for modeling multi-phase flow \cite{Boyer:99.1,Aland:10.1,Teigen:11.1}, for variational image processing and segmentation \cite{chan2001active,gangwar2018robust}, or for modeling phase transition and segregation processes \cite{Anders:12.1,Teigen:09.1,Liu:15.1,Zhao:15.1}.

A crucial aspect of a reduced basis method is the affine decomposition of the linear and bilinear forms. Affine decomposition enables the precomputation of reduced-basis stiffness matrices such that the final reduced order model is completely independent of the size of the original high-fidelity model. However, the (bi)linear forms resulting from the diffuse representation of the internal geometry do not satisfy this criterion. 
This is a common issue for many relevant parametrized PDEs, such as those that feature nonlinearities \cite{Chaturantabut2010,Chaturantabut2012,Drohmann2012,Amsallem2012}, complex material laws \cite{Ryckelynck2009,Kerfriden2011,Fauque2018}, or in general, spatially varying model coefficients \cite{Barrault2004,Grepl2007}. 
For those cases, the discrete empirical interpolation method (DEIM) may be used \cite{Barrault2004,Maday2015}. With DEIM, all non-affine parameter dependent fields are replaced by a low-dimensional interpolation on DEIM modes. For our diffuse representation, this method solves the problem of non-affine parameter dependence, but the interpolation of a domain indicator function is likely to produce Gibbs-type oscillations at regions with high gradients (i.e., the diffuse interfaces). Oscillations necessarily imply negative values and values larger than one, which in turn could produce nonphysical domain representations and unstable system matrices. In this paper, we mitigate this issue by introducing a non-negativity preserving version of DEIM.

Our article is structured as follows: in \cref{sec:Diffuse}, we derive the diffusely coupled Stokes/Darcy equations. We show that all three Beavers-Joseph-Saffman conditions can be treated naturally within a diffuse interface framework. In \cref{sec:DEIM}, we propose a non-negativity preserving variation of the discrete empirical interpolation method for the dimensional reduction of the phase-field geometry representation and show its effectiveness for three benchmark problems. In \cref{sec:RB}, we use the same three benchmark problems to study the relation between the required number of phase-field DEIM modes and the number of reduced basis functions in the reduced order model. In \cref{sec:SolutionMining}, we apply our methodology to efficiently estimate the subsurface flow characteristics of an in-situ leach (ISL) mining site that consists of a large number of repetitive hexagonal units. 
In \cref{sec:Conclusion}, we summarize our work and draw conclusions.

\section{Diffusely coupled Stokes/Darcy equations on parametrically defined domains}
\label{sec:Diffuse}

We consider an incompressible fluid moving through a partially porous medium at velocities that are sufficiently small to neglect the convective components in the material derivative. Steady state equilibrium of the fluid is then described by:
\begin{subequations}\label{StatEq}
\begin{alignat}{3}
-\nabla\cdot  \B{\sigma} &= \B{f} \qquad &&\text{in }\Omega\subset \mathbb{R}^d\,,\\
\nabla \cdot \B{u} &= 0 \qquad &&\text{in }\Omega\,,
\end{alignat}
\end{subequations}
with $\B{\sigma}$ the Cauchy stress tensor, $\B{u}$ the fluid velocity, $\B{f}$ the body force, and $\Omega$ the $d$-dimensional spatial domain.

\subsection{Weak formulation of the coupled Stokes/Darcy equations}
Our porous medium in $\Omega$ contains voids and cracks where the flow is unobstructed. We denote the union of these (potentially disconnected) subdomains $\Omega_S$, and in these subdomains we use the constitutive relation of an incompressible Newtonian fluid:
\begin{subequations}\label{ConstuStokes}
\begin{alignat}{3}
\B{\sigma} &= 2 \mu \nabla^s \B{u} - p \B{I} \qquad &&\text{in }\Omega_S\,, \\
 \B{f} &= \B{0}\qquad &&\text{in }\Omega_S\,,
\end{alignat}
\end{subequations}
with $\mu$ the viscosity and $p$ the pressure field. Hence, in $\Omega_S$, the governing equations are the Stokes equations.

In the remaining domain, $\Omega_D=\Omega\setminus\bar{\Omega}_S$, the interaction between the porous medium and the fluid produces a flow resistance that is assumed to be sufficiently high such that the viscous effects in the Newtonian fluid may be neglected. The closing relations become:
\begin{subequations}\label{ConstuDarcy}
\begin{alignat}{3}
\B{\sigma} &= - p \B{I} \qquad &&\text{in }\Omega_D \,,\\
\B{f} &= -\mu \B{\kappa}^{-1} \B{u} \qquad &&\text{in }\Omega_D\,,
\end{alignat}
\end{subequations}
where $\B{\kappa}$ is the second order permeability tensor. These assumptions lead to the Darcy equations in $\Omega_D$.

The interface that couples the Stokes and Darcy domains is denoted $\Gamma = \bar{\Omega}_S\cap\bar{\Omega}_D$. 
On $\Gamma$, the unit normal vectors $\B{n}_S$ and $\B{n}_D$ point out of the Stokes and Darcy domains, respectively. Across the interface, the solution fields are coupled due to the required balance of mass and balance of linear momentum. Additionally, the porous material introduces a shearing resistance to the fluid on the Stokes domain. These considerations lead to the so-called Beavers-Joseph-Saffman coupling conditions:
\begin{subequations}\label{BJS}
\begin{alignat}{3}
&\B{u}_S\cdot \B{n}_S = -\B{u}_D \cdot \B{n}_D  \qquad &&\text{on }\Gamma\,,\\
& (2\mu \nabla^s \B{u}_S \, \B{n}_S)\cdot \B{n}_S - p_S  = p_D  \qquad &&\text{on }\Gamma\,, \\
&(2\mu \nabla^s \B{u}_S \, \B{n}_S)\cdot \B{T} = -\alpha \, \B{T} \B{u}_S  \qquad &&\text{on }\Gamma\,,
\end{alignat}
\end{subequations}
where subscripts $S$ and $D$ refer to the solution fields in the Stokes and Darcy domains, respectively. The tensor $\B{T} = (\B{I}-\B{n}_S\otimes\B{n}_S)$ is the tangential projector and $\alpha$ is the shear resistance parameter.

On the boundary of $\Omega$, denoted $\partial\Omega$, we permit the following boundary conditions:
\begin{subequations}\label{BCs}
\begin{alignat}{3}
&\B{u}\cdot\B{n} = 0 \qquad && \text{on }\partial\Omega_u \,,\\
&\B{\sigma}\B{n}\cdot \B{T} = \B{0} \qquad && \text{on }\partial\Omega_u \,,\\
&\B{\sigma}\B{n} = t_n \B{n}\qquad &&\text{on }\partial\Omega_p \,.
\end{alignat}
\end{subequations}
The first two conditions represent no-inflow with free-slip and the last condition is a pressure condition.
These specific boundary conditions are chosen because they are valid essential and natural conditions for both the Stokes and the Darcy equations. The inflow condition is chosen homogeneous for simplicity, but it may just as well be set to a non-zero value.

We obtain a weak formulation by multiplying \cref{StatEq} by test functions $\B{v}$ and $q$, integrating over the domain $\Omega$, substituting the constitutive relations \cref{ConstuStokes} and \cref{ConstuDarcy} on their respective domains, performing integration by parts, and by substituting the coupling and boundary conditions of \cref{BJS,BCs}. This leads to the statement:
\begin{subequations}\label{Weak}
\begin{align}
& \text{Find }\B{u},p\in \B{H}_{0}(\Omega_S,\Omega_D,\partial\Omega_u)\times L^2(\Omega) \text{ s.t. } \forall\, \B{v},q\in \B{H}_0(\Omega_S,\Omega_D,\partial\Omega_u) \times L^2(\Omega): \nonumber\\
&\int\limits_{\Omega_S} 2 \mu \nabla^s \B{u} : \nabla^s \B{v}\dO  + \int\limits_{\Omega_D} \mu \B{\kappa}^{-1} \B{u}\cdot \B{v} \dO - \int\limits_{\Omega}  p \, \nabla \cdot \B{v}\dO + \int\limits_\Gamma \alpha \, \B{T}\B{u}_S \cdot \B{v}_S \dS =  \int\limits_{\partial\Omega_p} t_n \B{n} \cdot \B{v} \dS\,, \\
&\int\limits_{\Omega} q\, \nabla \cdot \B{u} \dO = 0 \,,
\end{align}
\end{subequations}
where the function space $\B{H}(\Omega_S,\Omega_D)$ is defined as:
\begin{align}
\begin{split}
    \B{H}_{0}(\Omega_S,\Omega_D,\partial\Omega_u) =& \{ \B{v}\in [L^2]^d : \norm{ \B{v} }^2_{H^1(\Omega_S)} + \norm{ \B{v} }^2_{H(\text{div},\Omega_D)} < \infty  , \, \B{u}\cdot\B{n} = 0 \text{ on }\partial\Omega_u  \} \,.
\end{split}
\end{align}
This space ensures sufficient regularity on the vector valued functions in the Stokes and Darcy domains. 

\subsection{Diffuse interface representation}
Next, we consider a parametrically defined domain $\Omega_S(\beta)$ (and thus also $\Omega_D=\Omega_D(\beta)=\Omega\setminus\Omega_S(\beta)$), where $\beta$ is a point in the parameter space $\mathbb{P}$. The parameter space $\mathbb{P}$ is finite dimensional and encodes in some sense the set of potential geometries of $\Omega_S$. For example, parameters in $\mathbb{P}$ could denote the number, position, size and/or orientation of voids and cracks through the domain $\Omega$. To handle this extensive geometric flexibility, we introduce a diffuse representation of the interface geometry in the weak statement. The strength of such an implicit interface representation is that the computational mesh does not have to fit the internal interface. 

\begin{figure}[t]
    \centering
    \includegraphics[ width=0.57\linewidth]{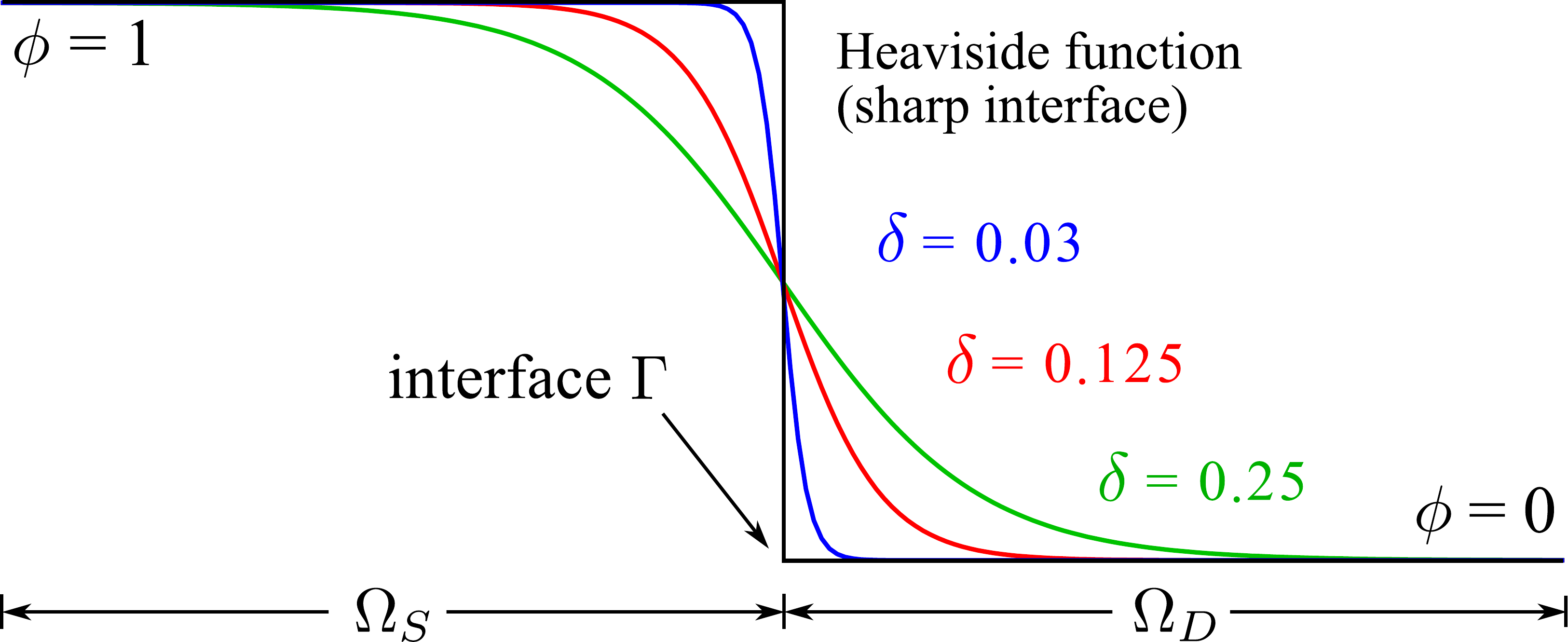}
    \caption{Example one-dimensional phase-fields $\phi$ for different interface widths $\delta$.}
    \label{fig:1dphase-field}
\end{figure}

Let $\phi(\beta)\in \mathcal{C}^1(\Omega)$ be a ``phase-field'' indicator function with values between 1 and 0. The function $\phi(\beta)$ tends to 1 in the Stokes domain and to 0 in the Darcy domain, and it monotonically decreases from 1 to 0 along straight lines that cross a thin region around the Stokes/Darcy interface. The thickness of the diffuse interface is characterized by the parameter $\delta$, as illustrated in \cref{fig:1dphase-field}. Based on this indicator function, volume integrals on $\Omega_S$ or $\Omega_D$ and surface integrals on $\Gamma$ can be approximated by volume integrals on $\Omega$ according to:\\[-0.5cm]
\begin{subequations}
\begin{align}
&\int_{\Omega_S (\beta)} g \dO \approx \int_{\Omega} g \, \phi(\beta) \dO \,, \\
&\int_{\Omega_D (\beta)} g \dO  \approx \int_{\Omega} g \, (1-\phi(\beta)) \dO\,, \\
&\int_{\Gamma (\beta)} g \dS \approx \int_{\Omega} g\, \norm{\nabla{\phi}(\beta)} \dO \,.
\end{align}
\end{subequations}
In the following, we indicate dependencies on $\beta$ only when relevant.

The normal vector $\B{n}_S$ may be approximated as $-\nabla\phi/\norm{\nabla{\phi}}$, such that:
\begin{align}
    \B{T} \approx \B{I} - \frac{\nabla\phi}{\norm{\nabla{\phi}}}\otimes \frac{\nabla\phi}{\norm{\nabla{\phi}}}\,.
\end{align}
Finally, we assume that the solutions $\B{u}$ are such that $\norm{\B{T}\B{u}_D} \ll \norm{\B{T}\B{u}_S}$, whereby $\B{T}\B{u}_S$ on $\Gamma$ can be approximated as twice the solution in the diffuse interface. 

By making use of these approximations in \cref{Weak}, the new weak statement becomes:
\begin{subequations}\label{WeakDif}
\begin{align}
& \text{Find }\B{u},p\in \B{H}_{0}(\phi,\partial\Omega_u)\times L^2(\Omega) \text{ s.t. } \forall\, \B{v},q\in \B{H}_0(\phi,\partial\Omega_u) \times L^2(\Omega): \nonumber\\
    &\int\limits_{\Omega} 2 \phi \mu \nabla^s \B{u} : \nabla^s \B{v} \dO + \!\int\limits_{\Omega} (1-\phi) \mu \B{\kappa}^{-1} \B{u}\cdot \B{v} \dO - \!\int\limits_{\Omega}  p \, \nabla \cdot \B{v} \dO+ \!\int\limits_\Omega  \B{\alpha} \B{u} \,  \cdot \B{v} \dO = \!\int\limits_{\partial\Omega_p} t_n \B{n} \cdot \B{v} \dS \,, \\
&\int\limits_{\Omega} q\, \nabla \cdot \B{u} \dO = 0 \,,
\end{align}
\end{subequations}
where we introduce the new second order tensor field $\B{\alpha}$, defined as:
\begin{align}
  \B{\alpha}  = \frac{2 \alpha}{\norm{\nabla{\phi}}}\, (\norm{\nabla\phi}^2\B{I} - \nabla\phi \otimes \nabla\phi)\,,
\end{align}
and the function space $\B{H}_{0}(\phi,\partial\Omega_u)$, defined as:
\begin{align}
    \B{H}_{0}(\phi,\partial\Omega_u) = \{ \B{v}\in [L^2]^d : \norm{ \phi\B{v} }^2_{H^1(\Omega)} + \norm{ (1-\phi)\B{v} }^2_{H(\text{div},\Omega)} < \infty  , \, \B{u}\cdot\B{n} = 0 \text{ on }\partial\Omega_u  \}\,.
\end{align}

\subsection{Relation with the Brinkman model}
Diffuse geometry representations are often plagued by non-physical behavior in the diffuse interface region. The only way to address this issue is to reduce the length-scale parameter $\delta$ and thereby reduce the interface width. As a result, significant adaptive mesh refinement is required to ensure that the induced error is below a tolerance threshold, which can easily lead to a prohibitive increase in computational cost. We would like to point out, however, that the weak formulation \eqref{WeakDif} corresponds to model equations with physical relevance even for a non-vanishing interface width.


Assuming that the obtained solution pair $(\B{u},p)$ is sufficiently smooth, we can perform integration by parts to arrive at the following statements:
\begin{subequations}\label{WeakBrink}
\begin{align} 
    &\!\!\!\!\int\limits_{\Omega} \! \Big\{\! -\nabla\cdot( 2 \phi \mu \nabla^s \B{u}  - p\B{I} ) + \big[ (1-\phi) \mu \B{\kappa}^{-1}\! +  \B{\alpha} \big] \B{u} \Big\}  \cdot \B{v} \dO  + \!\! \int\limits_{\partial\Omega_p} \! (2 \phi \mu \nabla^s \B{u}\,\B{n}  - p\,\B{n} -t_n \B{n}) \cdot \B{v}\dS = 0 \,, \hspace{-0.2cm}\\
&\!\!\!\!\int\limits_{\Omega} q\, \nabla \cdot \B{u} \dO = 0 \,,
\end{align}
\end{subequations}
which hold for all $ \B{v}\in \B{H}(\phi) $ and all $ q\in L^2(\Omega) $. The corresponding strong form equations are:
\begin{subequations}\label{BrinkStrong}
\begin{alignat}{3}
 -\nabla\cdot( 2\tilde{\mu} \nabla^s \B{u}  - p\B{I} ) &= -\tilde{\B{\kappa}}^{-1}\! \B{u} \qquad &&\text{in }\Omega\,,\\
\nabla \cdot \B{u} &= 0 \qquad &&\text{in }\Omega\,,\\
2\tilde{\mu} \nabla^s \B{u}\,\B{n}  - p\,\B{n}  & = t_n \B{n} \qquad &&\text{on }\partial\Omega_p\,,
\end{alignat}
\end{subequations}
with $\tilde{\mu}=\phi\mu$ and $\tilde{\B{\kappa}}^{-1} = (1-\phi) \mu \B{\kappa}^{-1}\! + \B{\alpha} $. These are the Brinkman equations \cite{Brinkman1949}, which are used frequently for describing a sufficiently fast moving fluid in porous media \cite{Popov2009,Juntunen2010,Mosharaf2019}. The field $\phi$, which was originally a domain indicator, is now a material parameter, and the third Beavers-Joseph-Saffman coupling condition \eqref{BJS} has emerged as an additional orthotropic contribution to the material porosity.

\subsection{Finite element approximation}
\label{ssec:FEspaces}
A finite element formulation based on the weak formulation \eqref{WeakDif} no longer requires a mesh that fits the interface. Instead, we are able to use the same mesh for computing the solution for all domain configurations in the parameter space $\mathbb{P}$. This is a crucial requirement for obtaining a reduced basis representation of the parametric problem later on.

The mixed nature of the equations, however, warrants careful selection of the finite element spaces. This is particularly challenging for the coupled Stokes/Darcy equations: the well-posedness requirements, i.e., the Ladyzhenskaya-Babu\v{s}ka-Brezzi (LBB) conditions, in the Stokes limit are different from those in the Darcy limit~\cite{Boffi2013}. Since our domain is parametric, but our mesh and approximation spaces are fixed, we seek a pressure/velocity pair of elements that is stable for both a pure Stokes problem and a pure Darcy problem. So, we require approximation spaces $\B{\mathcal{V}}^h\times\mathcal{Q}^h$ that are subspaces of the relevant function spaces for all domain configurations:
\begin{subequations}
\begin{alignat}{3}
    &\B{\mathcal{V}}^h &&\subset \B{H}_{0}(\phi(\beta)) \quad \forall\,\beta\in\mathbb{P} \,,\\ 
    &\mathcal{Q}^h &&\subset L^2(\Omega)\,,
\end{alignat}
\end{subequations}
and satisfy the LBB conditions in the limiting cases of both the Stokes and the Darcy equations.
One such pair is the combination of linear nodal elements for the construction of $\mathcal{Q}^h$ and so-called MINI elements for the construction of $\B{\mathcal{V}}^h$ \cite{Mardal2002,Juntunen2010}. The MINI element is a simplicial element and consist of linear interpolation functions enriched with a bubble function \cite{Arnold1984}. We use this combination of elements throughout the remainder of this article. For all simulations, we make use of the FEniCS finite element library \cite{AlnaesBlechta2015a} for computing the large stiffness matrices and the PETSc library for manipulating those matrices \cite{petsc-web-page}.

\section{Non-negativity preserving discrete empirical interpolation of the phase-field}
\label{sec:DEIM}

To ensure efficient formation of the stiffness matrix in the reduced order model, it is imperative that the bilinear and linear forms under consideration depend on the parameter $\beta$ in an affine sense. That is, we require:
\begin{subequations}
\begin{alignat}{3}
    & B((\B{u},p),(\B{v},q);\beta) = B_0((\B{u},p),(\B{v},q)) + \sum\limits_{i=1}^{N_B} \theta^B_i(\beta) B_i((\B{u},p),(\B{v},q)) \,, \label{affinebilinear} \\
    & L((\B{v},q);\beta) = L_0((\B{v},q)) + \sum\limits_{i=1}^{N_L} \theta^L_i(\beta) L_i((\B{v},q)) \,,
\end{alignat}
\end{subequations}
where $B(\cdot,\cdot\,;\beta)$ and $L(\cdot\,;\beta)$ are the bilinear and linear forms corresponding to \cref{WeakDif}. The (bi)linear forms $B_i(\cdot,\cdot)$ and $L_i(\cdot)$ are parameter \textit{independent}, and the parameter dependency of $B(\cdot,\cdot\,;\beta)$ and $L(\cdot\,;\beta)$ is captured through multiplication with the scalar-valued functions $\theta^B_i(\beta)$ and $\theta^L_i(\beta)$. 

\subsection{The discrete empirical interpolation method}
\label{ssec:DEIM}
In our case, the complex parameter dependency of the phase-field $\phi(\beta)$ does not naturally permit a decomposition of the bilinear form as shown in \cref{affinebilinear}. 
To remedy this deficiency, the affine decomposition can be approximated with the discrete empirical interpolation method (DEIM). DEIM produces the following low-dimensional representation of a field:
\begin{align}\label{DEIMrep}
    &f(\B{x};\beta) \approx \sum\limits_{i=1}^{N_f} \theta^f_i(\beta) \tilde{f}_i(\B{x})\,,
\end{align}
where $\tilde{f}_i(\B{x})$ are the parameter independent functions in space (the DEIM interpolation modes) and $\theta^f_i(\beta)$ are the corresponding parameter dependent weighting functions. In summary, the DEIM algorithm for obtaining such a decomposition is a two-step procedure:
\begin{enumerate}
    \item The manifold $\mathcal{M}_{f} = \{ f(\B{x};\beta) : \beta\in\mathbb{P} \}$ is explored and the modes that best represent the complete manifold in some norm are identified. In practice, a discrete representation of the manifold is obtained as a snapshot matrix of solution vectors of the projection of the functions $f(\B{x};\beta)$ onto a finite element space for a sampling of the parameter space. The modes $\tilde{f}_i(\B{x})$ follow from a principal orthogonal decomposition or singular value decomposition of the snapshot matrix.
    \item The functions $\theta^{f}_i(\beta)$ are designed such that \cref{DEIMrep} induces an interpolation of $f(\B{x})$ by the modes $\tilde{f}_i(\B{x})$ at a selection of points. These points are iteratively determined: the $n$-th point is the location where the interpolation of the $n$-th mode by the first $n-1$ modes produces the largest error. When the interpolation locations $\hat{\B{x}}_j$ are defined for $j=1,\cdots, N_f$, the weighting values $\theta^{f}_i(\beta)$ are computed by sampling $f(\hat{\B{x}}_j;\beta)$ and by solving the interpolation system of equations:
    \begin{align}\label{DEIMinterp}
        \sum\limits_{i=1}^{N_f} \theta^{f}_i(\beta) \, \tilde{f}_i(\hat{\B{x}}_j) = f(\hat{\B{x}}_j;\beta)\quad \text{for }j = 1,\cdots, N_f
    \end{align}
\end{enumerate}
The reader is referred to \cite{Barrault2004,Grepl2007} for more details. 

\subsection{A non-negativity preserving alteration}
The following properties of the phase-field $\phi$ and tensor field $\B{\alpha}$ guarantee that the stiffness matrix is positive semi-definite:
\begin{subequations}\label{fieldreqs}
\begin{alignat}{3}
&\phi \geq 0 \,,\\
&(1-\phi) \geq 0 \,,\\
&(\B{\alpha}\, \B{u})\cdot \B{u}  \geq 0\,\,\forall\, \B{u}\in\mathbb{R}^{d}\,.
\end{alignat}
\end{subequations}
These conditions are naturally satisfied by the true phase-field. To guarantee stability of the reduced basis method, they must also be satisfied by the DEIM approximation of these fields after the reconstruction strategy of \cref{DEIMrep}. This is challenging due to the nature of the possible fields $\phi(\beta)$: they are mostly constants of 1 or 0 except at thin regions where they have steep gradients. The non-local support of the DEIM interpolation functions would produce Gibbs-like oscillations in the neighborhood of the interface. 

To remedy this issue, we propose to rewrite the fields \eqref{fieldreqs} in the following form:
\begin{subequations}\label{fieldrews}
\begin{alignat}{3}
&\phi = \xi^2 \,,\\
&(1-\phi) := \psi = \zeta^2 \,,\\
&\B{\alpha} = \sum\limits_{i=1}^{d-1} \B{t}_i\otimes\B{t}_i\,.
\end{alignat}
\end{subequations}
We will now approximate the new fields $\xi$, $\zeta$ and $\B{t}_i$ with \cref{DEIMrep} and the DEIM algorithm. The affine representation of the original fields then follows by expanding the summation multiplication:
\begin{subequations}\label{nonnegreconstruction}
\begin{alignat}{3}
    \phi(\B{x};\beta) \approx &\Big(\sum\limits_{i=1}^{N_\xi} \theta^\xi_i(\beta) \tilde{\xi}_i(\B{x}) \Big)\Big(\sum\limits_{j=1}^{N_\xi} \theta^\xi_j(\beta) \tilde{\xi}_j(\B{x}) \Big) =: \!\!\! \sum\limits_{k=1}^{\frac{N_\xi (N_\xi+1)}{2} } \!\!\!\theta^\phi_k(\beta) \tilde{\phi}_k(\B{x})\,,
\label{nonnegreconstruction_phi} \\
    \psi(\B{x};\beta) \approx &\Big(\sum\limits_{i=1}^{N_\zeta} \theta^\zeta_i(\beta) \tilde{\zeta}_i(\B{x}) \Big)\Big(\sum\limits_{j=1}^{N_\zeta} \theta^\zeta_j(\beta) \tilde{\zeta}_j(\B{x}) \Big) =: \!\!\! \sum\limits_{k=1}^{\frac{N_\zeta (N_\zeta+1)}{2} } \!\!\!\theta^\psi_k(\beta) \tilde{\psi}_k(\B{x})\,,
 \label{nonnegreconstruction_psi} \\
    \B{\alpha}(\B{x};\beta) \approx &\Big(\sum\limits_{i=1}^{N_t} \theta^t_i(\beta) \tilde{\B{t}}_i(\B{x}) \Big)\otimes\Big(\sum\limits_{j=1}^{N_t} \theta^t_j(\beta) \tilde{\B{t}}_j(\B{x}) \Big) =: \!\!\! \sum\limits_{k=1}^{\frac{N_t (N_t+1)}{2} } \!\!\!\theta^\alpha_k(\beta) \tilde{\B{\alpha}}_k(\B{x})\,, \label{nonnegreconstruction_t}
\end{alignat}
\end{subequations}
where, in \cref{nonnegreconstruction_t}, $d$ is taken $2$ for simplicity of notation.

Since the final approximations are constructed from squares, these series expansions will satisfy the requirements of \cref{fieldreqs}. The downside of this approach is the increase of the cost for  the summations involved.

\subsection{Benchmark problems}
In the following, we investigate the relation between the quality of the approximation and the required number of modes for three benchmark problems. For each of these problems, we plot the first four DEIM modes $\tilde{\phi}_k$ and $\tilde{\B{t}}_k$, and we show reconstructions of the true fields $\phi$ and $\B{t}$ for representative parameter points. 
We also qualitatively compare our non-negativity preserving DEIM reconstruction to a standard DEIM approximation of that same field.

\subsubsection{Circular domain with an inwards spiraling channel}
\label{ssec:B1}

\begin{figure}[!b]
    \centering
    \subfloat[Schematic drawing.]{\includegraphics[trim=0 0 0 0,clip, width=0.5\linewidth]{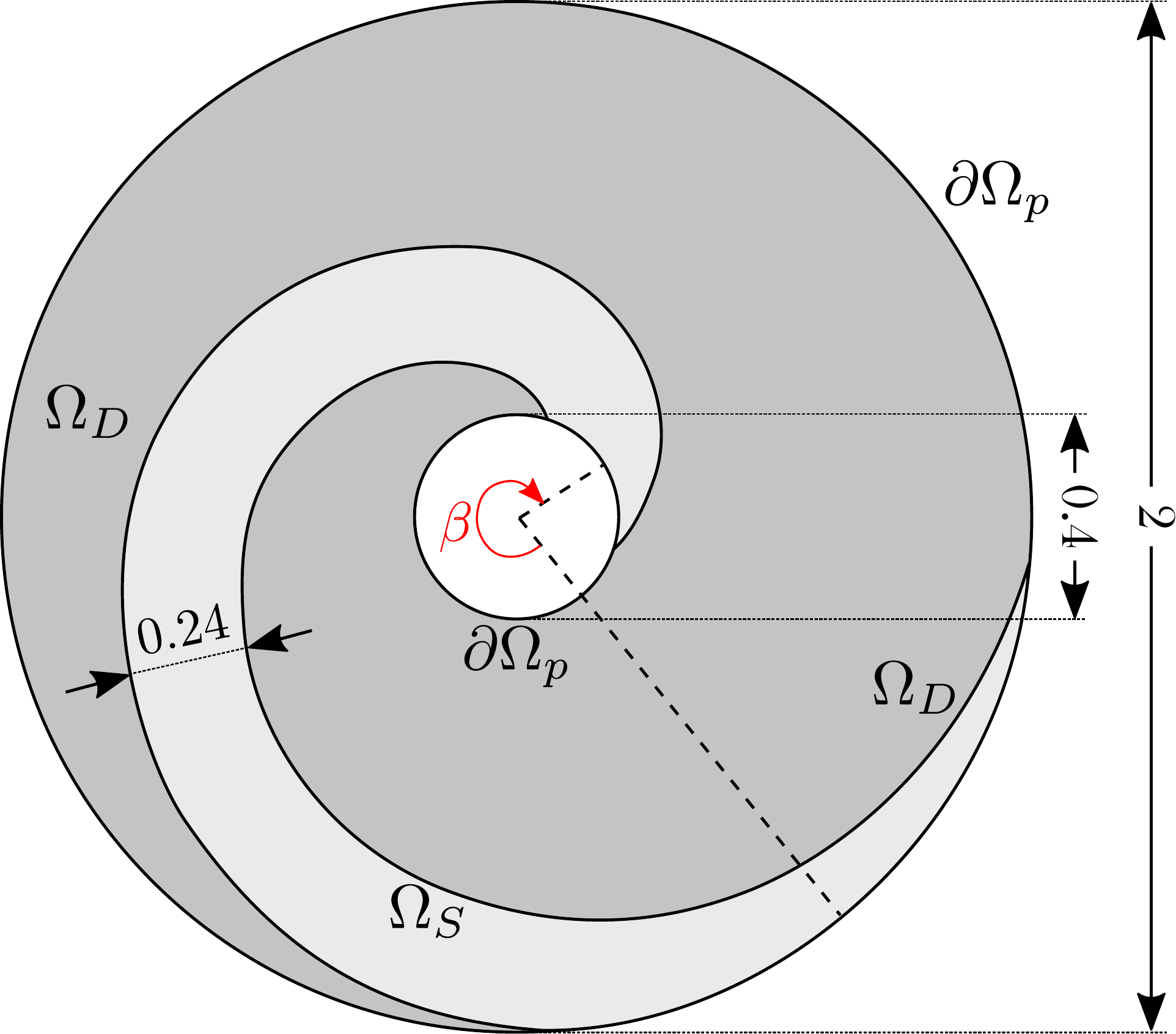} \label{fig:B1Overviewa}}\hspace{0.1cm}
    \subfloat[High fidelity mesh.]{\includegraphics[trim=0 0 0 0,clip, width=0.44\linewidth]{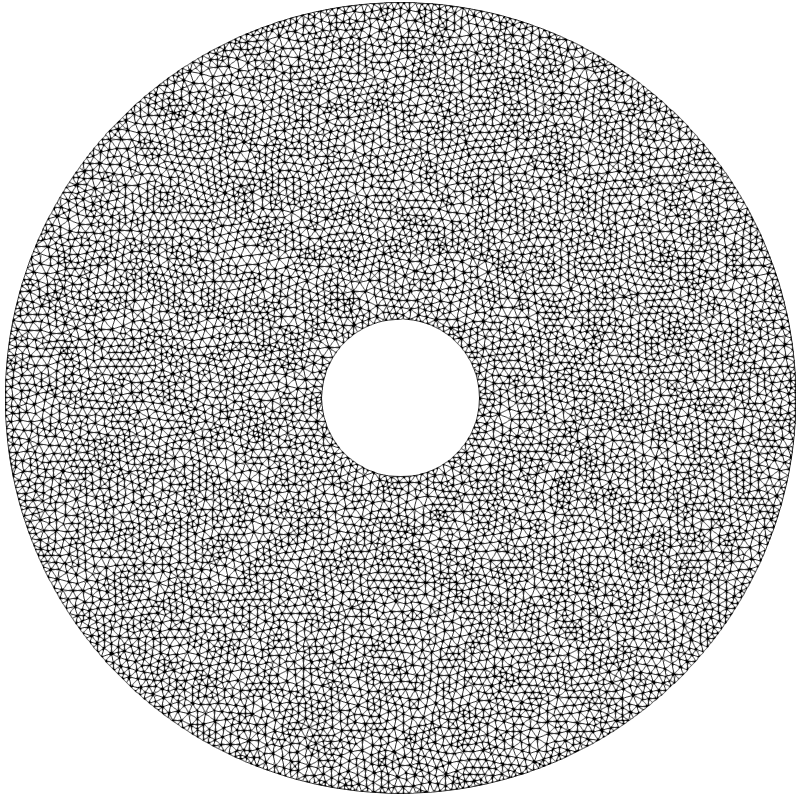} \label{fig:B1Overviewb}}
    \caption{Benchmark 1: a spiraling channel.}
    \label{fig:B1Overview}
\end{figure}

The first benchmark problem consists of a circular domain of porous material with a circular void in the center. The phase-field, representing Stokes flow, is spiraling inwards. The parameter space \mbox{$\mathbb{P}=[180^\circ,540^\circ]$} is one-dimensional and denotes the number of full rotations of the spiral. The material parameters are set to $\kappa=5\cdot 10^{-5} $m$^2$, $\nu=0.5$Pa$\cdot$s and $\alpha=10$Pa$\cdot$s$\cdot$m$^{-2}$, and the boundary pressures $t_n$ are set to 1000Pa and 0Pa at the exterior and interior circles respectively. The problem set-up is illustrated in \cref{fig:B1Overviewa}, and \Cref{fig:B1Overviewb} shows the high-fidelity mesh that we use for all subsequent computations.

\begin{figure}[!b]
    \centering
    \subfloat[First mode, $\tilde{\xi}_1(\B{x})$.\mbox{\hspace{-1cm}}]{\includegraphics[trim=70 40 93 50, clip, width=0.46\linewidth]{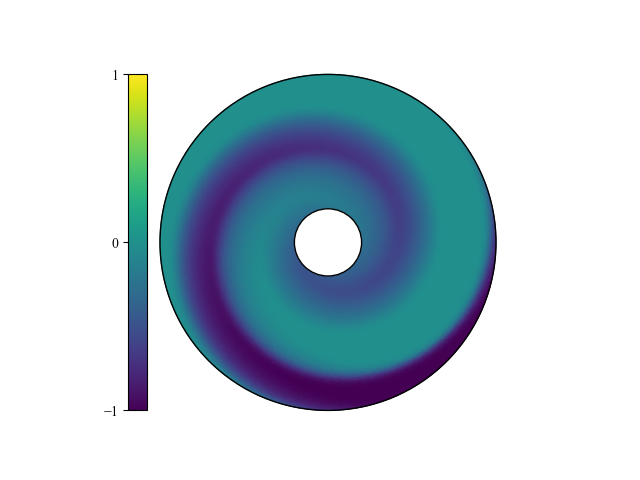} }\hspace{.5cm}
    \subfloat[Second mode, $\tilde{\xi}_2(\B{x})$.]{\includegraphics[trim=108 40 93 50, clip,  width=0.4\linewidth]{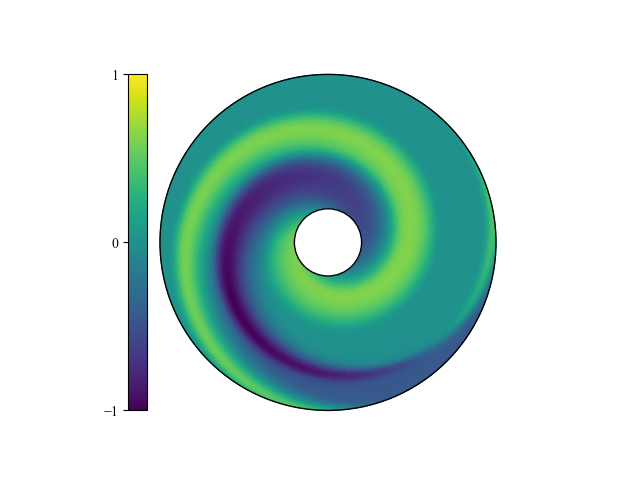} }\\
    \hspace{.75cm}\subfloat[Third mode, $\tilde{\xi}_3(\B{x})$.]{\includegraphics[trim=108 40 93 50, clip, width=0.4\linewidth]{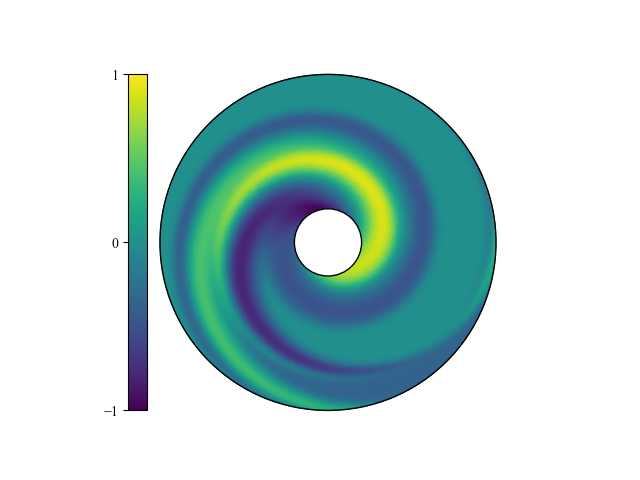}  }\hspace{.5cm}
    \subfloat[Fourth mode, $\tilde{\xi}_4(\B{x})$.]{\includegraphics[trim=108 40 93 50, clip, width=0.4\linewidth]{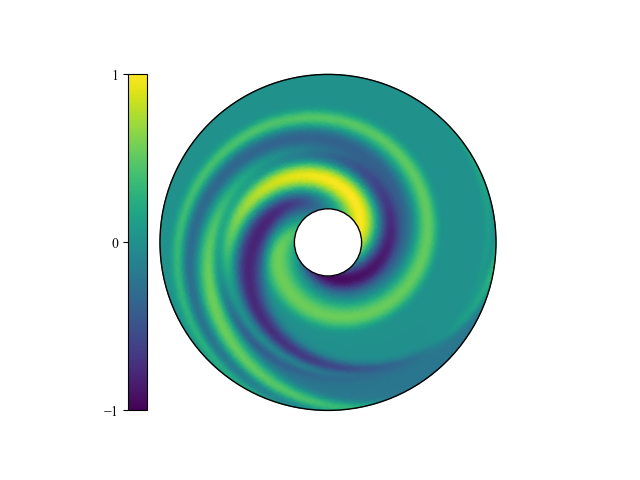} }
    \caption{Benchmark 1: first four DEIM modes of the field $\xi(\beta)$.}
    \label{fig:B1phimodes}
\end{figure}

First, we obtain the DEIM approximation modes of the fields $\xi(\beta)$, $\zeta(\beta)$ and $\B{t}(\beta)$. We use a equidistantly spaced discrete approximation of the parameter space with 1001 samples for the parameter sampling: $\mathbb{P}^h=\{ 180^\circ, (180+\frac{360}{1000})^\circ, \cdots, 540^\circ \}$. For each $\beta\in\mathbb{P}^h$ we construct discrete approximations $\xi^h(\beta)$,  $\zeta^h(\beta)$ and $\B{t}^h(\beta)$ as interpolants of the true fields on the high-fidelity mesh. We collect the solution vectors corresponding to the interpolated fields in a matrix, and performing a singular value decomposition on this snapshot matrix. The first four left-singular-vectors of the fields $\xi(\beta)$ are shown in \cref{fig:B1phimodes}.

\begin{figure}[!t]
    \centering 
    \subfloat[True field.\mbox{\hspace{-1cm}}]{\includegraphics[trim=55 40 103 50, clip, width=0.46\linewidth]{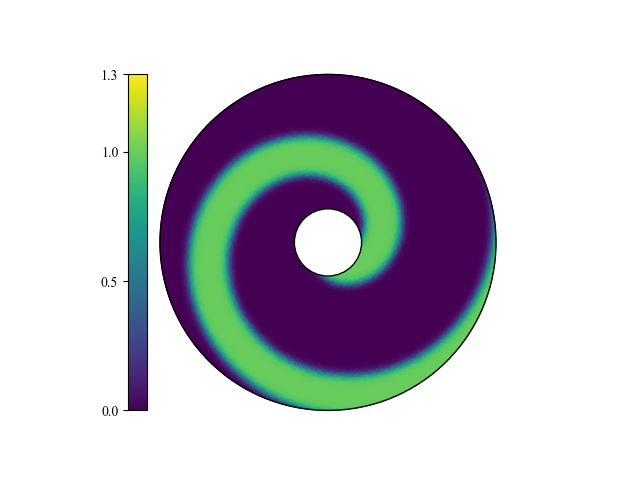} \label{fig:B1phireconstructiona}}\hspace{.5cm}
    \subfloat[Reconstruction with five modes.]{\includegraphics[trim=108 40 93 50, clip, width=0.4\linewidth]{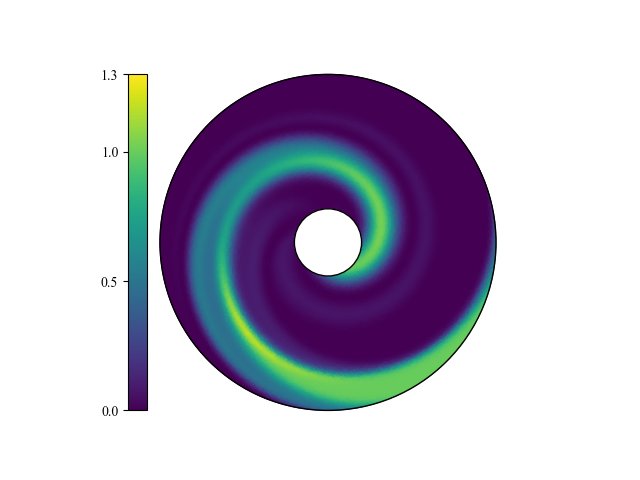} }\\
    \hspace{1.2cm}\subfloat[Reconstruction with ten modes.]{\includegraphics[trim=108 40 93 30, clip, width=0.4\linewidth]{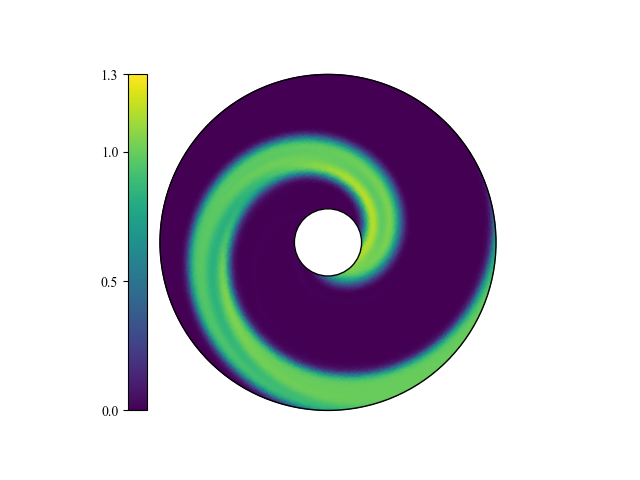}  }\hspace{.3cm}
    \subfloat[Reconstruction with twenty modes.]{\includegraphics[trim=108 40 93 30, clip, width=0.4\linewidth]{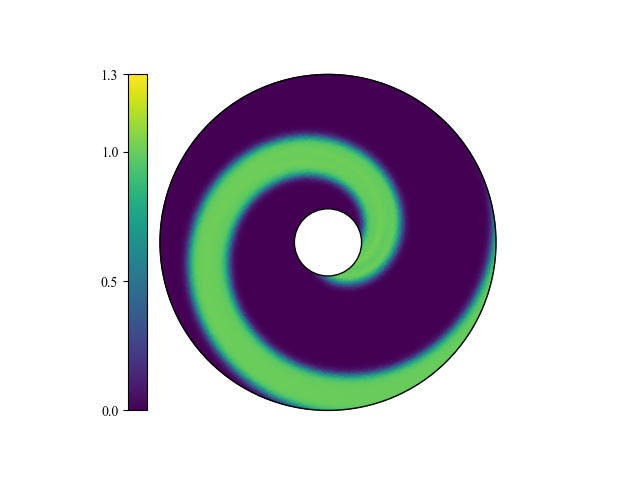}\label{fig:B1phireconstructiond}}
    \caption{Benchmark 1: non-negativity preserving DEIM reconstruction of $\phi(360^\circ)$.}
    \label{fig:B1phireconstruction}
\end{figure}

As an example, we use these interpolation modes to reconstruct the field $\phi$ at the parameter point $\beta=360^\circ$. The original field $\phi(360^\circ)$ is shown in \Cref{fig:B1phireconstructiona}. The other subfigures show the reconstructions with the non-negativity preserving DEIM method for different numbers of modes. As may be observed from the color bar, the reconstruction of $\phi$ is non-negative everywhere, as guaranteed by our reconstruction strategy of \cref{nonnegreconstruction}. When fewer than ten modes are used for the reconstruction then there is still a significant mismatch between the approximated and the true phase-field. From ten onward, the reconstruction clearly resembles the true field, and the difference between the approximated field and the true field is minor when the reconstruction makes use of twenty modes.

\begin{figure}[!b]
    \centering \vspace{-2cm}
    \subfloat[Reconstruction with five modes.\mbox{\hspace{-1cm}}]{\includegraphics[trim=55 40 93 50, clip, width=0.46\linewidth]{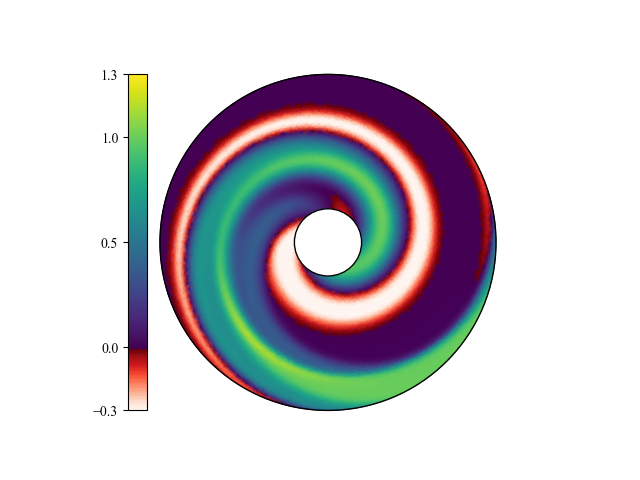}  }\hspace{0.5cm}
    \subfloat[Reconstruction with ten modes.]{\includegraphics[trim=108 40 93 50, clip, width=0.4\linewidth]{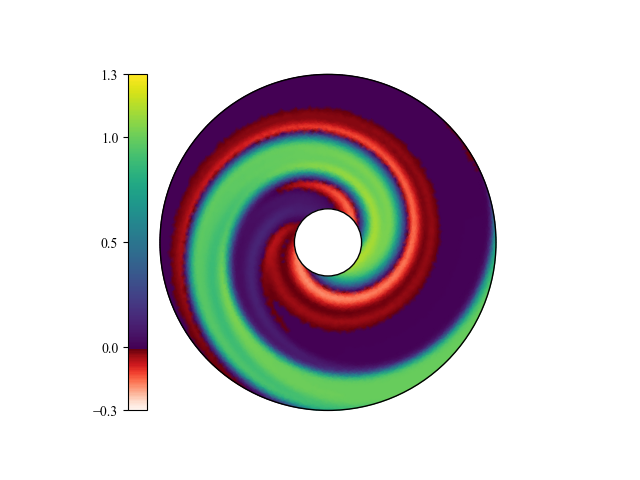} }
    \caption{Benchmark 1: classical DEIM reconstruction of $\phi(360^\circ)$.\\[-0.8cm]}
    \label{fig:B1phireconstructionNegative}
    \vspace{-0.5cm}
\end{figure}

$$$$\newpage
To demonstrate the importance of such a non-negativity preserving reconstruction, we employ the standard DEIM method for the reconstruction of the same field $\phi(360^\circ)$ with the same interpolation modes. The result for five and ten interpolation modes are shown in \cref{fig:B1phireconstructionNegative}. As the color bar indicates, there are banded regions around the interface in the approximated phase-field that are negative. Such undershoots would result in a negative diffusive term, potentially rendering the weak formulation \eqref{WeakDif} unstable.

\begin{figure}[!b]
    \centering
    \subfloat[First mode, $\tilde{\B{t}}_1(\B{x})$.\mbox{\hspace{-1cm}}]{\includegraphics[trim=50 40 93 30, clip, width=0.49\linewidth]{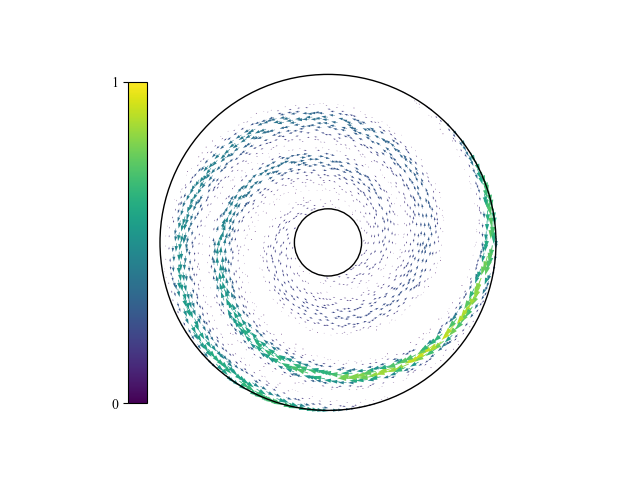} }\hspace{.25cm}
    \subfloat[Second mode, $\tilde{\B{t}}_2(\B{x})$.]{\includegraphics[trim=108 40 93 30, clip,  width=0.4\linewidth]{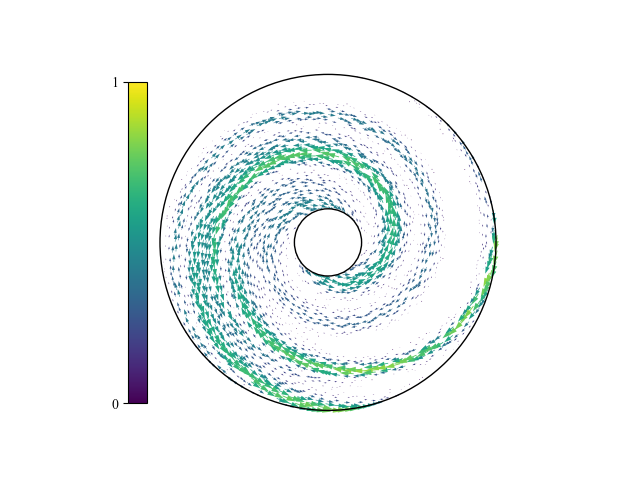} }\\
    \hspace{1.5cm}\subfloat[Third mode, $\tilde{\B{t}}_3(\B{x})$.]{\includegraphics[trim=108 40 93 30, clip, width=0.4\linewidth]{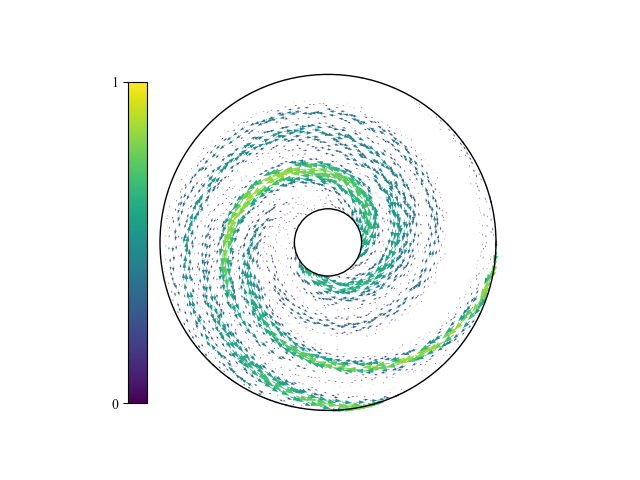}  }\hspace{.25cm}
    \subfloat[Fourth mode, $\tilde{\B{t}}_4(\B{x})$.]{\includegraphics[trim=108 40 93 30, clip, width=0.4\linewidth]{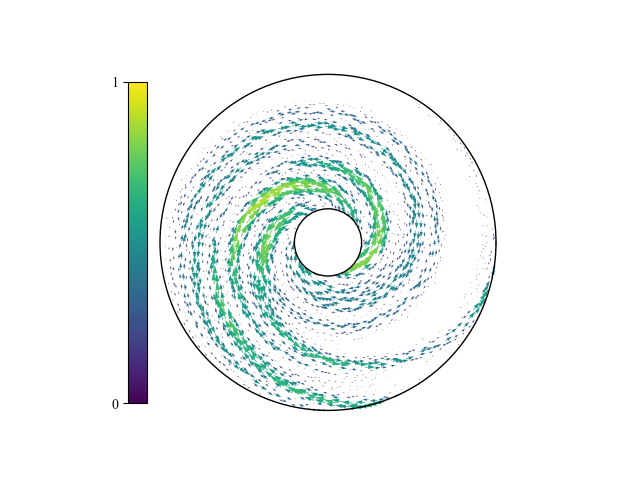} }
    \caption{Benchmark 1: first four DEIM modes of the field ${\B{t}}(\beta)$.}
    \label{fig:B1tmodes}
\end{figure}

\begin{figure}[!t]
    \centering
    \subfloat[True field.\mbox{\hspace{-1cm}}]{\includegraphics[trim=70 40 93 30, clip, width=0.46\linewidth]{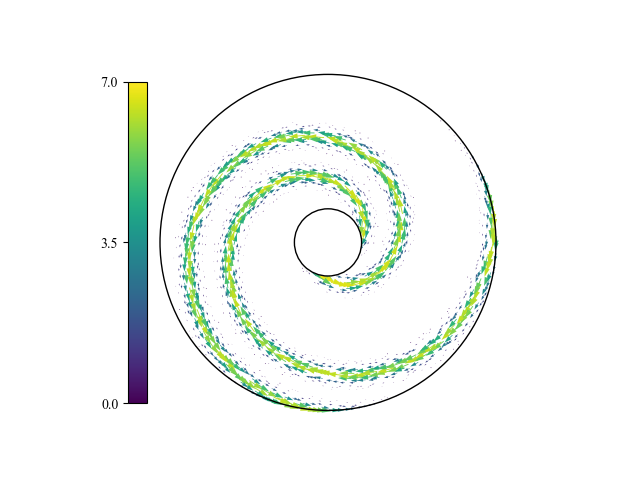} \label{fig:B1treconstructiona}}\hspace{.25cm}
    \subfloat[Reconstruction with five modes.]{\includegraphics[trim=108 40 93 30, clip, width=0.4\linewidth]{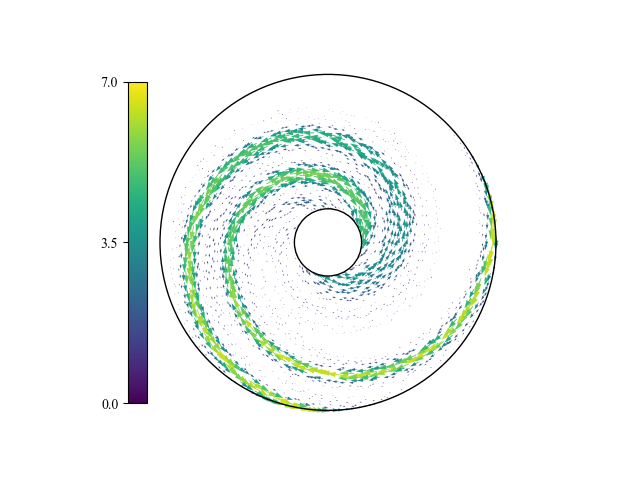} }\\
    \hspace{.75cm}\subfloat[Reconstruction with ten modes.]{\includegraphics[trim=108 40 93 30, clip, width=0.4\linewidth]{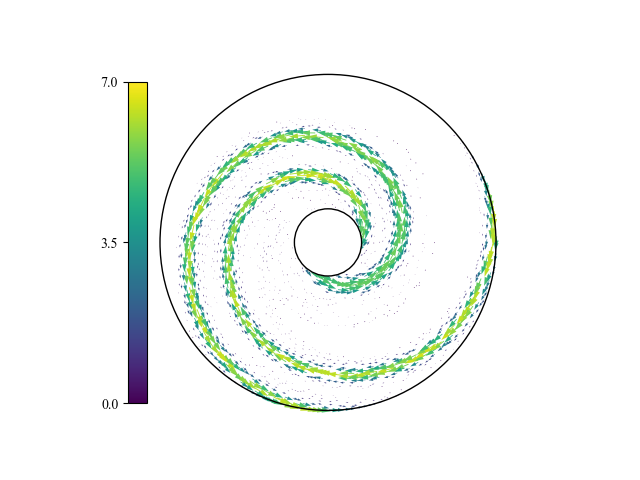}  }\hspace{.25cm}
    \subfloat[Reconstruction with twenty modes.]{\includegraphics[trim=108 40 93 30, clip, width=0.4\linewidth]{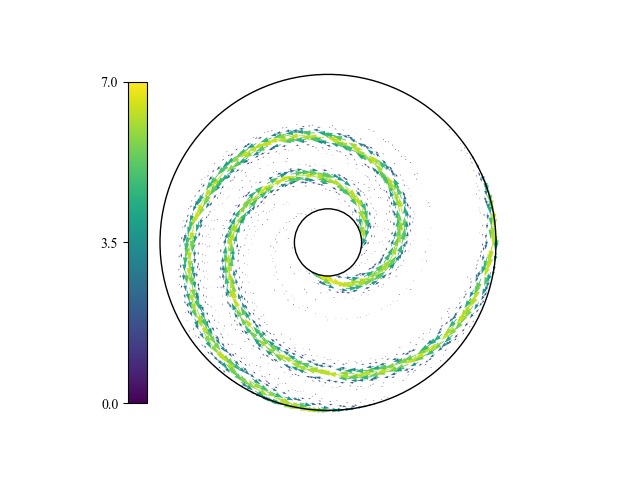} \label{fig:B1treconstructiond}}
    \caption{Benchmark 1: DEIM reconstruction of $\B{t}(360^\circ)$.}
    \label{fig:B1treconstruction}
\end{figure}

Next we focus on the approximation of the field $\B{t}(\beta)$. The first four interpolation modes are shown in \cref{fig:B1tmodes}, and an example reconstruction for the case $\beta=360^\circ$ with various numbers of modes is shown in \cref{fig:B1treconstruction}. We observe a similar trend as for the approximation of $\phi$: ten interpolation modes are required to yield an approximation that resembles the true field, and for twenty interpolation modes the approximation and the true field are qualitatively difficult to distinguish.

The singular values of the singular value decomposition of the different snapshot matrices can be leveraged to quantify the approximation power of a certain number of modes of either one of the fields $\psi$, $\xi$ and $\B{t}$ with respect to the entire (discrete) parameter set. This notion is made rigorous by the Eckart-Young-Mirsky theorem, which relates the singular values to the optimal relative error (in the Frobenius norm) that can be achieved when representing the entire snapshot matrix with the first $n$ interpolation modes
\cite{Eckart1936}. Based on this theorem, we consider the following quality measure:
\begin{align}\label{errormeasure}
    \epsilon(n) = \sqrt{\frac{\sum\limits_{i=n}^{\text{dim}(\mathbb{P}_h)} \sigma_i^2 }{\sum\limits_{i=1}^{\text{dim}(\mathbb{P}_h)} \sigma_i^2 }} \,
\end{align}
where $\sigma_i$ are the singular values.

\begin{figure}[!t]
    \centering
    \includegraphics[trim=0 0 0 0, clip,width=0.6\linewidth]{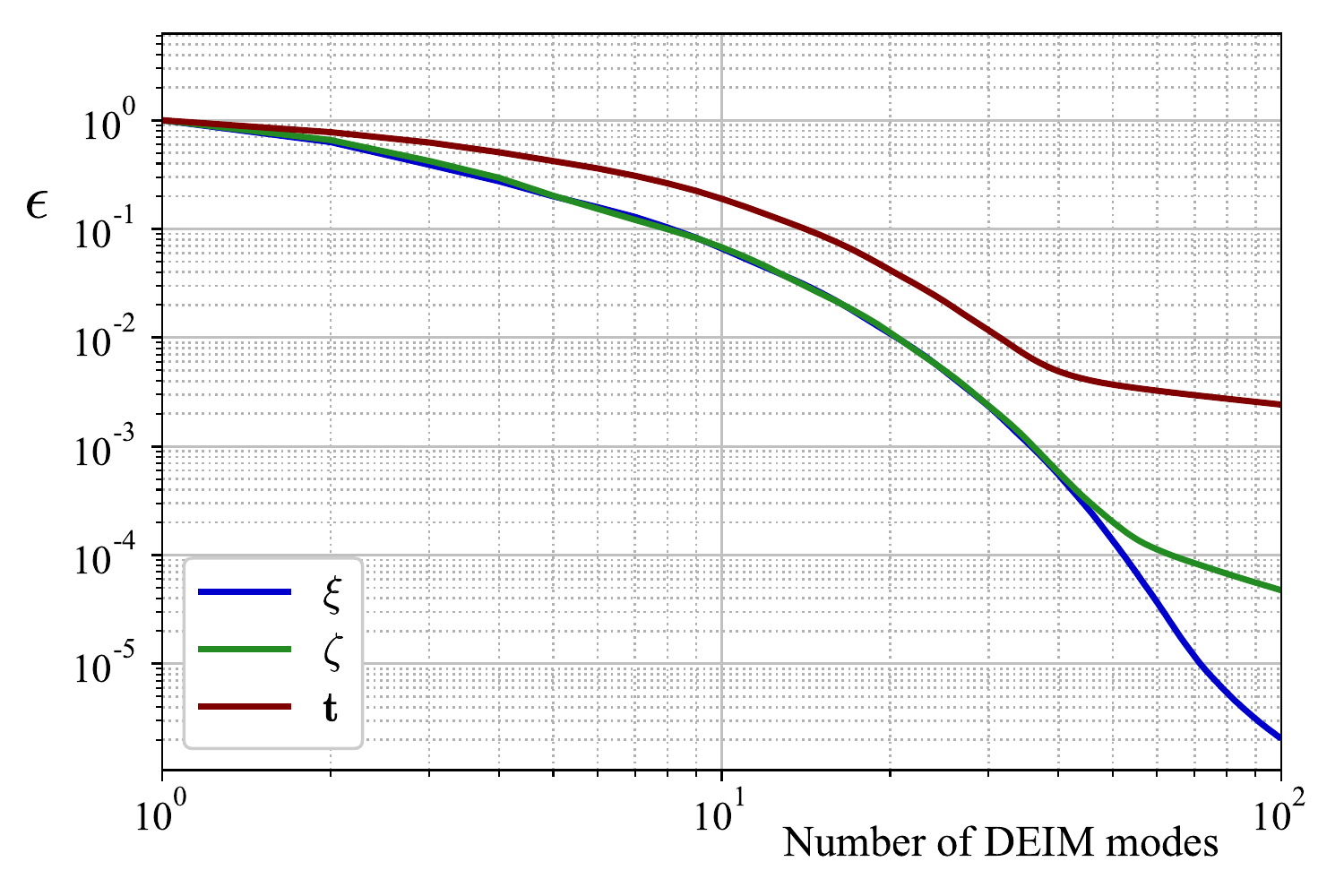}
    \caption{Benchmark 1: convergence of $\epsilon$ as defined in \cref{errormeasure} with the number of modes for $\xi$, $\zeta$ and $\B{t}$.\\[-0.5cm]}\vspace{-0.4cm}
    \label{fig:B1singularvalues}
\end{figure}

\Cref{fig:B1singularvalues} shows $\epsilon(n)$ for the different fields $\xi$, $\psi$ and $\B{t}$. The overlap of the lines  $\xi$ and $\psi$ for the majority of the 
for graph indicates that the approximation of $\xi$ and $\psi$ may be expected to perform equally well. This is to be expected since they are closely related. The graphs also shows that $\B{t}$ is the most challenging field to represent accurately with a finite-dimensional linear approximation. Still, with a total of $\scriptstyle{\sim}$30 modes, the measure for the relative error drops below 1\%, which is more than sufficient for applications where interface geometries are described by a phase-field. The 20-mode reconstruction of $\phi$, shown in \cref{fig:B1phireconstructiond}, corresponds to an $\epsilon_\xi$-value of approximately 0.01, and the 20-mode reconstruction of $\B{t}$, shown in \cref{fig:B1treconstructiond}, corresponds to an $\epsilon_t$-value of 0.03.
The difference between $\epsilon(n)$ for the fields $\xi$, $\psi$ and $\B{t}$ implies that a computationally optimal reduced basis method will require different numbers of interpolation modes to represent the different phase-field quantities. We explore this concept in more depth in \cref{sec:RB}.

\subsubsection{Rotating channel}
\label{ssec:B2}

\begin{figure}[!b]
\vspace{-0.2cm}
    \centering
    \subfloat[Schematic drawing.]{\includegraphics[trim=40 50 20 50,clip, width=0.56\linewidth]{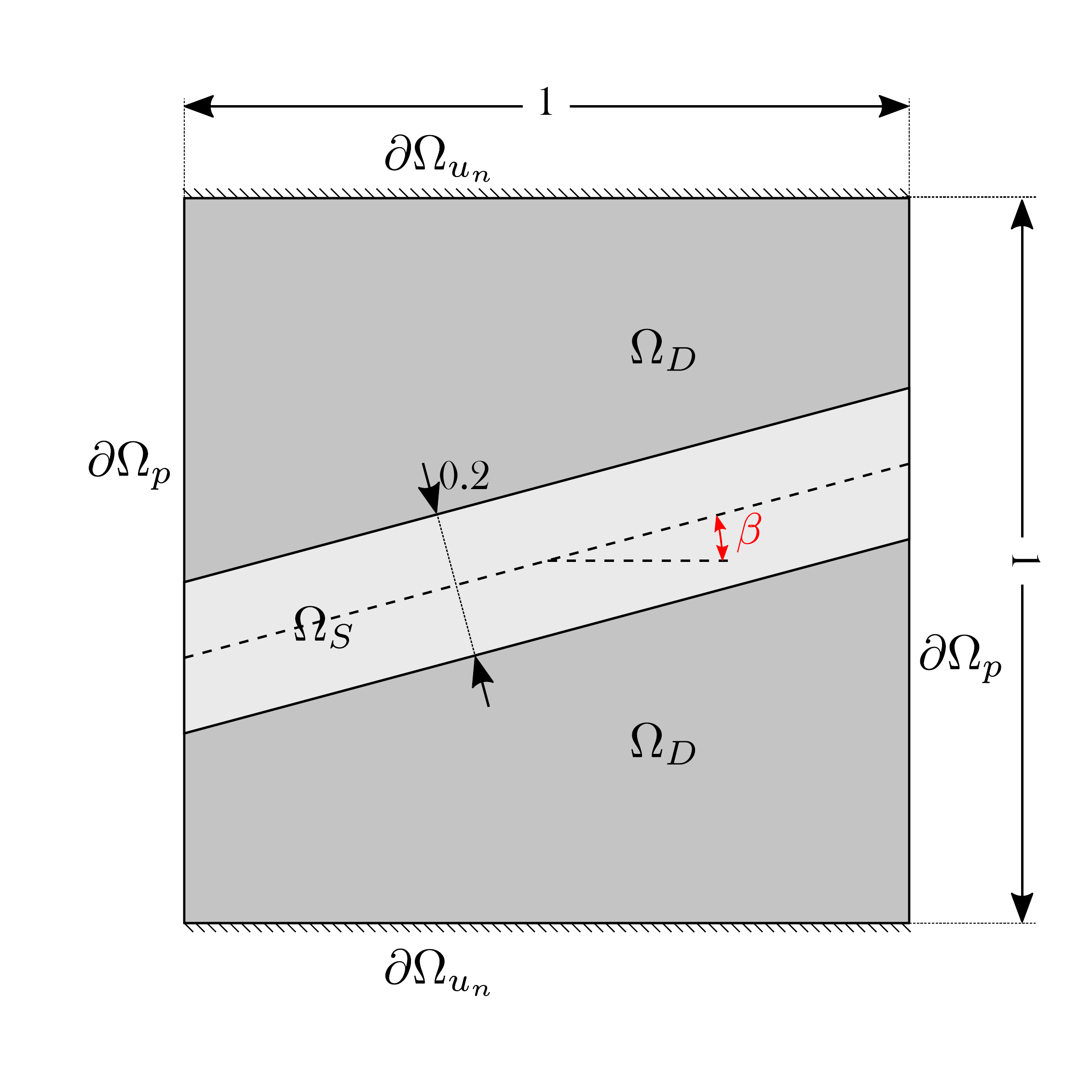}\label{fig:B2Overviewa}}\hspace{0.1cm}
    \subfloat[High fidelity mesh.]{\includegraphics[trim=0 -55 0 0,clip, width=0.42\linewidth]{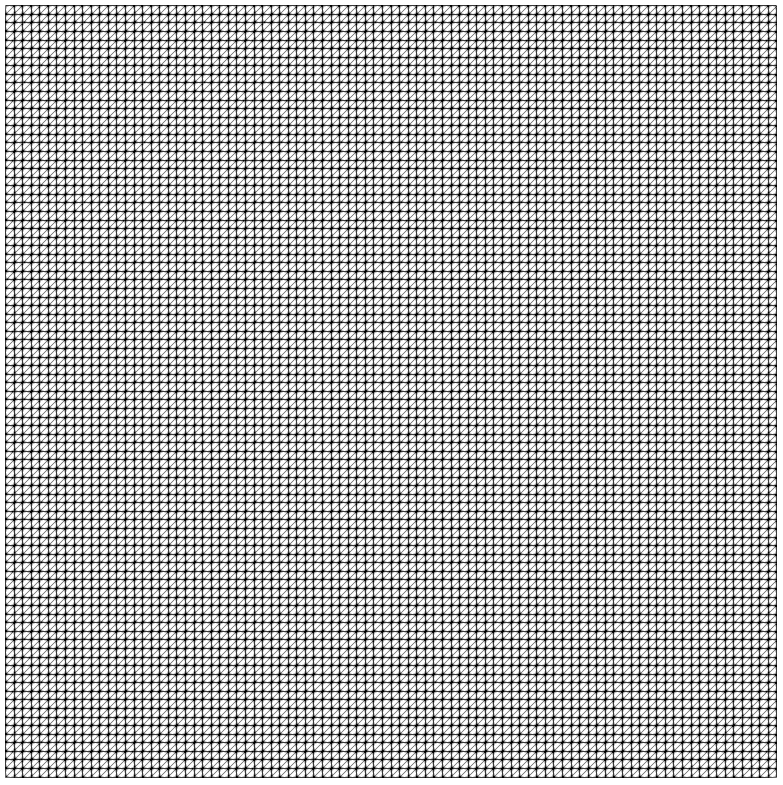}\label{fig:B2Overviewb}}
    \centering
    \caption{Benchmark 2: a rotating channel.\\[-0.7cm]}\vspace{-0.5cm}
    \label{fig:B2Overview}
\end{figure}

Next, we consider an angled straight channel of Stokes flow that is embedded in a square domain of porous material. The parameter is the angle of rotation of the channel, such that the parameter space $\mathbb{P}=[0^\circ,180^\circ]$ is one-dimensional. We set $\kappa=5\cdot 10^{-5} $m$^2$, $\nu=0.5$Pa$\cdot$s and $\alpha=10$Pa$\cdot$s$\cdot$m$^{-2}$, and $t_n=1000$ at the left boundary and $t_n=0$ at the right boundary. The problem is illustrated in \cref{fig:B2Overviewa}, and the corresponding high-fidelity mesh is illustrated in \cref{fig:B2Overviewb}. We may anticipate that the DEIM approximation is very ineffective for this type of problem: for different parameter points $\beta$, i.e., different angles of the rotating channel, the phase-field is largely uncorrelated. Widely varying entries in the $\phi(\beta)$-manifold means that the Kolmogorov $n$-width of the manifold is large and that linear dimensional reduction requires many modes to produce accurate approximations. The purpose of this benchmark is to investigate the impact of this problem. We are interested in the required number of DEIM modes such that we are still able to obtain acceptable reconstructions of the internal geometry.

\begin{figure}[!b]
    \centering
    \subfloat[First mode, $\tilde{\xi}_1(\B{x})$.\mbox{\hspace{-1cm}}]{\includegraphics[trim=55 30 100 50, clip, width=0.535\linewidth]{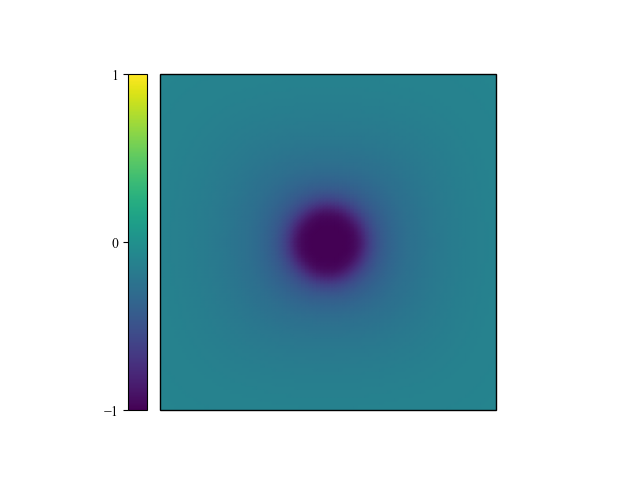} }
    \subfloat[Second mode, $\tilde{\xi}_2(\B{x})$.]{\includegraphics[trim=110 30 100 50, clip,  width=0.44\linewidth]{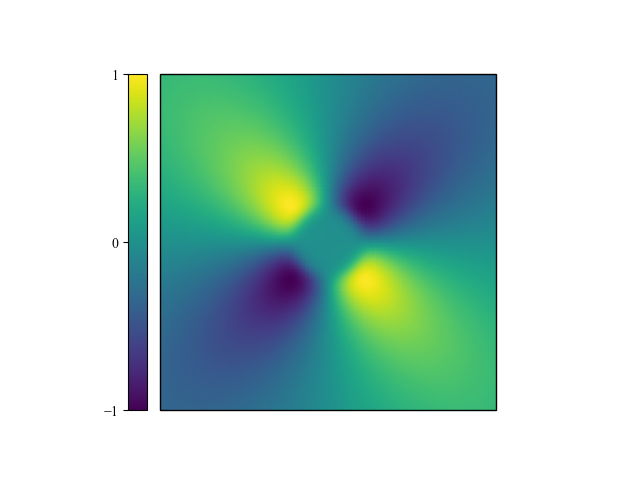} }\\\hspace{1.5cm}
    \subfloat[Third mode, $\tilde{\xi}_3(\B{x})$.]{\includegraphics[trim=110 30 100 50, clip, width=0.44\linewidth]{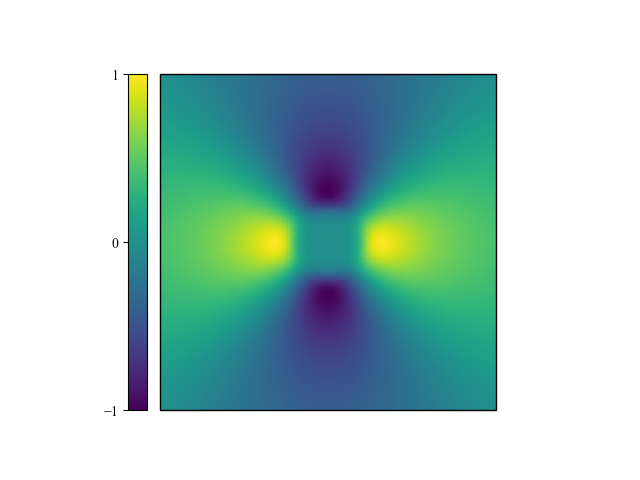}  }
    \subfloat[Fourth mode, $\tilde{\xi}_4(\B{x})$.]{\includegraphics[trim=110 30 100 50, clip, width=0.44\linewidth]{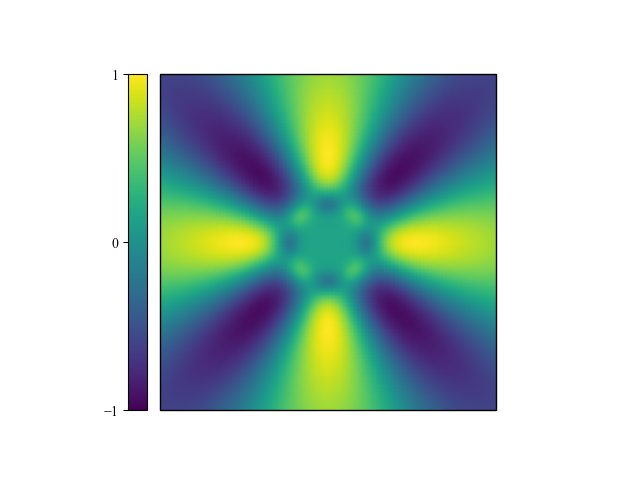} }
    \caption{First four DEIM modes of the field $\xi(\beta)$ of benchmark 3.}
    \label{fig:B2phimodes}
\end{figure}

\begin{figure}[t!]
\vspace{-0.1cm}
    \centering
    \subfloat[True field.\mbox{\hspace{-1cm}}]{\includegraphics[trim=55 40 100 50, clip, width=0.49\linewidth]{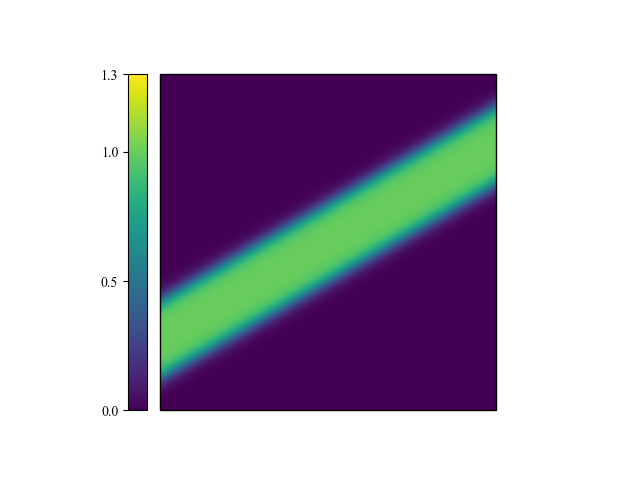}}
    \subfloat[Reconstruction with ten modes.]{\includegraphics[trim=110 40 102 55, clip, width=0.40\linewidth]{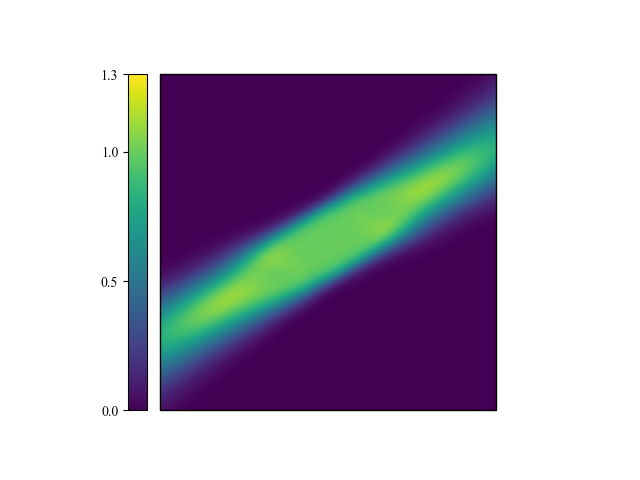} }\\ \hspace{1.3cm}
    \subfloat[Reconstruction with twenty modes.]{\includegraphics[trim=110 40 102 55, clip, width=0.40\linewidth]{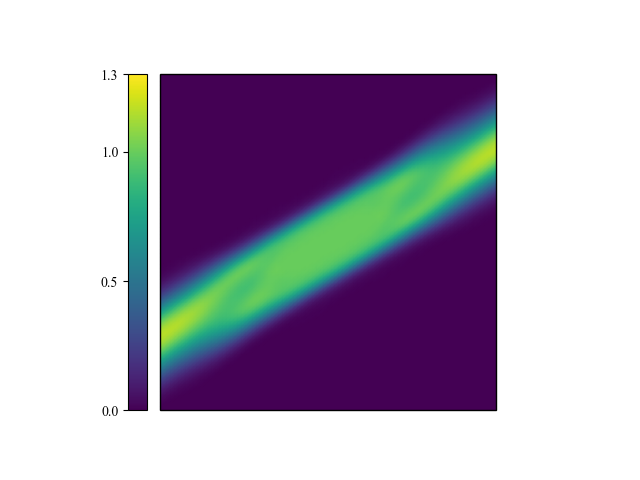} 
    \label{fig:B2phireconstructionc} }
    \subfloat[Reconstruction with forty modes.]{\includegraphics[trim=114 40 98 55, clip, width=0.40\linewidth]{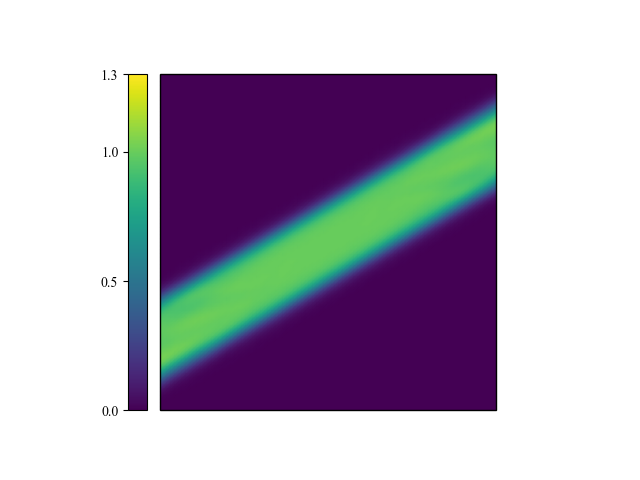} \label{fig:B2phireconstructiond}}
    \caption{Benchmark 2: non-negativity preserving DEIM reconstruction of $\phi(30^\circ)$.}
    \label{fig:B2phireconstruction}
\end{figure}

\begin{figure}[!b]
\vspace{-1.8cm}
    \centering
    \subfloat[Reconstruction with ten modes.\mbox{\hspace{-1cm}}]{\includegraphics[trim=55 40 100 50, clip, width=0.49\linewidth]{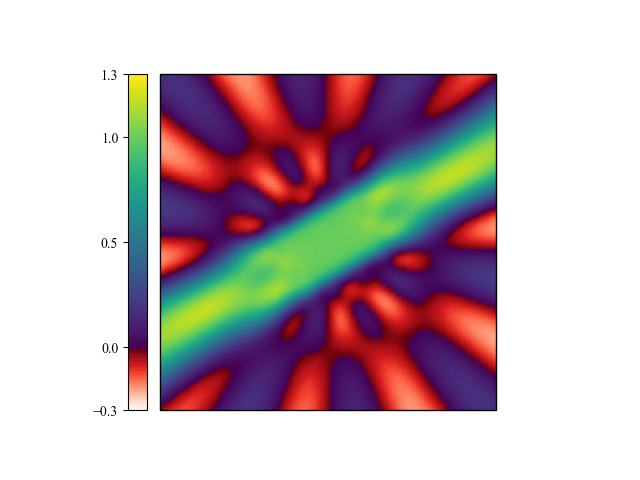}  }
    \subfloat[Reconstruction with twenty modes.]{\includegraphics[trim=110 40 102 50, clip, width=0.40\linewidth]{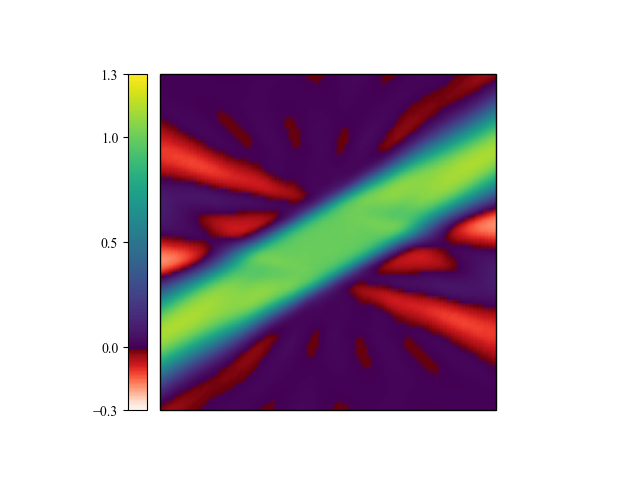} }
    \caption{Benchmark 2: classical DEIM reconstruction of $\phi(30^\circ)$.\\[-1.0cm]}
    \label{fig:B2phireconstructionNegative}
\vspace{-0.6cm}
\end{figure}

\begin{figure}[H]
    \centering
    \subfloat[First mode, $\tilde{\B{t}}_1(\B{x})$.\mbox{\hspace{-1cm}}]{\includegraphics[trim=55 30 105 40, clip, width=0.54\linewidth]{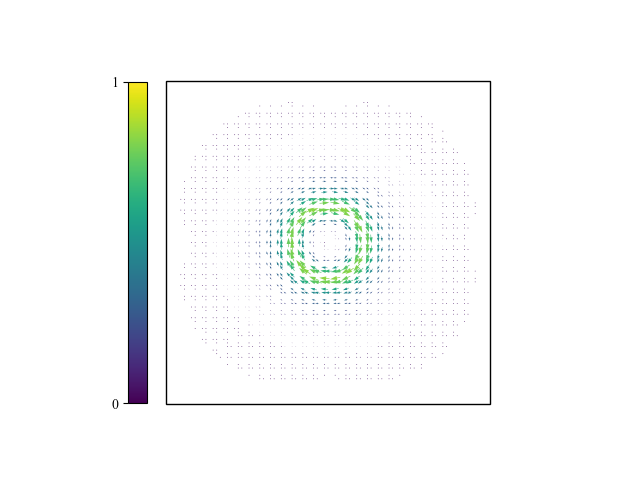} }
    \subfloat[Second mode, $\tilde{\B{t}}_2(\B{x})$.]{\includegraphics[trim=110 30 105 40, clip,  width=0.44\linewidth]{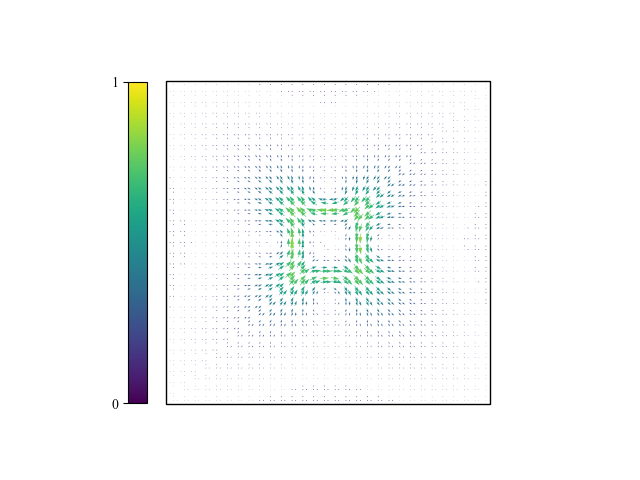} }\\\hspace{1.5cm}
    \subfloat[Third mode, $\tilde{\B{t}}_3(\B{x})$.]{\includegraphics[trim=110 30 105 40, clip, width=0.44\linewidth]{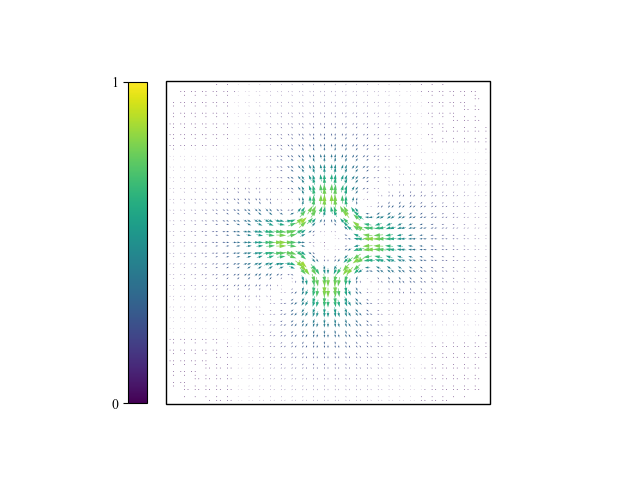}  }
    \subfloat[Fourth mode, $\tilde{\B{t}}_4(\B{x})$.]{\includegraphics[trim=110 30 105 40, clip, width=0.44\linewidth]{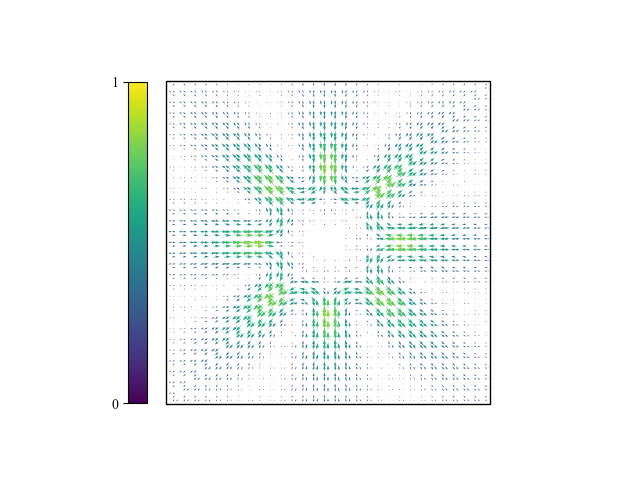} }
    \caption{Benchmark 2: first four DEIM modes of the field ${\B{t}}(\beta)$.}
    \label{fig:B2tmodes}
\end{figure}

\begin{figure}[H]
    \centering
    \subfloat[True field.\mbox{\hspace{-1cm}}]{\includegraphics[trim=55 30 105 40, clip, width=0.54\linewidth]{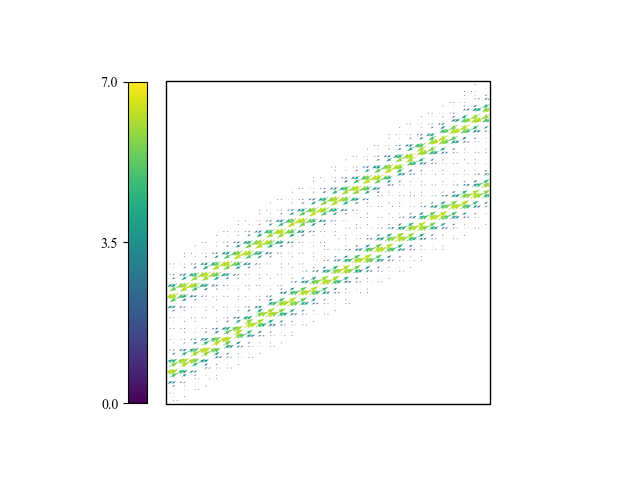}}
    \subfloat[Reconstruction with ten modes.]{\includegraphics[trim=110 30 105 40, clip, width=0.44\linewidth]{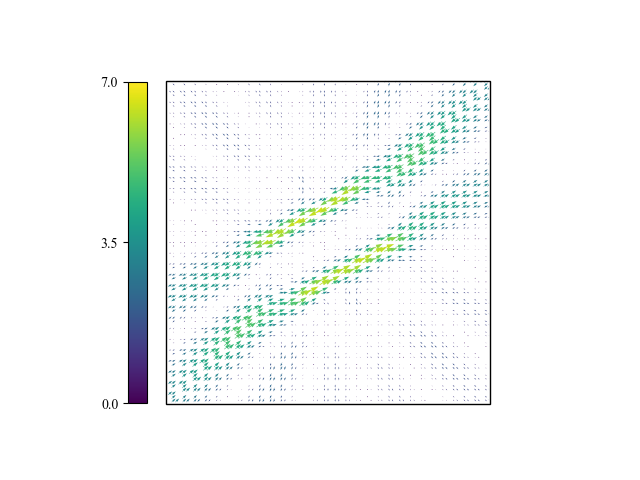} }\\\hspace{1.5cm}
    \subfloat[Reconstruction with twenty modes.]{\includegraphics[trim=110 30 105 40, clip, width=0.44\linewidth]{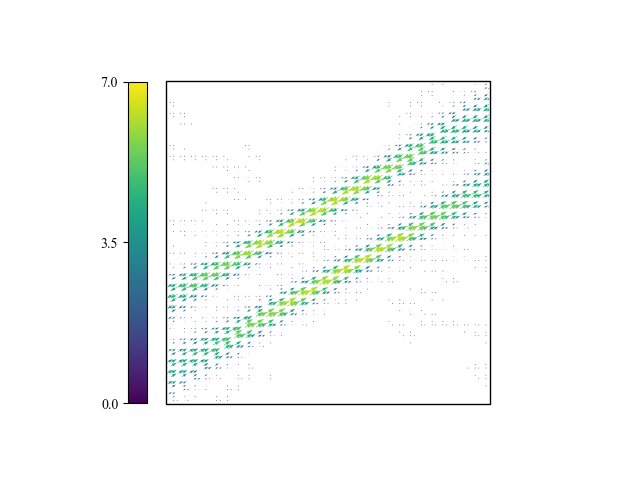}  \label{fig:B2treconstructionc}}
    \subfloat[Reconstruction with forty modes.]{\includegraphics[trim=110 30 105 40, clip, width=0.44\linewidth]{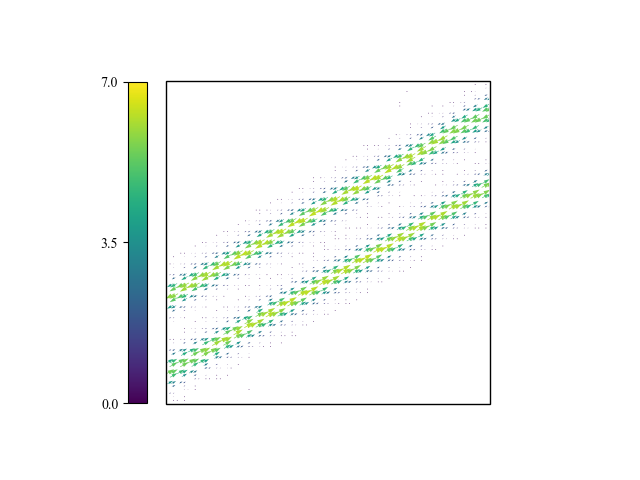} \label{fig:B2treconstructiond}}
    \caption{Benchmark 2: DEIM reconstruction of $\B{t}(30^\circ)$.}
    \label{fig:B2treconstruction}
\end{figure}

We again sample the parameter space with 1001 equidistantly spaced samples, producing $\mathbb{P}^h = \{ 0^\circ, \frac{180}{1000}^\circ, \cdots ,180^\circ \}$ as a discrete approximation of the full parameter space. \Cref{fig:B2phimodes} shows the DEIM modes of the field $\xi$. Examples the reconstruction of the field $\phi(30^\circ)$ are shown in \cref{fig:B2phireconstruction} for the non-negativity preserving DEIM approximation and in \cref{fig:B2phireconstructionNegative} for the standard DEIM approximation. We again observe large patches with negative phase-field values in \cref{fig:B2phireconstructionNegative}, and no regions with negative values in \cref{fig:B2phireconstruction}.
\Cref{fig:B2tmodes} shows the DEIM modes for $\B{t}$, and example reconstructions of $\B{t}(30^\circ)$ are shown in \cref{fig:B2treconstruction}. The reconstructed fields still exhibit quite significant errors away from the center of the channel. Where 20 modes were sufficient to approximate an example phase-field in the previous benchmark, \cref{fig:B2phireconstructionc,fig:B2treconstructionc} still shows significant deviations from the true fields. The approximation becomes sufficiently accurate only when 40 modes are used for the reconstruction, as illustrated in \cref{fig:B2phireconstructiond,fig:B2treconstructiond}.

\begin{figure}[H]
    \centering
    \includegraphics[trim=0 0 0 0, clip,width=0.6\linewidth]{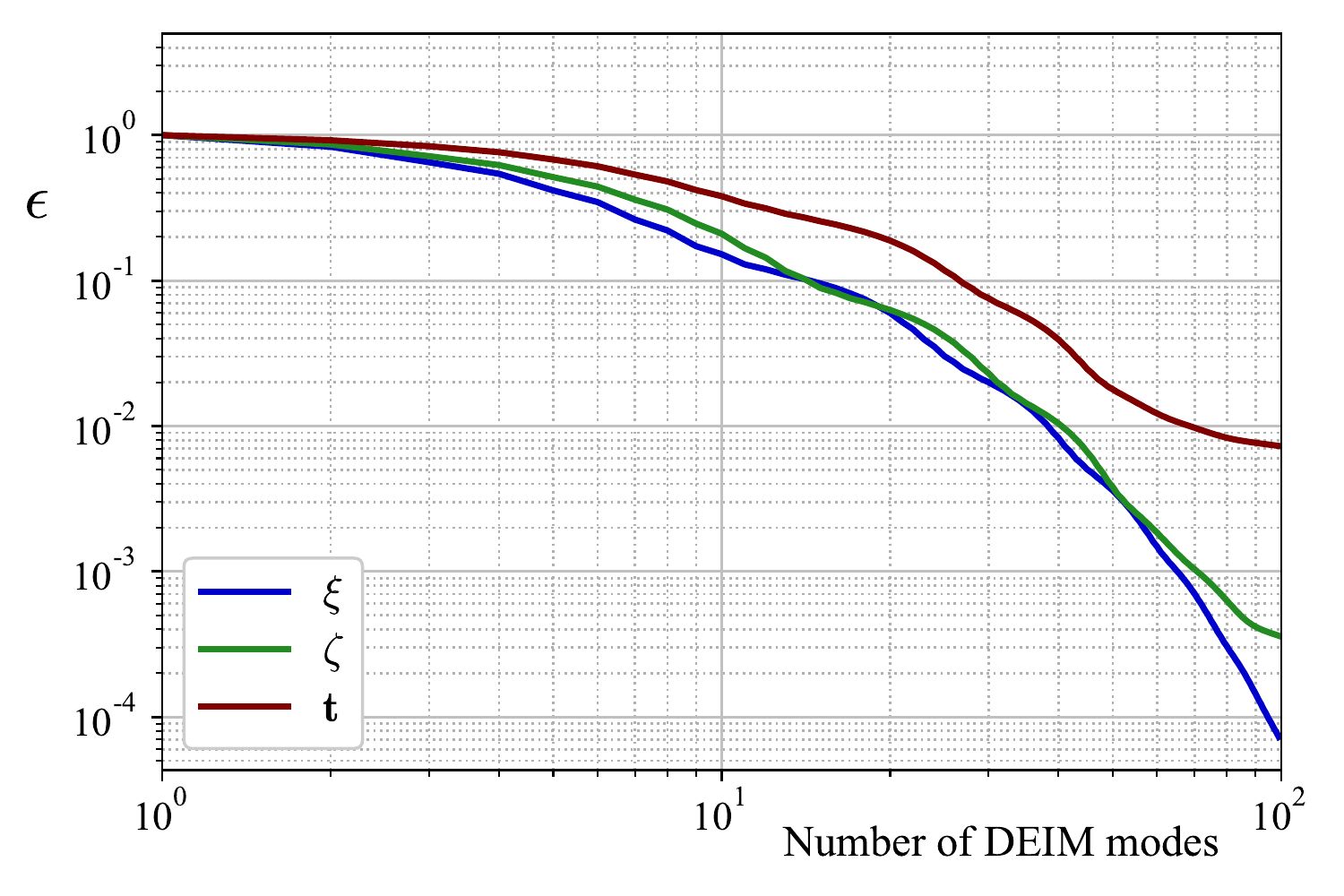}
    \caption{Benchmark 2: convergence of $\epsilon$ for $\xi$, $\zeta$ and $\B{t}$ with the number of modes.\\[-0.2cm]}\vspace{-0.2cm}
    \label{fig:B2singularvalues}
\end{figure}

To quantify these qualitative observations, we plot the error measure $\epsilon$ for each of the three fields in \cref{fig:B2singularvalues}.
The slow decay of the $\epsilon$-values in this graph confirm that the DEIM representation of the geometry for this benchmark is less effective than it was for the first benchmark problem. For example, the graph indicates that $\scriptstyle{\sim}$70 modes are required to reach an $\epsilon$-values of $\scriptstyle{\sim}$1\% for the $\B{t}$-field, whereas $\scriptstyle{\sim}$30 modes were required to reach the same threshold for benchmark problem 1 (as shown in \cref{fig:B1phireconstructiond}). The acceptable reconstructions of \cref{fig:B2phireconstructiond,fig:B2treconstructiond} again correspond to $\epsilon_\xi \approx 0.01$ and $\epsilon_t \approx 0.03$, implying that these are appropriate target values.

\subsubsection{Multiple holes}
\label{ssec:B3}

\begin{figure}[!b]
    \centering
    \subfloat[Schematic drawing.]{\includegraphics[trim=0 0 0 0,clip, width=0.51\linewidth]{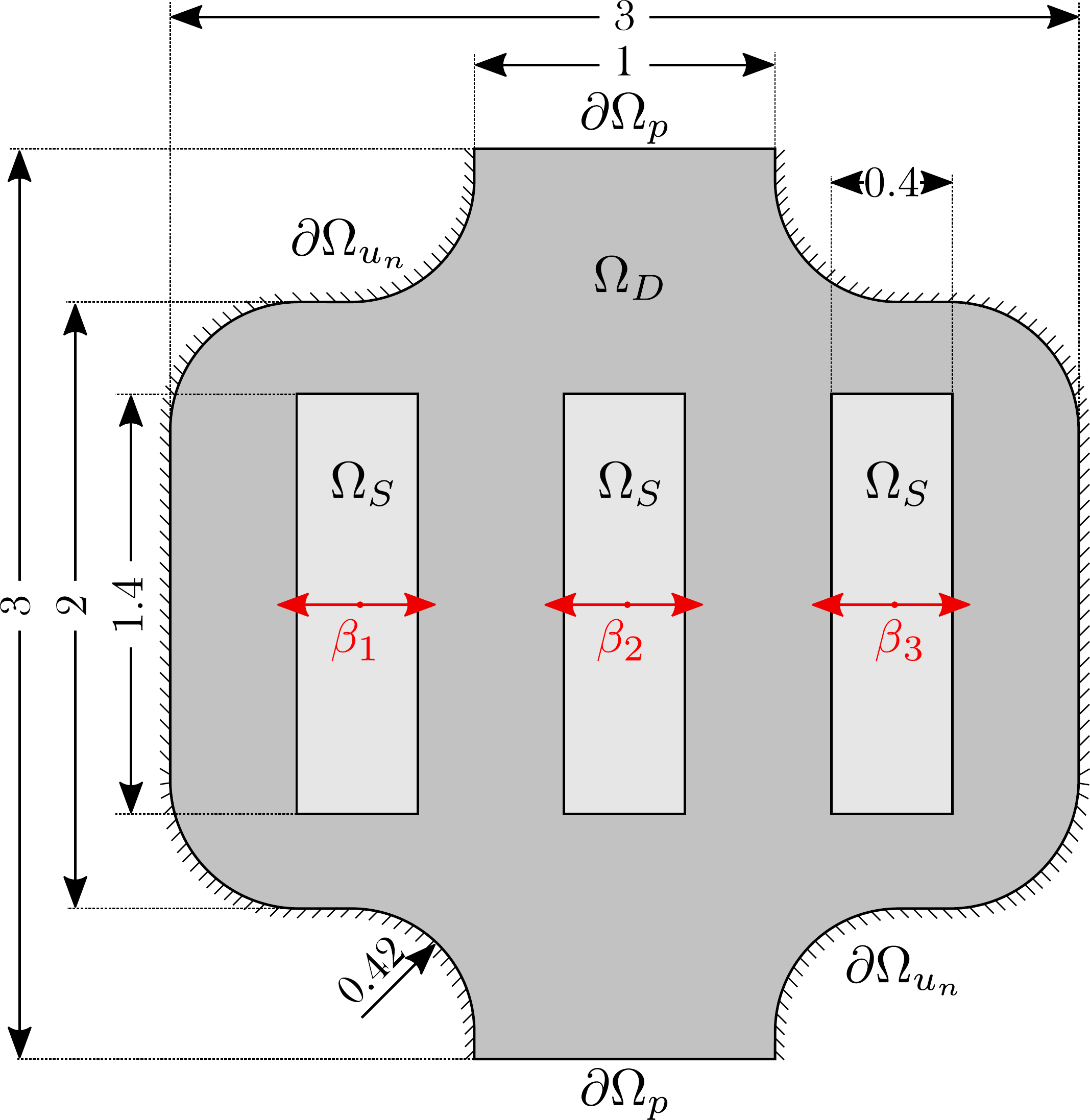}}\hspace{0.1cm}
    \subfloat[High fidelity mesh.]{\includegraphics[trim=0 -30 0 0,clip,width=0.44\linewidth]{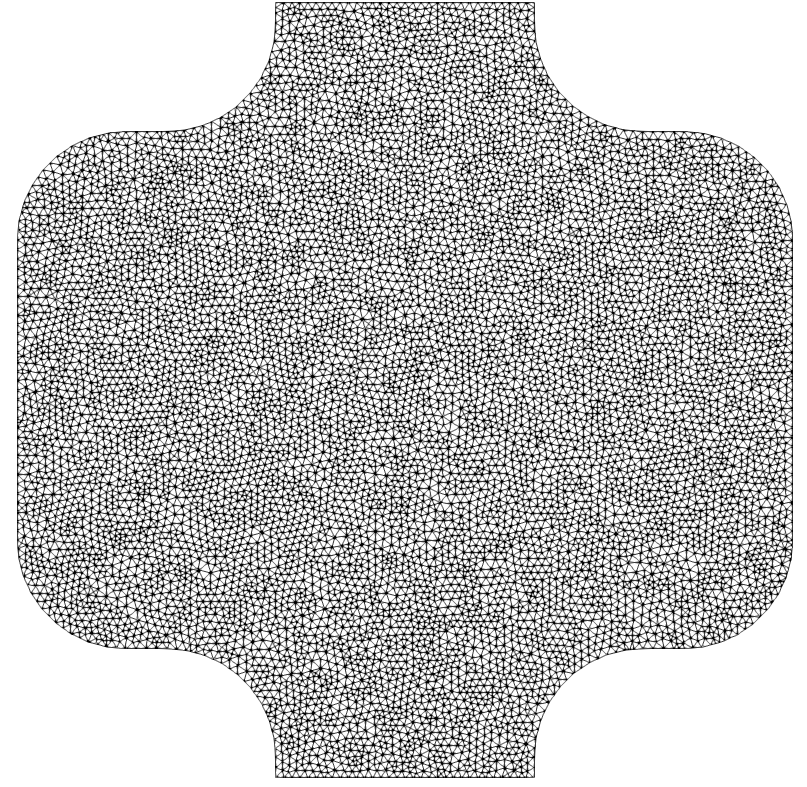}\label{fig:B3Overviewb}} 
    \caption{Benchmark 2: a filtering device.}
    \label{fig:B3Overview}
\end{figure}

As a third benchmark, we consider a block of porous material with three rectangular holes. For the materials parameters we use $\kappa=5\cdot 10^{-4} $m$^2$, $\nu=0.5$Pa$\cdot$s and $\alpha=0.01$Pa$\cdot$s$\cdot$m$^{-2}$, and we set the boundary pressures to 1000Pa and 0Pa at bottom and top boundaries respectively. The location of the three Stokes-flow governed holes is variable. Each hole is permitted to displace by a distance of a single width compared to the base configuration as illustrated in \cref{fig:B3Overview}. The three holes can thus partially overlap to form two, or even just one single hole. 
This benchmark illustrates the capability of describing internal geometry variations with changing topology. Using our phase-field representation of the geometry, we can naturally change the topological properties of the Stokes and Darcy domains.
As emphasized in the introduction, this is not the case for the standard approach to incorporating geometric flexibility in a reduced basis method, which is based on a reinterpretation of the geometric variability as a mapping from a reference domain \cite{quarteroni2015,Hesthaven2016}. 

\begin{figure}[!b]
    \centering
    \subfloat[First mode, $\tilde{\xi}_1(\B{x})$.\mbox{\hspace{-1cm}}]{\includegraphics[trim=70 30 60 30, clip, width=0.515\linewidth]{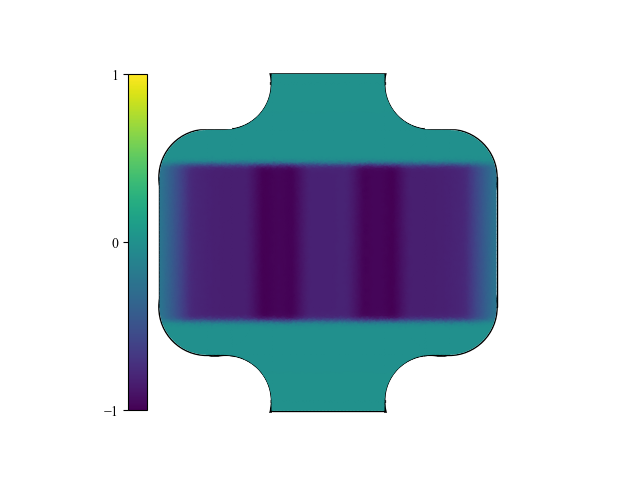} }
    \subfloat[Second mode, $\tilde{\xi}_2(\B{x})$.]{\includegraphics[trim=108 30 90 30, clip,  width=0.41\linewidth]{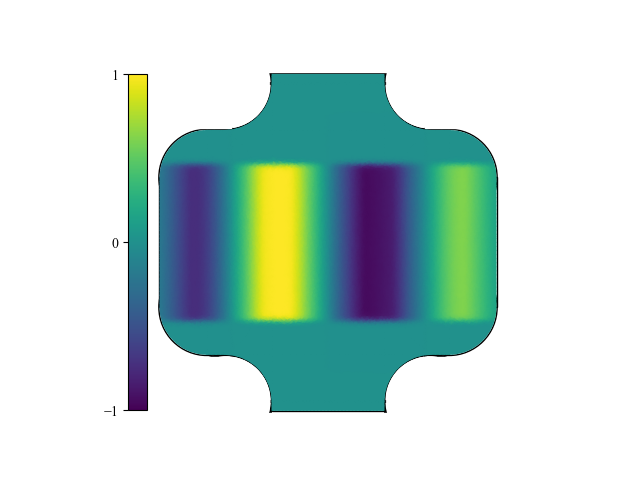} }\\
    \hspace{1cm}\subfloat[Third mode, $\tilde{\xi}_3(\B{x})$.]{\includegraphics[trim=108 30 90 30, clip, width=0.41\linewidth]{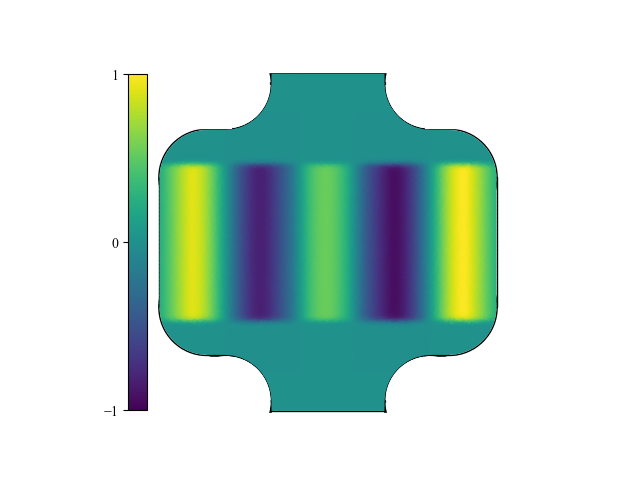}  }\hspace{0.65cm}
    \subfloat[Fourth mode, $\tilde{\xi}_4(\B{x})$.]{\includegraphics[trim=108 30 90 30, clip, width=0.41\linewidth]{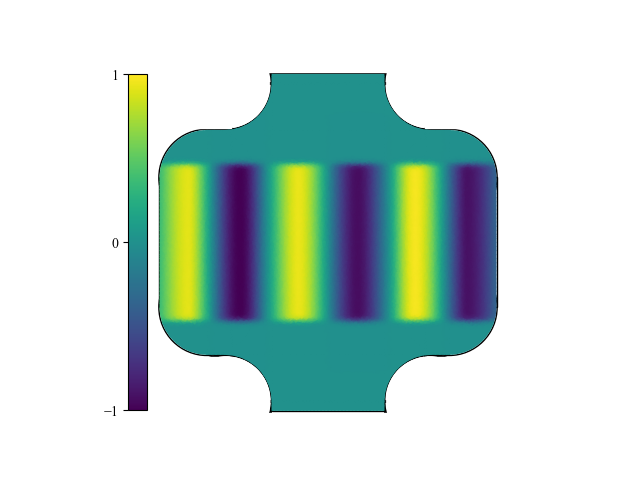} }
    \caption{Benchmark 3: first four DEIM modes of the field $\xi(\beta)$.}
    \label{fig:B3phimodes}
\end{figure}

\begin{figure}[!t]
    \centering
    \subfloat[True field.\mbox{\hspace{-1cm}}]{\includegraphics[trim=60 40 60 50, clip, width=0.51\linewidth]{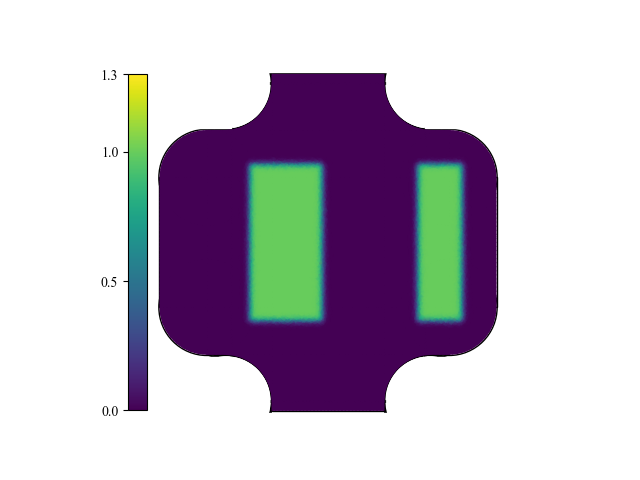}}
    \subfloat[Reconstruction with ten modes.]{\includegraphics[trim=108 40 90 50, clip, width=0.39\linewidth]{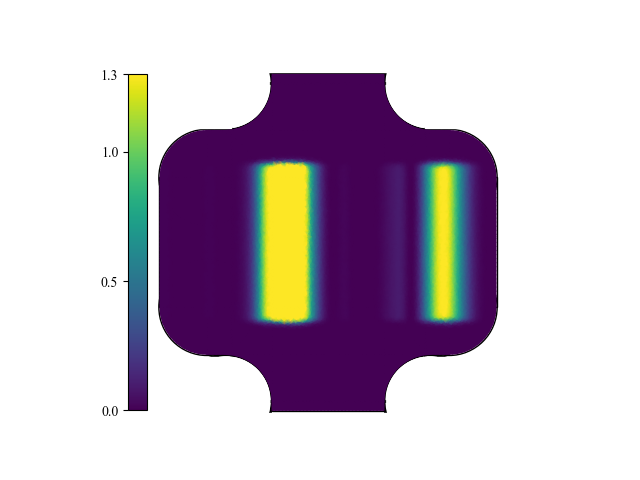} }\\
    \hspace{1.25cm}\subfloat[Reconstruction with twenty modes.]{\includegraphics[trim=108 40 90 40, clip, width=0.39\linewidth]{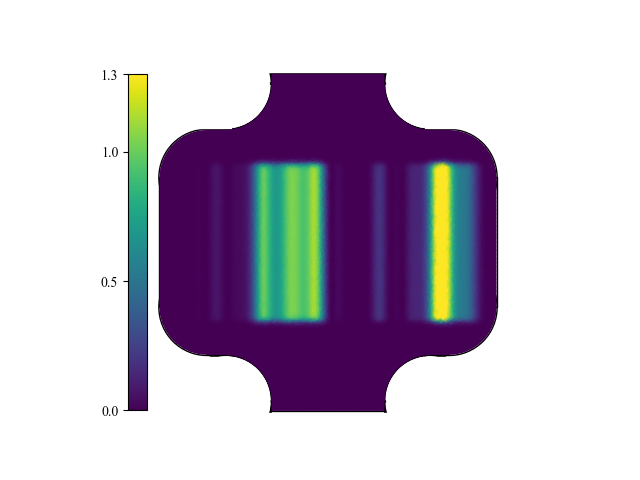}  }\hspace{0.5cm}
    \subfloat[Reconstruction with forty modes.]{\includegraphics[trim=108 40 90 40, clip, width=0.39\linewidth]{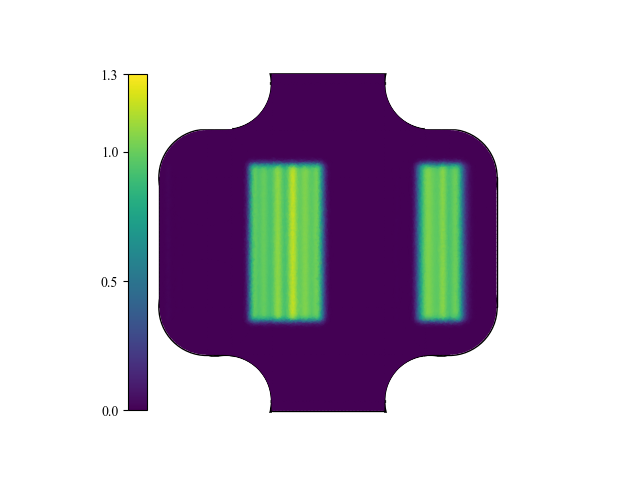}\label{fig:B3phireconstructiond} }
    \caption{Benchmark 3: non-negativity preserving DEIM reconstruction of $\phi([0.35,-0.25,0.1])$.}
    \label{fig:B3phireconstruction}
\end{figure}

\begin{figure}[!b]
\vspace{-0.75cm}
    \centering
    \subfloat[Reconstruction with ten modes.]{\includegraphics[trim=60 40 60 40, clip, width=0.51\linewidth]{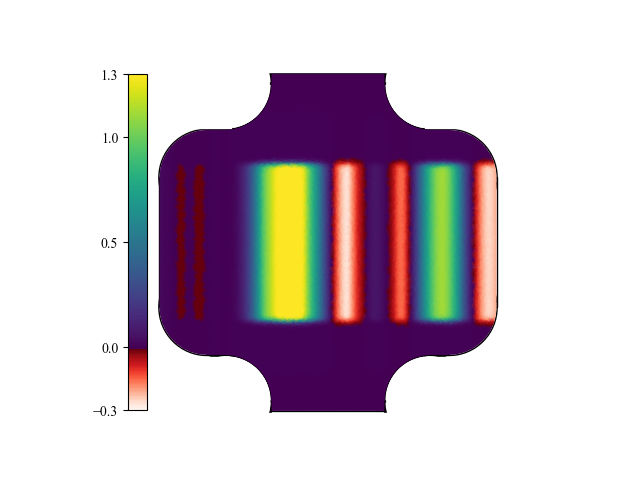}  }\hspace{0.2cm}
    \subfloat[Reconstruction with twenty modes.]{\includegraphics[trim=108 40 90 40, clip, width=0.39\linewidth]{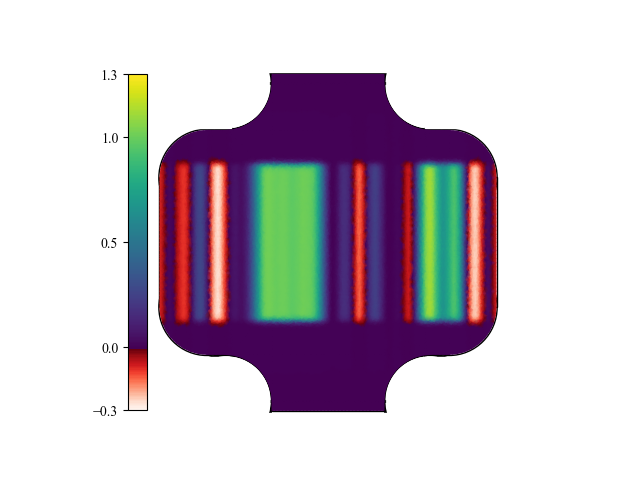} }
    \caption{Benchmark 3: classical DEIM reconstruction of $\phi([0.35,-0.25,0.1])$.\\[-0.4cm]}
    \label{fig:B3phireconstructionNegative}
    \vspace{-0.0cm}
\end{figure}

\begin{figure}[!b]
    \centering
    \subfloat[First mode, $\tilde{\B{t}}_1(\B{x})$.\mbox{\hspace{-1cm}}]{\includegraphics[trim=90 50 120 70, clip, width=0.55\linewidth]{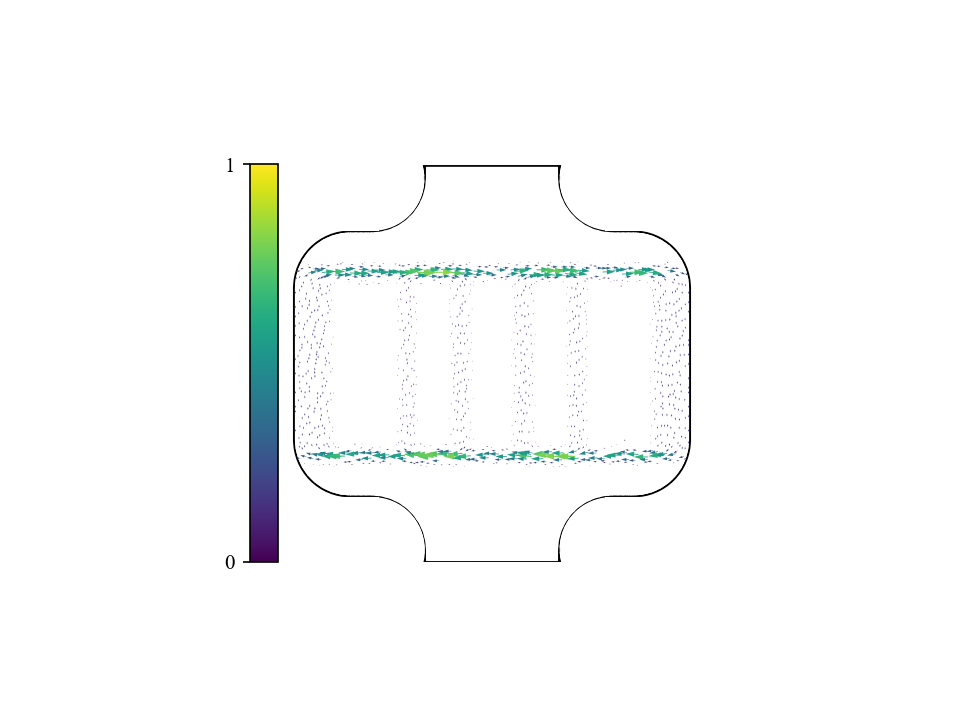} }
    \subfloat[Second mode, $\tilde{\B{t}}_2(\B{x})$.]{\includegraphics[trim=140 50 120 70, clip, width=0.44\linewidth]{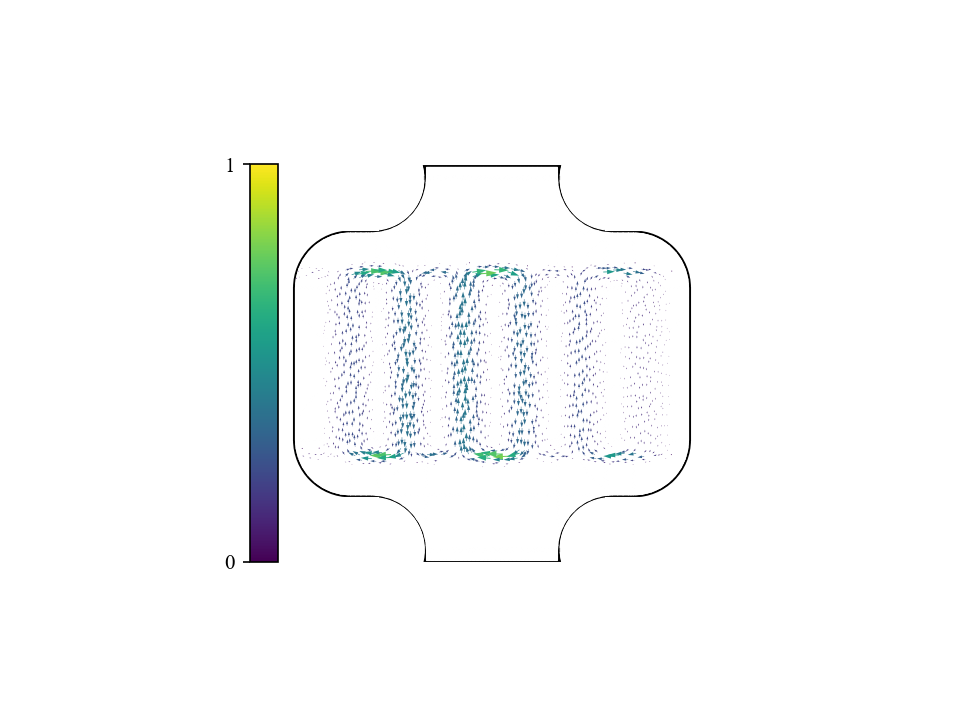} }\\\hspace{1.5cm}
    \subfloat[Third mode, $\tilde{\B{t}}_3(\B{x})$.]{\includegraphics[trim=140 50 120 70, clip, width=0.44\linewidth]{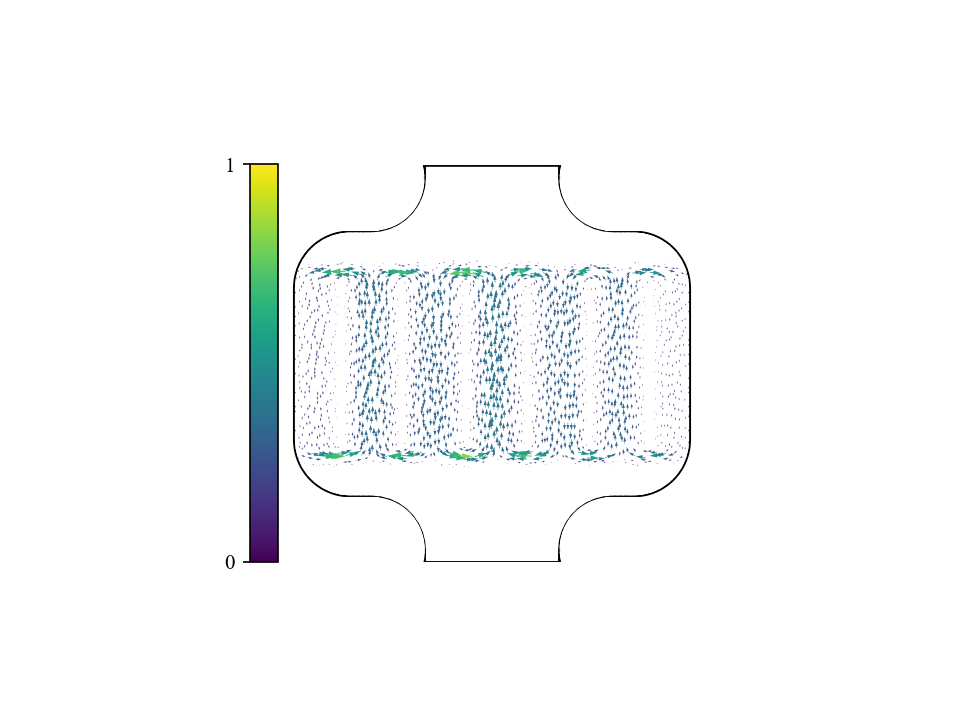}  }
    \subfloat[Fourth mode, $\tilde{\B{t}}_4(\B{x})$.]{\includegraphics[trim=140 50 120 70, clip, width=0.44\linewidth]{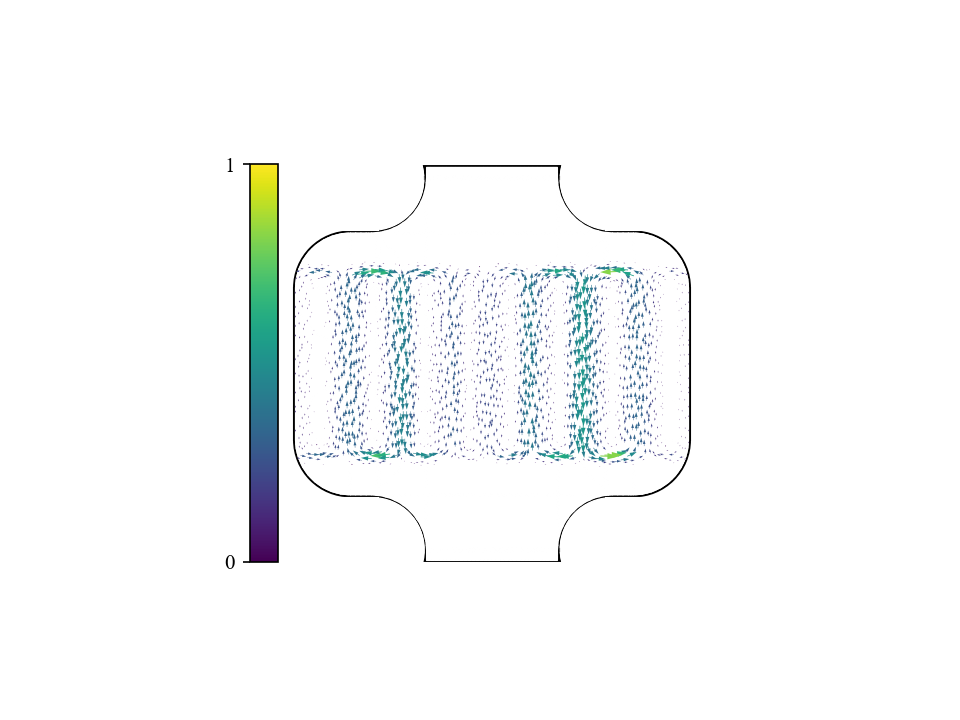} }
    \caption{Benchmark 3: first four DEIM modes of the field ${\B{t}}(\beta)$.}
    \label{fig:B3tmodes}
\end{figure}

\begin{figure}[!t]
    \centering
    \subfloat[True field.\mbox{\hspace{-1cm}}]{\includegraphics[trim=80 70 120 70, clip, width=0.51\linewidth]{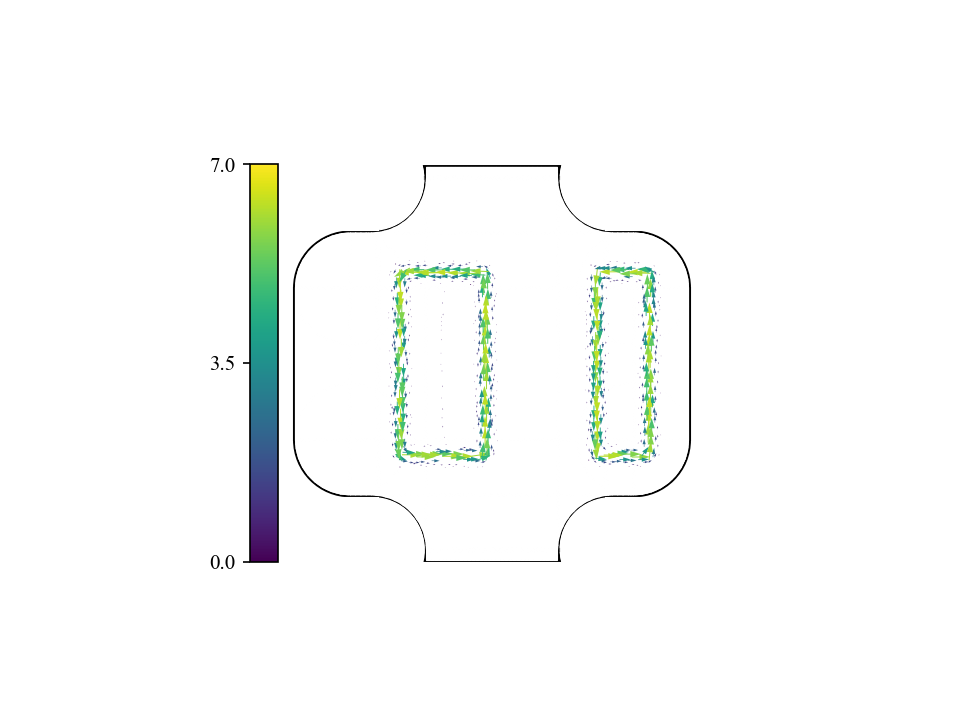}}
    \subfloat[Reconstruction with ten modes.]{\includegraphics[trim=140 70 120 80, clip, width=0.39\linewidth]{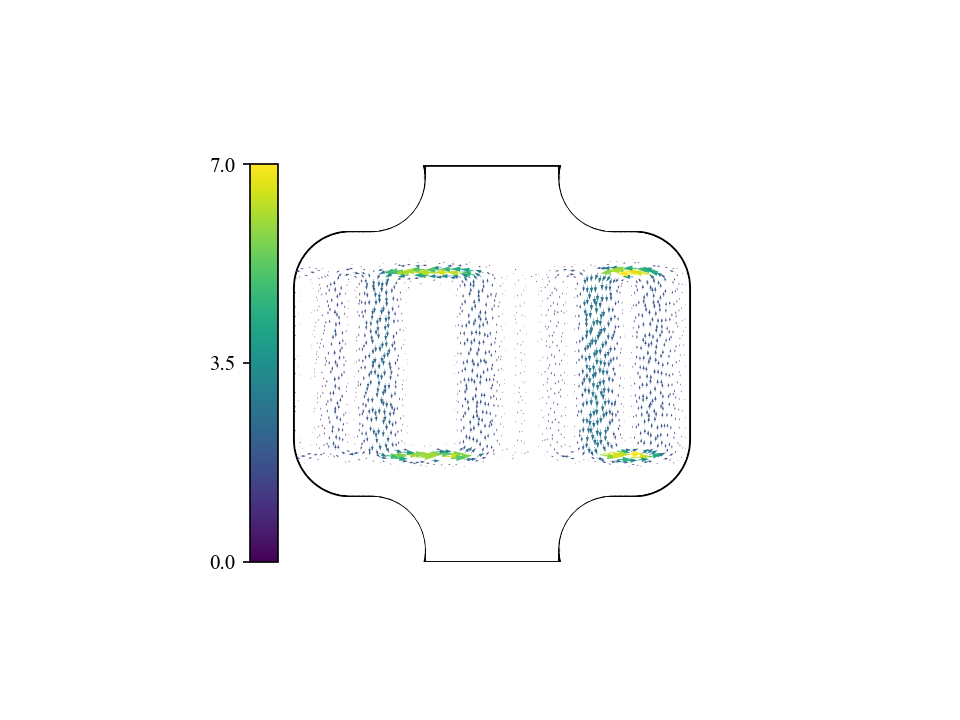} }\\\hspace{1.5cm}
    \subfloat[Reconstruction with twenty modes.]{\includegraphics[trim=140 70 120 80, clip, width=0.39\linewidth]{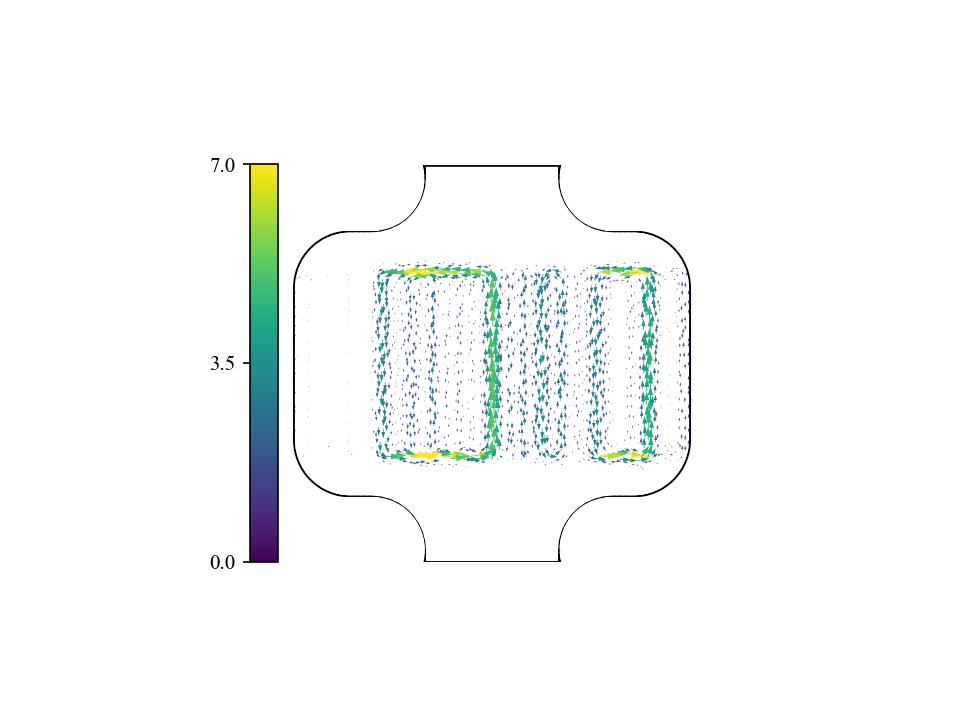}  }
    \subfloat[Reconstruction with fifty modes.]{\includegraphics[trim=140 70 120 80, clip, width=0.39\linewidth]{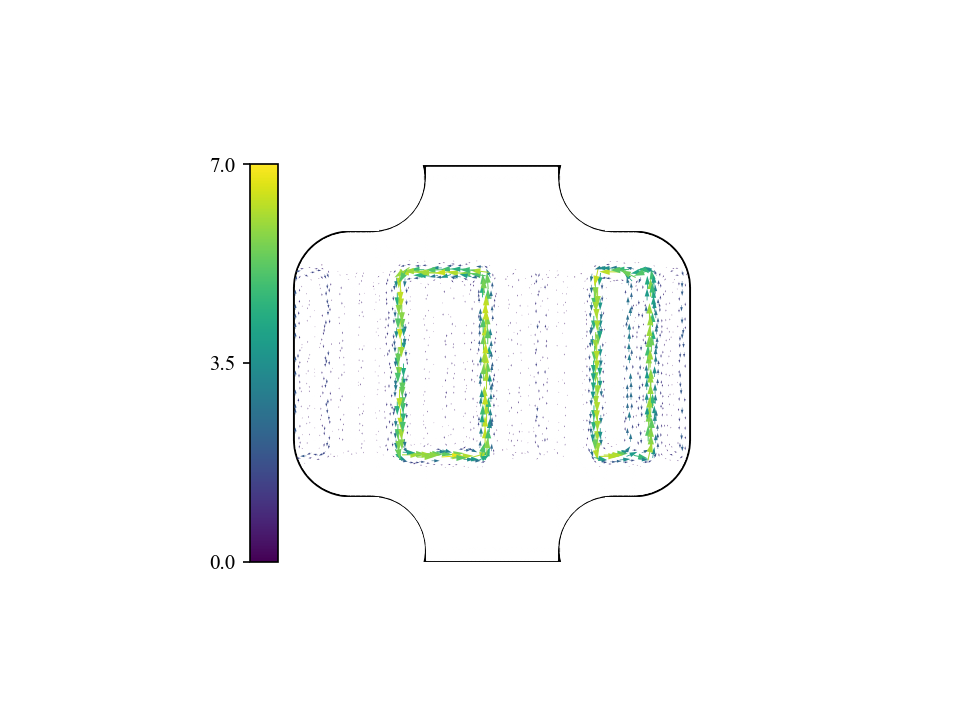} }
    \caption{Benchmark 3: DEIM reconstruction of $\B{t}([0.35,-0.25,0.1])$.}
    \label{fig:B3treconstruction}
\end{figure}

 The parameter space $\mathbb{P}$ of this benchmark problem is three-dimensional. We construct the sampling set $\mathbb{P}^h$ as a Cartesian product of 21 uniformly distributed samples along each axis, yielding a set of 9261 samples: $\mathbb{P}_h = \{ -0.5 , -0.45, \cdots, 0.5 \} \times \{ -0.5 , -0.45, \cdots, 0.5 \}\times \{ -0.5 , -0.45, \cdots, 0.5 \}$. With this sampling set, we obtain the DEIM modes for the fields $\xi$ and $\B{t}$ depicted in \cref{fig:B3phimodes,fig:B3tmodes}. \Cref{fig:B3phireconstruction,fig:B3phireconstructionNegative} then show reconstructions of the field $\phi([0.35,-0.25,0.1])$ with the non-negativity preserving DEIM and the standard DEIM approaches, and \cref{fig:B3treconstruction} shows the reconstruction of $\B{t}([0.35,-0.25,0.1])$. Note that this parameter point produces an overlap of the left two rectangular Stokes domains. 


\begin{figure}[!b]
\vspace{-0.5cm}
    \centering
    \includegraphics[trim=0 0 0 0, clip,width=0.6\linewidth]{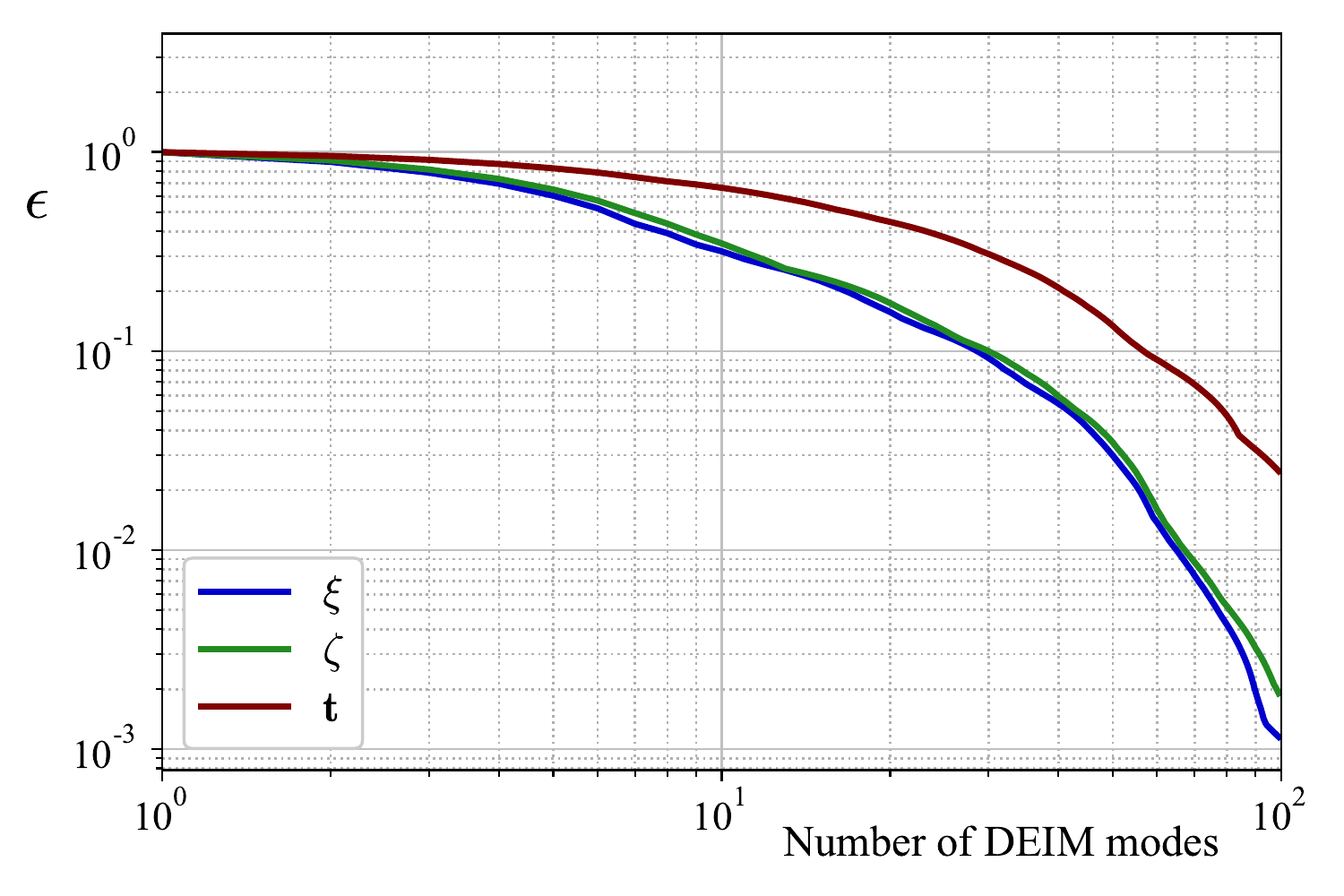}
    \caption{Benchmark 3: convergence of $\epsilon$ with increasing number of modes for $\xi$, $\zeta$ and $\B{t}$.\\[-0.2cm]}
    \label{fig:B3singularvalues}
    \vspace{-0.1cm}
\end{figure}

The approximations of $\phi$ and $\B{t}$ still show defects even when 40 interpolation modes are used. 
Naturally, this is due to the higher dimensionality of the parameter space. We can determine roughly how many modes would be required to obtain adequate approximations by investigating the error measure $\epsilon(n)$. This measure is plotted in \cref{fig:B3singularvalues} for each of the three fields. When compared to \cref{fig:B1singularvalues,fig:B2singularvalues}, we indeed observe a significantly larger overall $\epsilon$. Where benchmarks 1 and 2 require 30 and 70 modes respectively to reach an $\epsilon_t$ of 1\%, this threshold is not reached for benchmark 3 within the first 100 modes. For benchmarks 1 and 2 $\epsilon$-values of 1\% and 3\% are sufficient to obtain an acceptable approximation of $\phi$ and $\B{t}$ respectively. For the current benchmark problem, these $\epsilon$-values are reached for $\scriptstyle{\sim}$65 and $\scriptstyle{\sim}$80 modes respectively. We consider this an excessive amount. In the next section, we investigate the importance of the accuracy of the DEIM reconstructions of $\phi(\beta)$, $\psi(\beta)$ and $\B{t}(\beta)$ on the fidelity of the reduced basis formulation.

\section{Reduced basis method of the coupled equations on the parametric domain}
\label{sec:RB}

\begin{figure}[!b]
    \centering
    \subfloat[First mode.\mbox{\hspace{-1cm}}]{\includegraphics[trim=50 45 100 50, clip, width=0.52\linewidth]{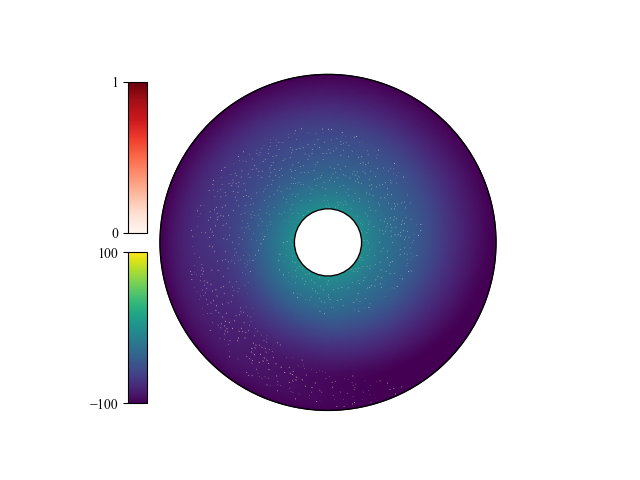} }\hspace{.0cm}
    \subfloat[Second mode.]{\includegraphics[trim=110 45 93 50, clip,  width=0.43\linewidth]{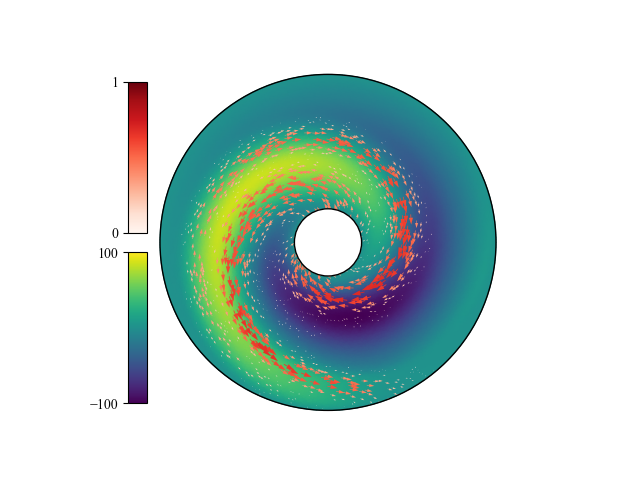} }\\\hspace{1.5cm}
    \subfloat[Third mode.]{\includegraphics[trim=110 45 93 50, clip, width=0.43\linewidth]{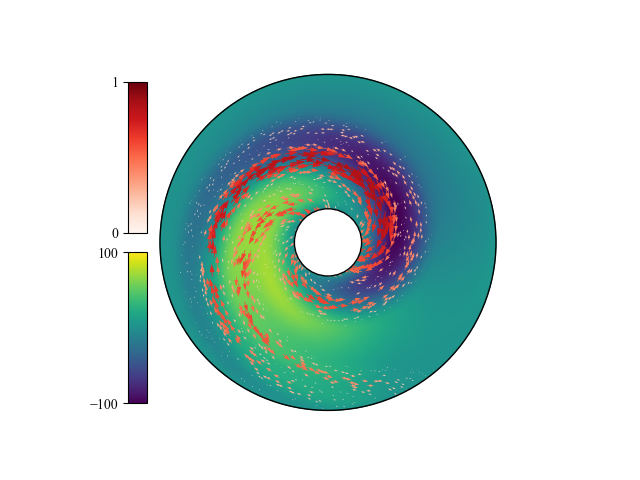}  }
    \subfloat[Fourth mode.]{\includegraphics[trim=110 45 93 50, clip, width=0.43\linewidth]{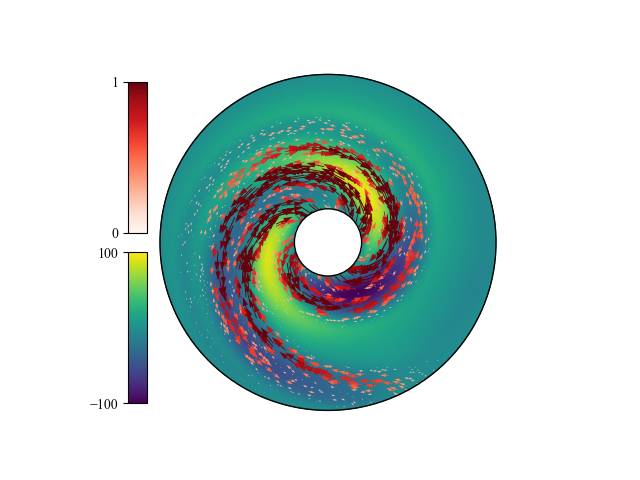} }
    \caption{First four reduced basis modes of benchmark 1.\\[-0.7cm]}
    \label{fig:B1solmodes}
    \vspace{-0.3cm}
\end{figure}

From \cref{sec:Diffuse}, we have the following parameter dependent finite element formulation:
\begin{align}\label{HFFE}
& \text{Find }\B{u}^h(\beta),p^h(\beta)\in \B{\mathcal{V}}^h\times \QQ^h \text{ s.t. } \forall\, \B{v},q\in \B{\mathcal{V}}^h\times \QQ^h: \nonumber\\
&  \quad  B((\B{u}^h,p^h),(\B{v},q);\beta) = L((\B{v},q))\,,
\end{align}
where the bilinear and linear forms are those of \cref{WeakDif}. The mixed finite element space $\B{\mathcal{V}}^h\times \QQ^h $ is that described in \cref{ssec:FEspaces} and could be kept independent of $\beta$. This approximation space is assumed to be sufficiently refined such that the finite element solutions $(\B{u}^h(\beta),p^h(\beta))$ can for all intents and purposes be considered the true solution $(\B{u}(\beta),p(\beta))$. The meshes illustrated in \cref{fig:B1Overviewb,fig:B2Overviewb,fig:B3Overviewb} satisfy this requirement for the previous three benchmark problems. 

We now wish to find a low-dimensional subspace $RB^r\subset \B{\mathcal{V}}_0^h\times \QQ^h$ (with $r = \dim( RB^r )$) that has high approximation power with respect to the full solution manifold $\mathcal{M}_{(u^h,p^h)} = \{ (\B{u}^h(\beta),p^h(\beta)) : \beta\in\mathbb{P}\}$. Once this low-dimensional subspace has been found, it may be used in place of $\B{\mathcal{V}}_0^h\times \QQ^h$ in \cref{HFFE} to produce a reduced order model for approximating any solution pair $(\B{u}^h,p^h)\approx (\B{u},p)$ in $\mathcal{M}_{(u^h,p^h)}$ for any $\beta\in\mathbb{P}$ at a very low computational cost.

\begin{figure}[!b]
    \centering
    \subfloat[First mode.\mbox{\hspace{-1cm}}]{\includegraphics[trim=60 30 80 30, clip, width=0.54\linewidth]{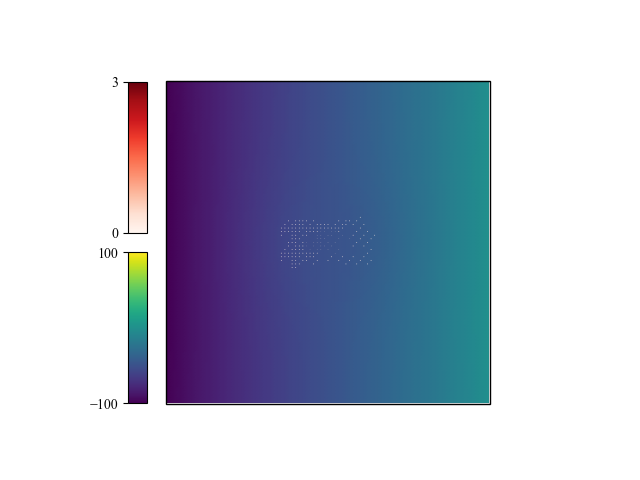} }
    \subfloat[Second mode.]{\includegraphics[trim=120 30 80 30, clip,  width=0.44\linewidth]{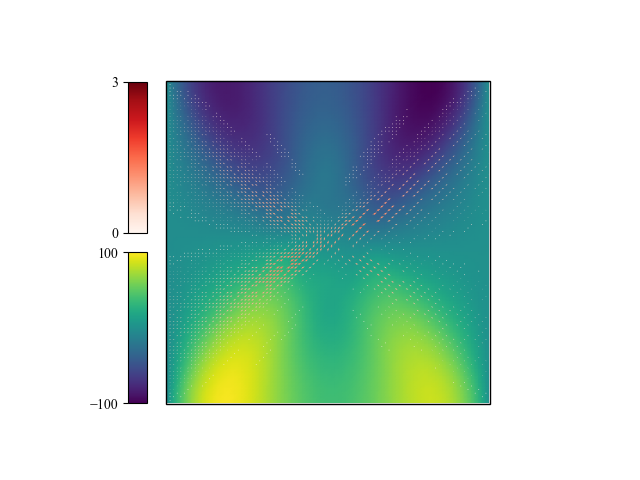} }\\\hspace{1.5cm}
    \subfloat[Third mode.]{\includegraphics[trim=120 30 80 30, clip, width=0.44\linewidth]{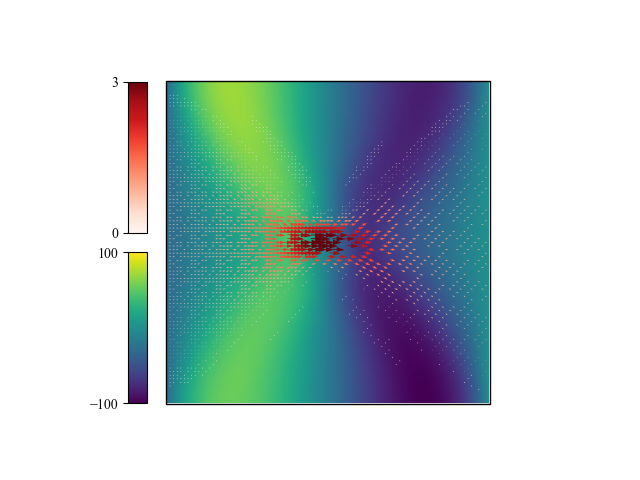}  }
    \subfloat[Fourth mode.]{\includegraphics[trim=120 30 80 30, clip, width=0.44\linewidth]{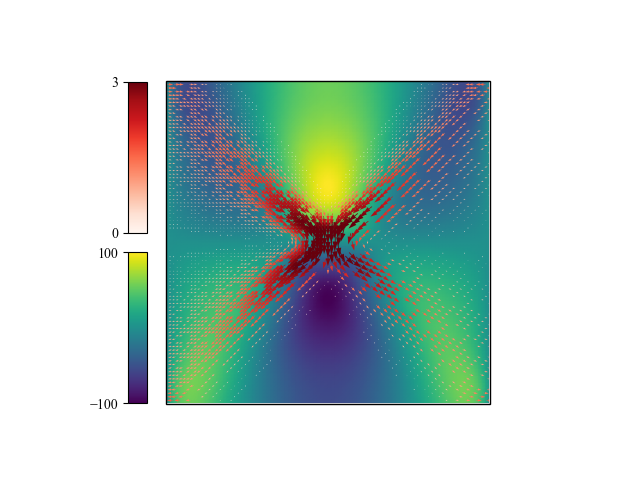} }
    \caption{First four reduced basis modes of benchmark 2.}
    \label{fig:B2solmodes}
\end{figure}

\begin{figure}[!t]
    \centering
    \subfloat[First mode.\mbox{\hspace{-1cm}}]{\includegraphics[trim=90 70 120 60, clip, width=0.54\linewidth]{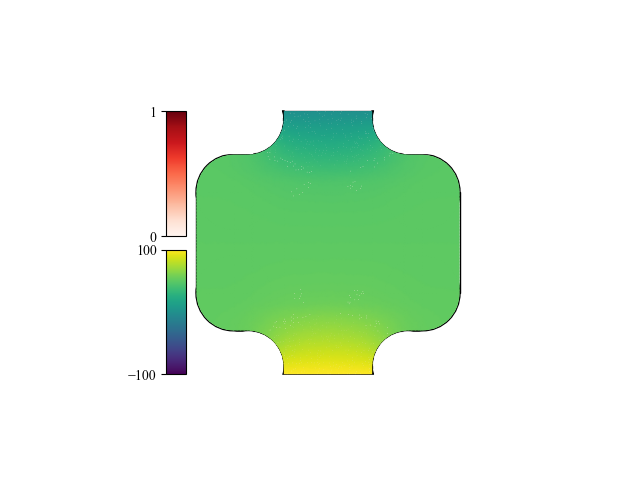} }
    \subfloat[Second mode.]{\includegraphics[trim=140 70 120 60, clip, width=0.43\linewidth]{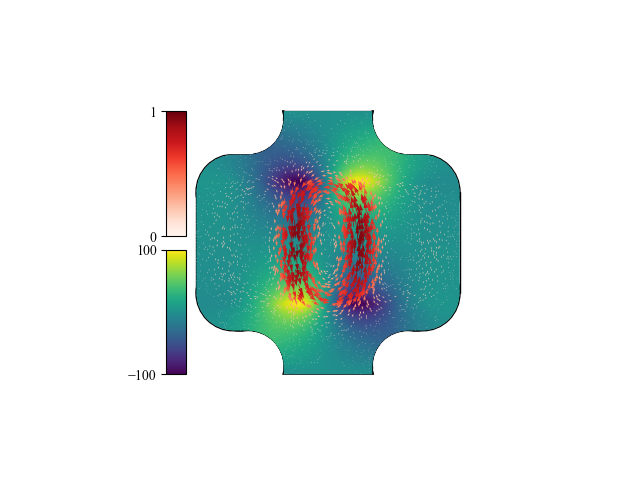} }\\\hspace{1.5cm}
    \subfloat[Third mode.]{\includegraphics[trim=140 70 120 60, clip, width=0.44\linewidth]{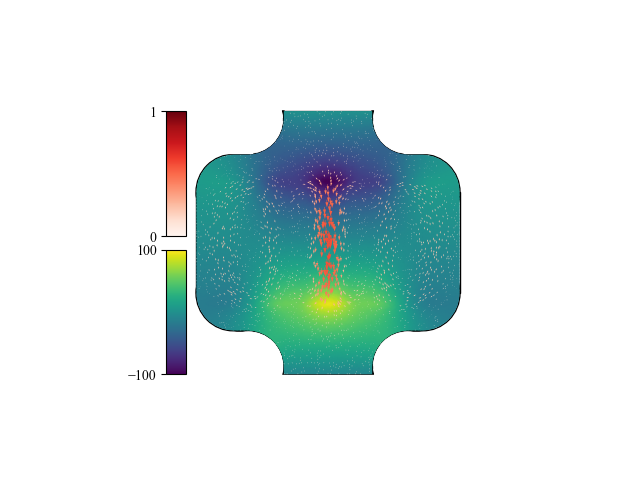}  }
    \subfloat[Fourth mode.]{\includegraphics[trim=140 70 120 60, clip, width=0.44\linewidth]{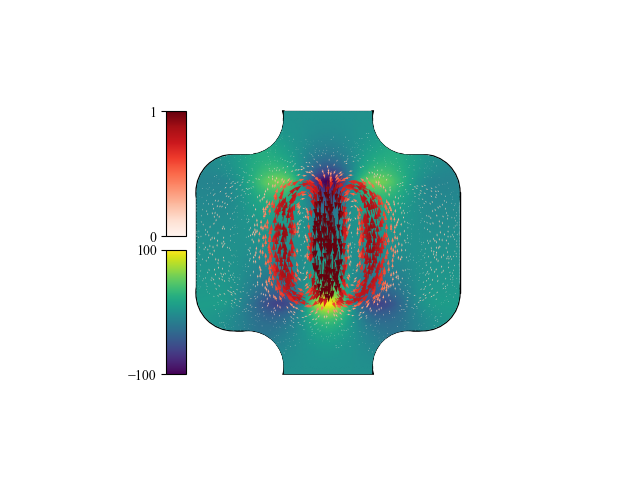} }
    \caption{First four reduced basis modes of benchmark 3.}
    \label{fig:B3solmodes}
\end{figure}

The procedure for determining $RB^r$ again relies on a singular value decomposition of a snapshot matrix, which is in this case constructed from solution vectors of high-fidelity finite element solutions. We perform this computation for all three benchmark problems. The snapshot matrices are filled with solutions to \cref{HFFE} on the meshes of \cref{fig:B1Overviewb,fig:B2Overviewb,fig:B3Overviewb} for the same parameter samplings described in \cref{ssec:B1,ssec:B2,ssec:B3}. \Cref{fig:B1solmodes,fig:B2solmodes,fig:B3solmodes} show the first four resulting reduced basis functions for each of the three benchmark problems. The figures show pressure fields overlapped by velocity fields. This is indicative of the nature of these functions: since we perform a singular value decomposition of a snapshot matrix of velocity-pressure pairs, each of the left singular vectors represents a combination of a velocity and a pressure field.

\begin{figure}[!t]
\vspace{0.5cm}
    \centering
    \includegraphics[trim=0 0 0 0, clip,width=0.6\linewidth]{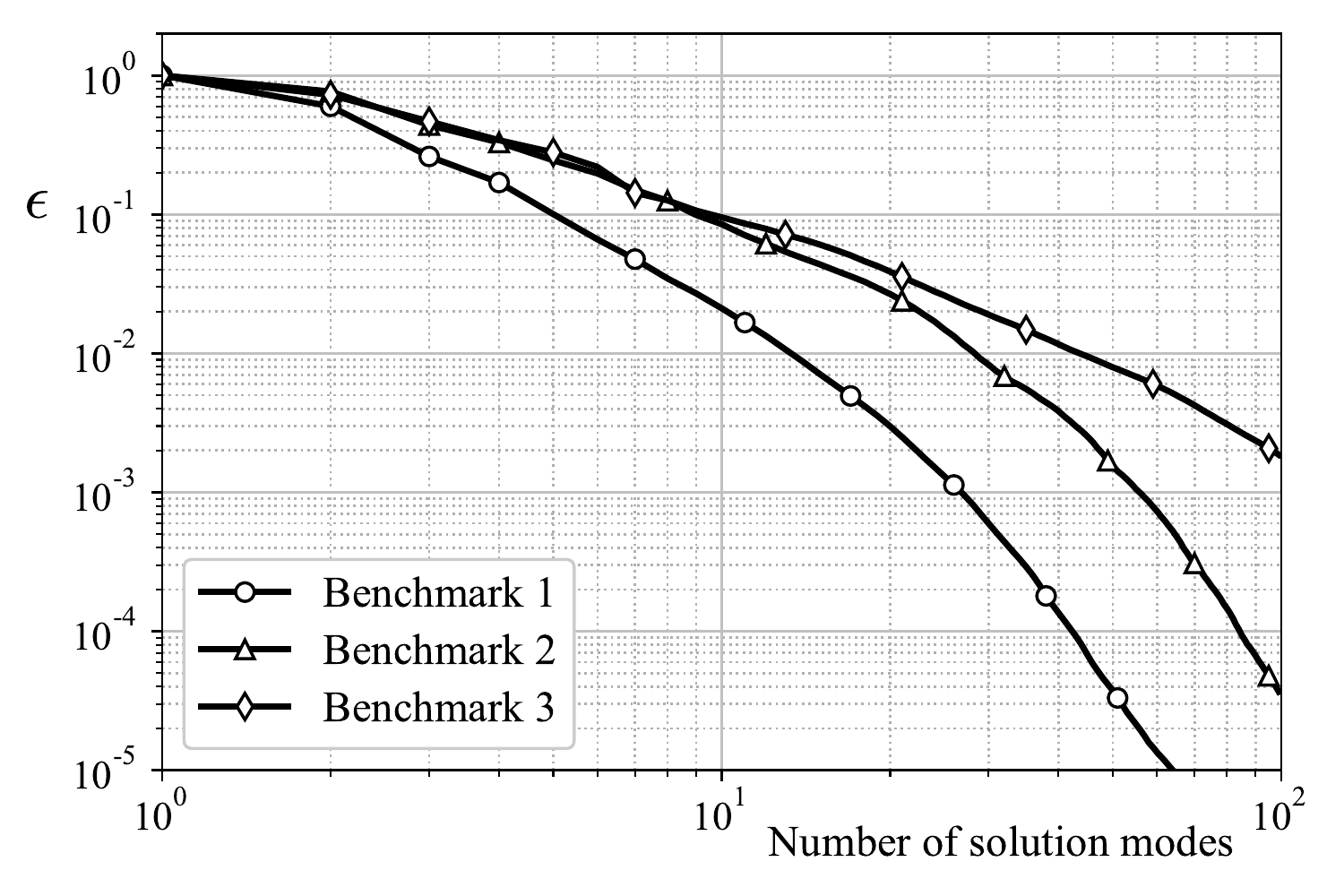}
    \caption{Convergence of $\epsilon$ for the solution fields with the number of modes.}
    \label{fig:Uepsilon}
\end{figure}

The singular values corresponding to the reduced basis functions again relate to the approximation power of the cumulative basis compared to the full (discrete) solution manifold. This is quantified by the measure $\epsilon$ from \cref{errormeasure}, which is plotted for all three benchmark problems in \cref{fig:Uepsilon}. We observe the same trends as we did for the DEIM geometry representations in \cref{sec:DEIM}: benchmark 1 performs best, followed by benchmark 2, and benchmark 3 is worst. Again, this can be motivated from the complexity of the problems, and the related Kolmogorov $n$-width of the solution manifold. The higher irregularity of benchmark 2 produces a higher Kolmogorov $n$-width than the better-behaved benchmark problem 1, but the higher dimensionality of benchmark 3 has a more significant impact. Nevertheless, all three benchmark problems achieve an $\epsilon$-value below 5\% for less than 20 solution modes.

In order to efficiently use these functions as basis functions in the reduced basis formulation of \cref{HFFE}, we require an affine decomposition of the (bi)linear forms into parameter independent bilinear forms multiplied by parameter dependent weights. This permits a precomputation of the corresponding stiffness matrices during an offline phase, such that a full integration of the reduced basis functions can be avoided during the online phase. Based on the DEIM representation of all parameter dependent fields from \cref{sec:DEIM}, an affine decomposition of the bilinear form follows as:
\begin{align}
\begin{split}
&    B((\B{u}^h,p^h),(\B{v},q);\beta) \approx B_0((\B{u}^h,p^h),(\B{v},q)) + \sum\limits_{i=1}^{N_\phi} \theta^\phi_i(\beta) B^{\phi}_i((\B{u}^h,p^h),(\B{v},q)) \\[-0.15cm]
&   \hspace{1cm} + \sum\limits_{i=1}^{N_\psi} \theta^\psi_i(\beta) B^{\psi}_i((\B{u}^h,p^h),(\B{v},q))  + \sum\limits_{i=1}^{N_\alpha} \theta^\alpha_i(\beta) B^{\alpha}_i((\B{u}^h,p^h),(\B{v},q)) \,,
    \end{split}
\end{align}
with:
\begin{subequations}
\begin{alignat}{3}
& B_0((\B{u}^h,p^h),(\B{v},q)) &&= - \int\limits_{\Omega}  p^h \, \nabla \cdot \B{v} \dO - \int\limits_{\Omega} q\, \nabla \cdot \B{u}^h \dO \,,\\
& B^{\phi}_i((\B{u}^h,p^h),(\B{v},q)) &&= \!\int\limits_{\Omega} \tilde{\phi}_i \, \mu \nabla^s \B{u}^h : \nabla^s \B{v} \dO \quad\quad\text{for }1 \leq i \leq N_\phi \,,\label{mat_xi}\\
& B^{\psi}_i((\B{u}^h,p^h),(\B{v},q)) &&= \!\int\limits_{\Omega} \tilde{\psi}_i \, \mu \B{\kappa}^{-1} \B{u}^h \cdot \B{v} \dO \hspace{1.2cm}\text{for }1 \leq i \leq N_\psi  \,,\label{mat_zeta}\\
& B^{\alpha}_i((\B{u}^h,p^h),(\B{v},q)) &&= \! \int\limits_{\Omega} \tilde{\B{\alpha}}_i\B{u}^h\cdot \B{v}\dO \hspace{1.95cm}\!\!\text{for }1 \leq i \leq N_\alpha  \,.\label{mat_t}\end{alignat}
\end{subequations}
With the functions from \cref{fig:B1solmodes,fig:B2solmodes,fig:B3solmodes} as approximation bases for the solution space, our final reduced order model reads:
\begin{align}\label{ROM}
& \text{Find }\underbar{u}^{RB}\in RB^r \text{ s.t. } \forall\, \underbar{v}\in RB^r: \nonumber\\
\begin{split}
&  \quad  B_0(\underbar{u}^{RB},\underbar{v}) + \sum\limits_{i=1}^{N_\phi} \theta^\phi_i(\beta) B^{\phi}_i(\underbar{u}^{RB},\underbar{v})  + \sum\limits_{i=1}^{N_\psi} \theta^\psi_i(\beta) B^{\psi}_i(\underbar{u}^{RB},\underbar{v}) + \sum\limits_{i=1}^{N_\alpha} \theta^\alpha_i(\beta) B^{\alpha}_i(\underbar{u}^{RB},\underbar{v}) = L(\underbar{v})\,, 
    \end{split}
\end{align}
where $\underbar{u}$ represents the solution tuple that combines the velocity and the pressure solution. The matrix representation of this problem becomes:
\begin{align}\label{ROM_mat}
& \text{Find }\hat{\underbar{u}}^{RB}\in \mathbb{R}^r \text{ s.t.} : \nonumber\\
\begin{split}
&  \quad \Big( \underbar{\underbar{K}}_0 + \sum\limits_{i=1}^{N_\phi}  \theta^\phi_i(\beta) \underbar{\underbar{K}}_i^\phi  + \sum\limits_{i=1}^{N_\psi}  \theta^\psi_i(\beta)\underbar{\underbar{K}}_i^\psi + \sum\limits_{i=1}^{N_\alpha} \theta^\alpha_i(\beta) \underbar{\underbar{K}}_i^\alpha \Big)\, \hat{\underbar{u}}^{RB}  = \underbar{F}\,, 
    \end{split}
\end{align}
where all $\underbar{\underbar{K}}\in \mathbb{R}^{r\times r}$, $ \underbar{F}\in \mathbb{R}^{r}$, and $\hat{\underbar{u}}^{RB}$ is the vector of coefficients for each of the reduced basis functions.

\begin{figure}[!b]
    \centering
    \subfloat[True solution.\mbox{\hspace{-1cm}}]{\includegraphics[trim=45 40 93 50, clip, width=0.54\linewidth]{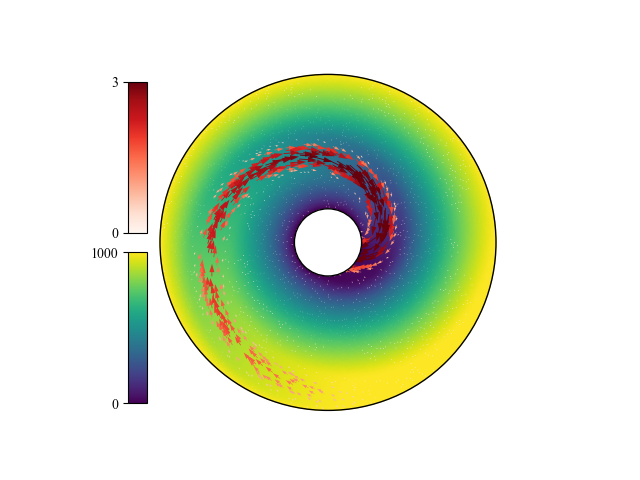} }
    \subfloat[Approximation with five modes.]{\includegraphics[trim=108 40 93 50, clip, width=0.44\linewidth]{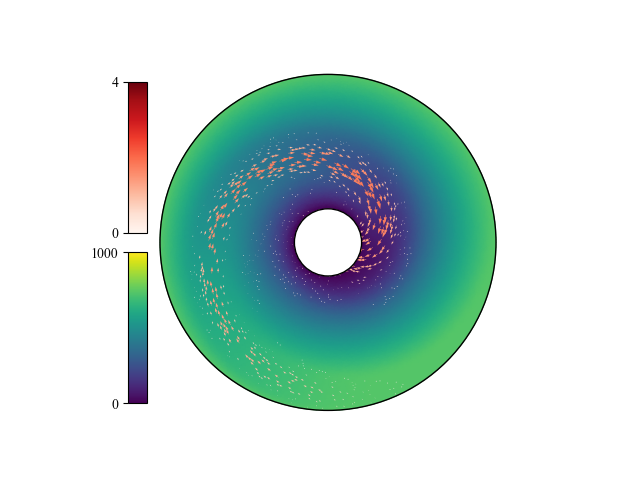} }\\
    \hspace{1.5cm}\subfloat[Approximation with ten modes.]{\includegraphics[trim=108 40 93 50, clip, width=0.44\linewidth]{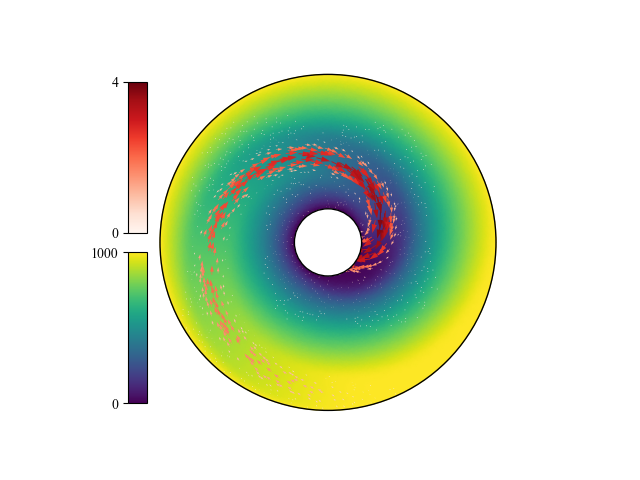}  }
    \subfloat[Approximation with twenty modes.]{\includegraphics[trim=108 40 93 50, clip, width=0.44\linewidth]{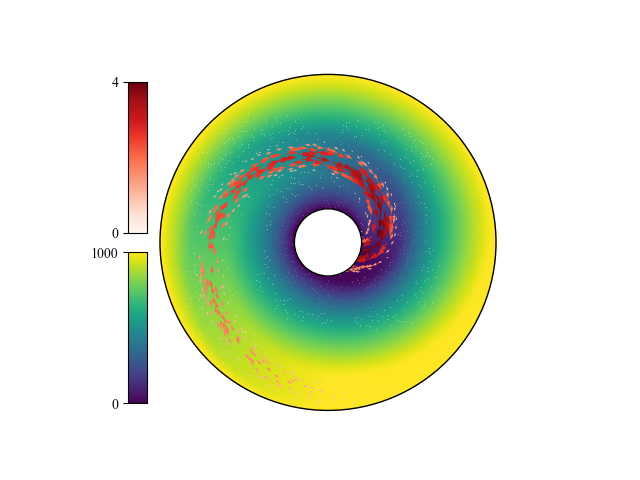} }
    \caption{Reduced order model simulation result for benchmark 1 with $\beta=360^\circ$.\\[-0.7cm]}
    \label{fig:B1RBmeth}
    \vspace{-0.3cm}
\end{figure}

This reduced order model is operated by taking the following steps:
\begin{enumerate}
    \item In the offline phase (i.e., before operation), the matrices corresponding to the bilinear forms $B^{\phi}_i$, $B^{\psi}_i$ and $B^{\alpha}_i$ are precomputed.
    \item Also in the offline phase, the interpolation matrix from \cref{DEIMinterp} corresponding to each of the weight vectors $\underline\theta^\xi(\beta)$, $\underline\theta^\zeta(\beta)$ and $\underline\theta^t(\beta)$ is inverted and stored.
    \item In the online phase (i.e., during operation), the weight vectors  $\underline\theta^\xi(\beta)$, $\underline\theta^\zeta(\beta)$ and $\underline\theta^t(\beta)$ are computed for the given parameter $\beta$.
    \item From those weights, we compute the weights $ \theta^\phi_i(\beta)$, $ \theta^\psi_i(\beta)$ and $ \theta^\alpha_i(\beta)$ for to the non-negativity preserving approximation of \cref{nonnegreconstruction}.
    \item The summation from \cref{ROM} is carried out and the resulting $r\times r$ system of equation is be solved.
\end{enumerate}
Subsequently, the obtained low-order approximation can be visualized in a postprocessing step by computing the weighted sum of the reduced basis functions and the computed coefficients.

\begin{figure}[!b]
    \centering
    \subfloat[True solution.\mbox{\hspace{-1cm}}]{\includegraphics[trim=50 50 100 50, clip, width=0.54\linewidth]{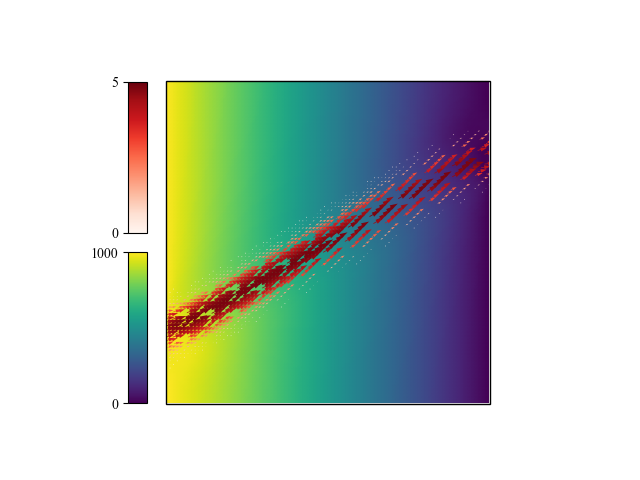} }
    \subfloat[Approximation with five modes.]{\includegraphics[trim=108 50 100 50, clip, width=0.44\linewidth]{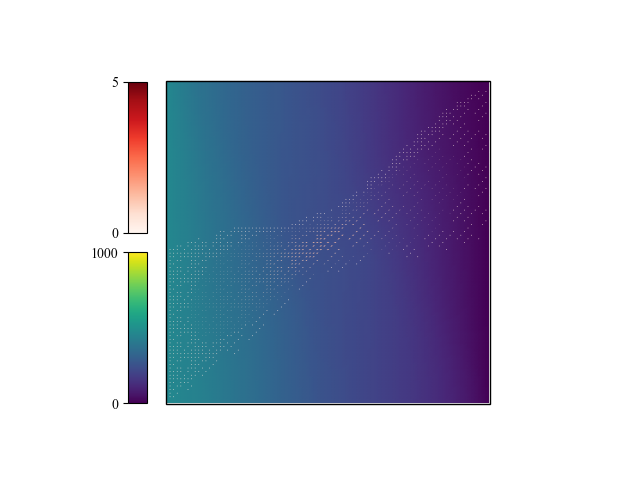} }\\
    \hspace{1.25cm}\subfloat[Approximation with ten modes.]{\includegraphics[trim=108 50 100 50, clip, width=0.44\linewidth]{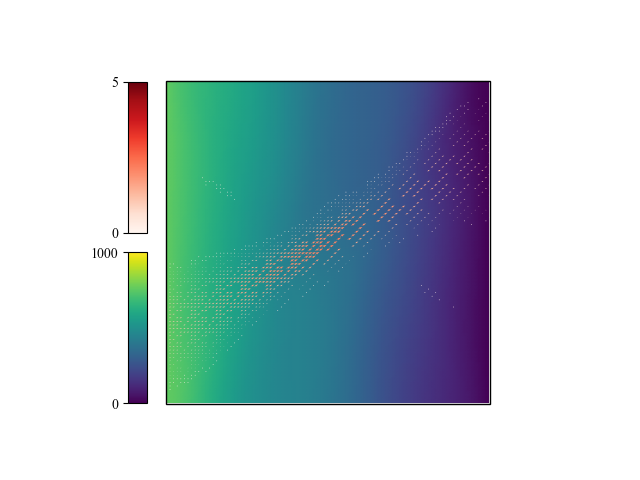}  }
    \subfloat[Approximation with twenty modes.]{\includegraphics[trim=108 50 100 50, clip, width=0.44\linewidth]{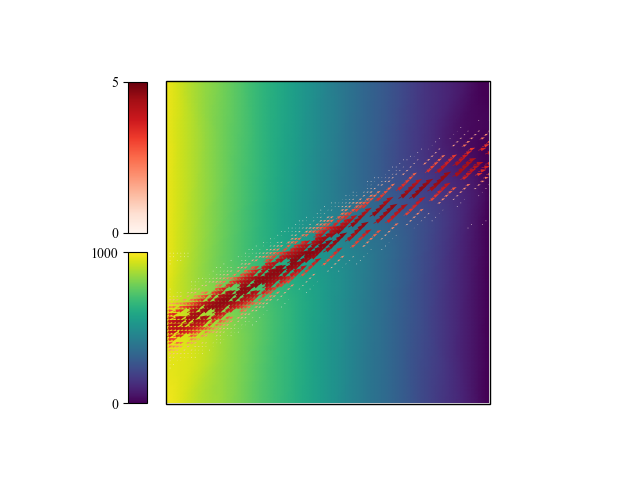} }
    \caption{Reduced order model simulation result for benchmark 2 with $\beta=30^\circ$.\\[-0.7cm]}
    \label{fig:B2RBmeth}
    \vspace{-0.3cm}
\end{figure}

\begin{figure}[H]
    \centering
    \subfloat[True solution.\mbox{\hspace{-1cm}}]{\includegraphics[trim=90 70 120 75, clip, width=0.54\linewidth]{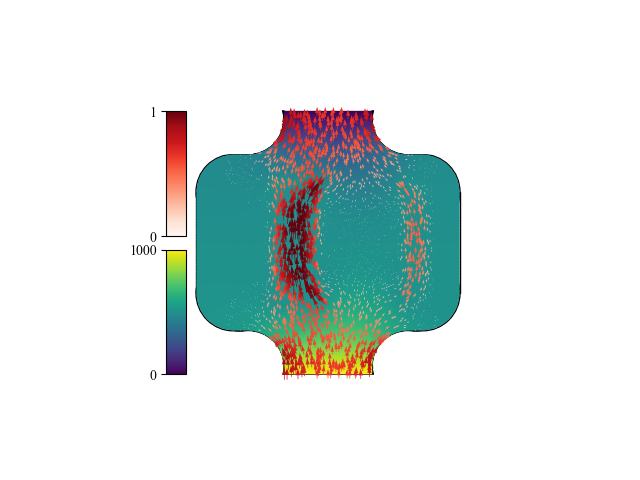} }\hspace{0cm}
    \subfloat[Approximation with five modes.]{\includegraphics[trim=140 70 120 75, clip, width=0.44\linewidth]{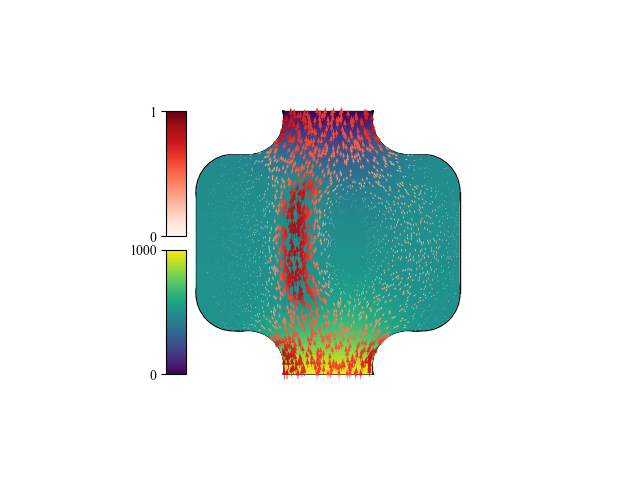} }\\\hspace{1.5cm}
    \subfloat[Approximation with ten modes.]{\includegraphics[trim=140 70 120 75, clip, width=0.44\linewidth]{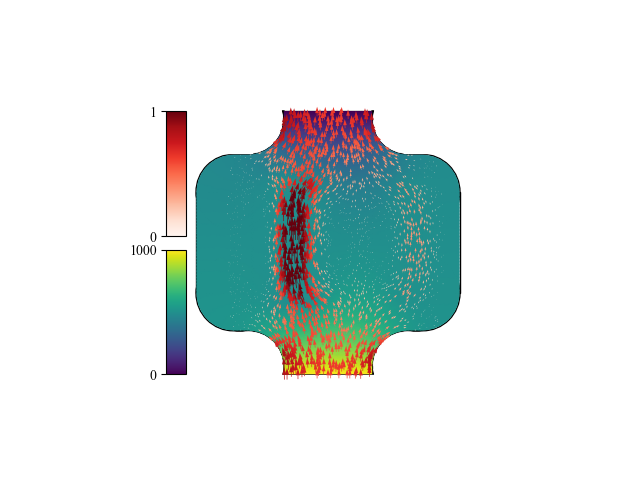}  }
    \subfloat[Approximation with twenty modes.]{\includegraphics[trim=140 70 120 75, clip, width=0.44\linewidth]{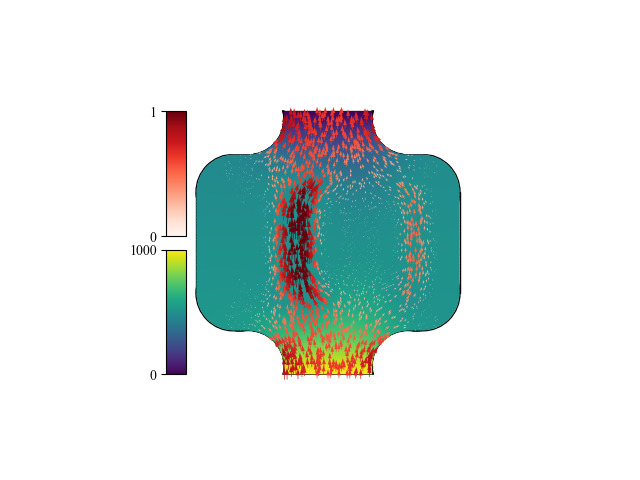} }
    \caption{Reduced order model simulation result for benchmark 3 with $\beta=[0.35,-0.25,0.1]$.}
    \label{fig:B3RBmeth}
\end{figure}

\Cref{fig:B1RBmeth,fig:B2RBmeth,fig:B3RBmeth} show example results of the complete reduced order model for all three benchmark problems for various numbers of reduced basis functions (i.e., various $r$). In all cases, 50 modes were used for the DEIM approximation, which we consider the maximal number for maintaining computational efficiency. 
For all three benchmarks, 20 solution modes are sufficient to produce reduced basis approximations that are sufficiently accurate in the ``eye-ball norm'' for these particular parameter points. Only upon careful inspection of the velocity magnitudes are we able to distinguish the approximations from the high-fidelity computations. Interestingly, benchmark 1 and benchmark 3 produce seemingly acceptable results already for ten solution modes, whereas benchmark 2 clearly does not.

The remaining open issue is the relation between the number of DEIM modes for each of the three fields, the number of solution modes in the reduced basis, and the fidelity of the reduced order model throughout the parameter space. Recall that $\epsilon(n)$ plotted in \cref{fig:B1singularvalues,fig:B2singularvalues,fig:B3singularvalues,fig:Uepsilon} represents the lowest possible relative error that may be achieved when approximating all solutions in the snapshot matrix with the first $n$ modes. This does not guarantee that the DEIM approximation and the reduced order model yield solutions close to that optimum. Additionally, the $\epsilon$ measure concerns an average over the parameter space. This says little about the maximum relative error, which is arguably the more important measure. Still, these computed $\epsilon$ values may guide the design of the reduced order model in terms of number of modes for the various reduced order approximations.

\begin{figure}[!b]
    \centering
    \subfloat[Naive choice of number of DEIM modes.]{\includegraphics[trim=0 0 0 0, clip,width=0.48\linewidth]{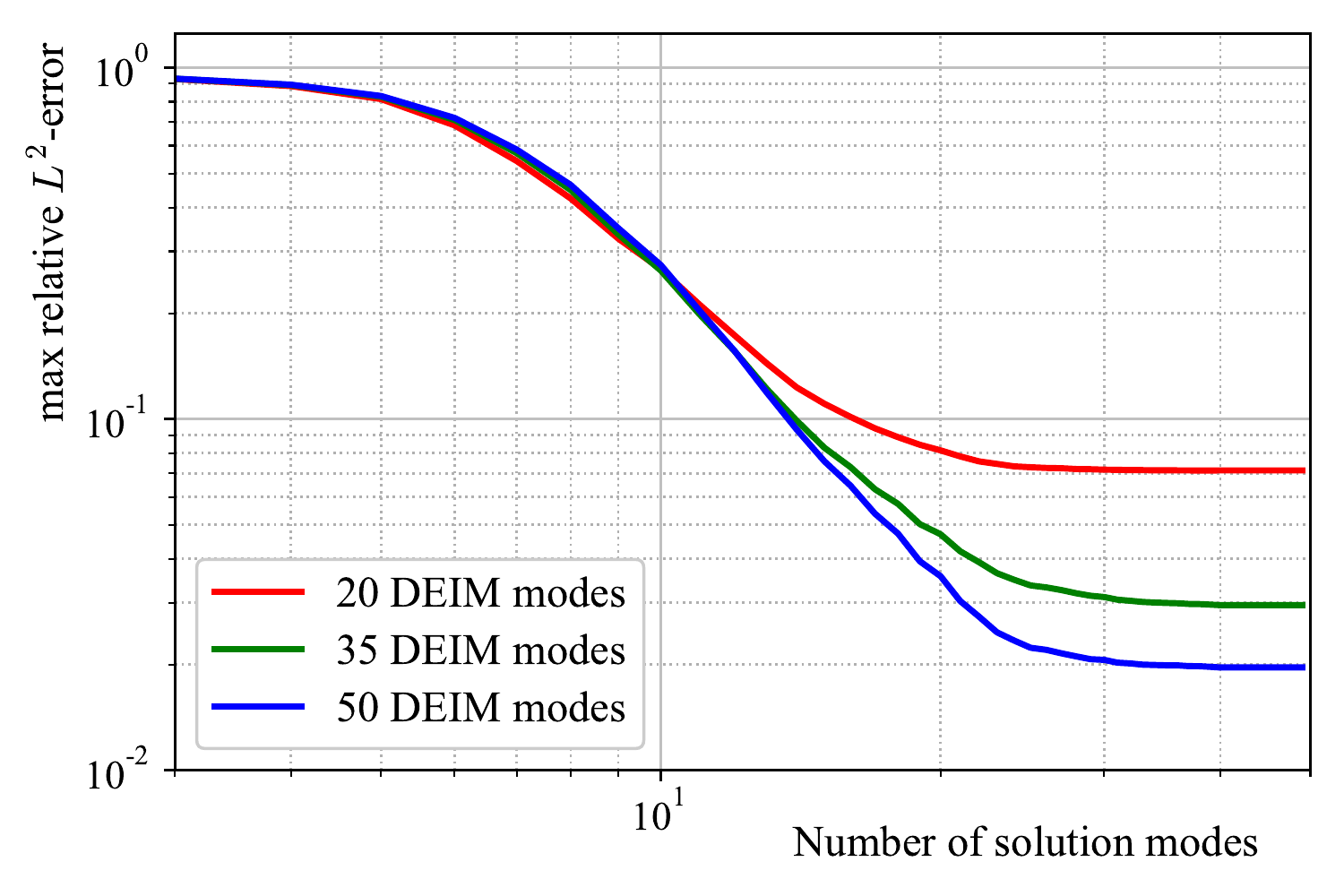}\label{fig:B1L2conva} }\hspace{.2cm}
    \subfloat[Guided choice of number of DEIM modes.]{\includegraphics[trim=0 0 0 0, clip,width=0.48\linewidth]{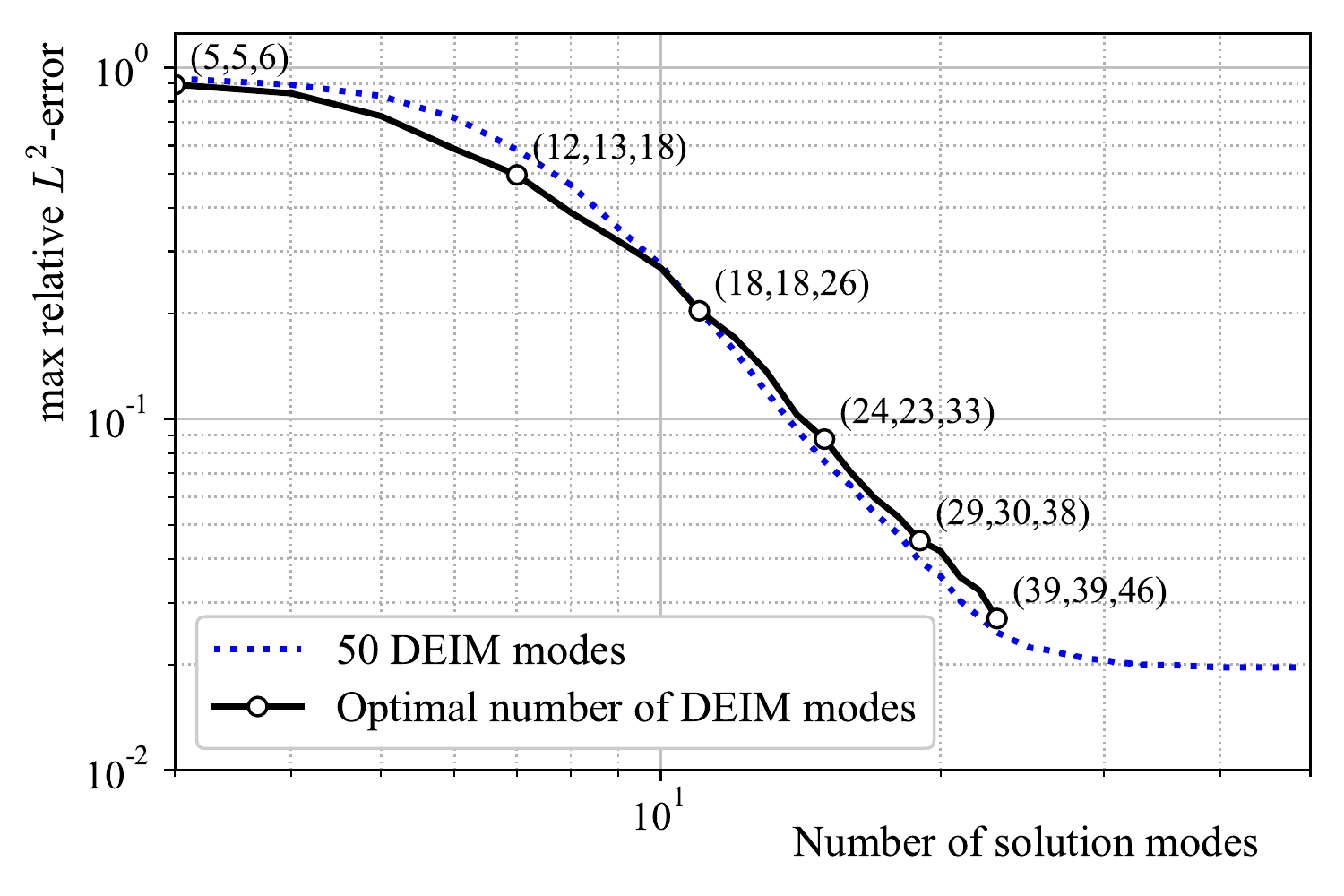} \label{fig:B1L2convb}}
    \caption{Convergence of the relative maximum $L^2$-error over the complete parameter space as a function of the number of reduced basis modes for different number of DEIM modes for benchmark 1.}
    \label{fig:B1L2conv}
\end{figure}

\begin{figure}[!t]
    \centering
    \subfloat[Naive choice of number of DEIM modes.]{\includegraphics[trim=0 0 0 0, clip,width=0.48\linewidth]{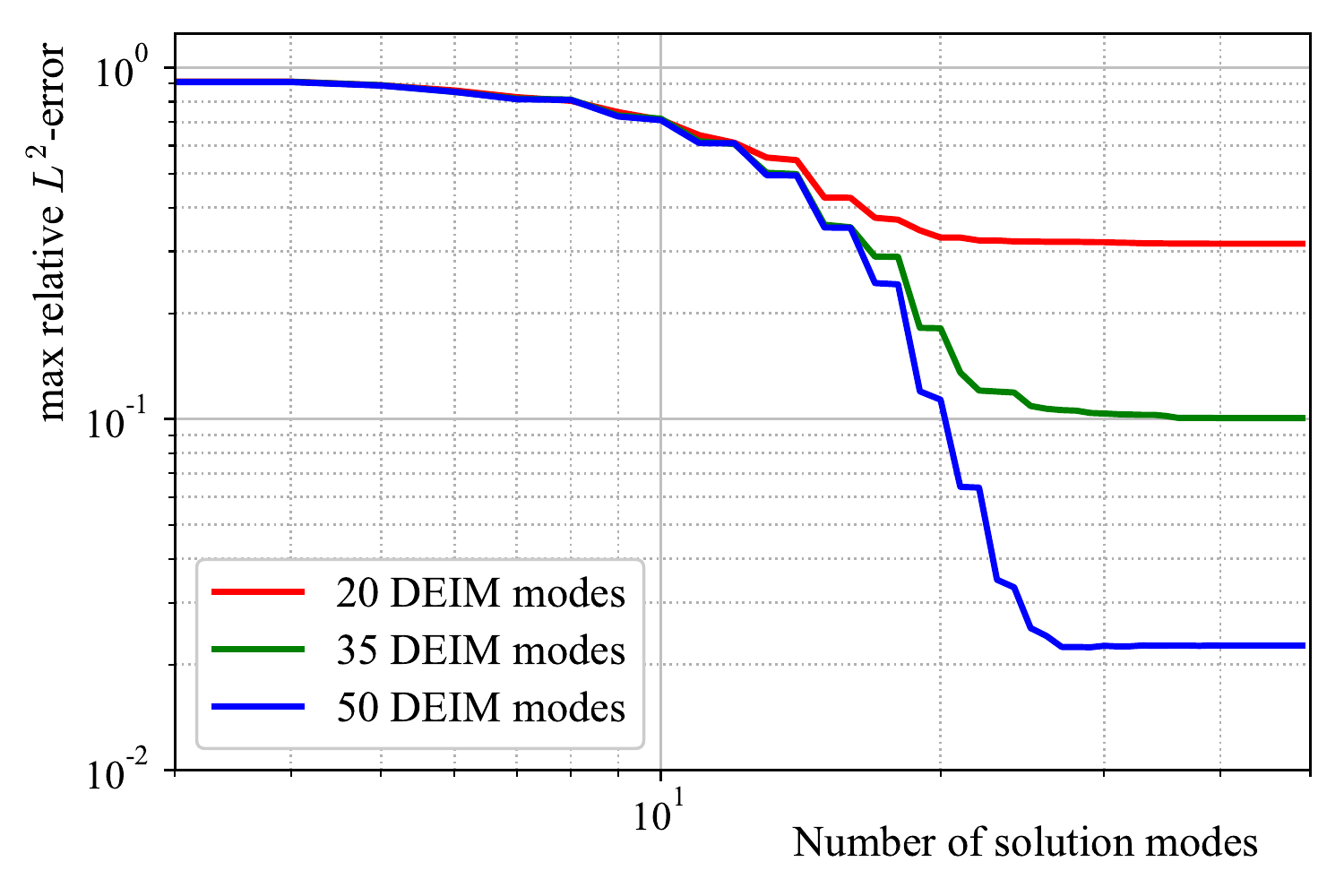} \label{fig:B3L2conva} }\hspace{.2cm}
    \subfloat[Guided choice of number of DEIM modes.]{\includegraphics[trim=0 0 0 0, clip,width=0.48\linewidth]{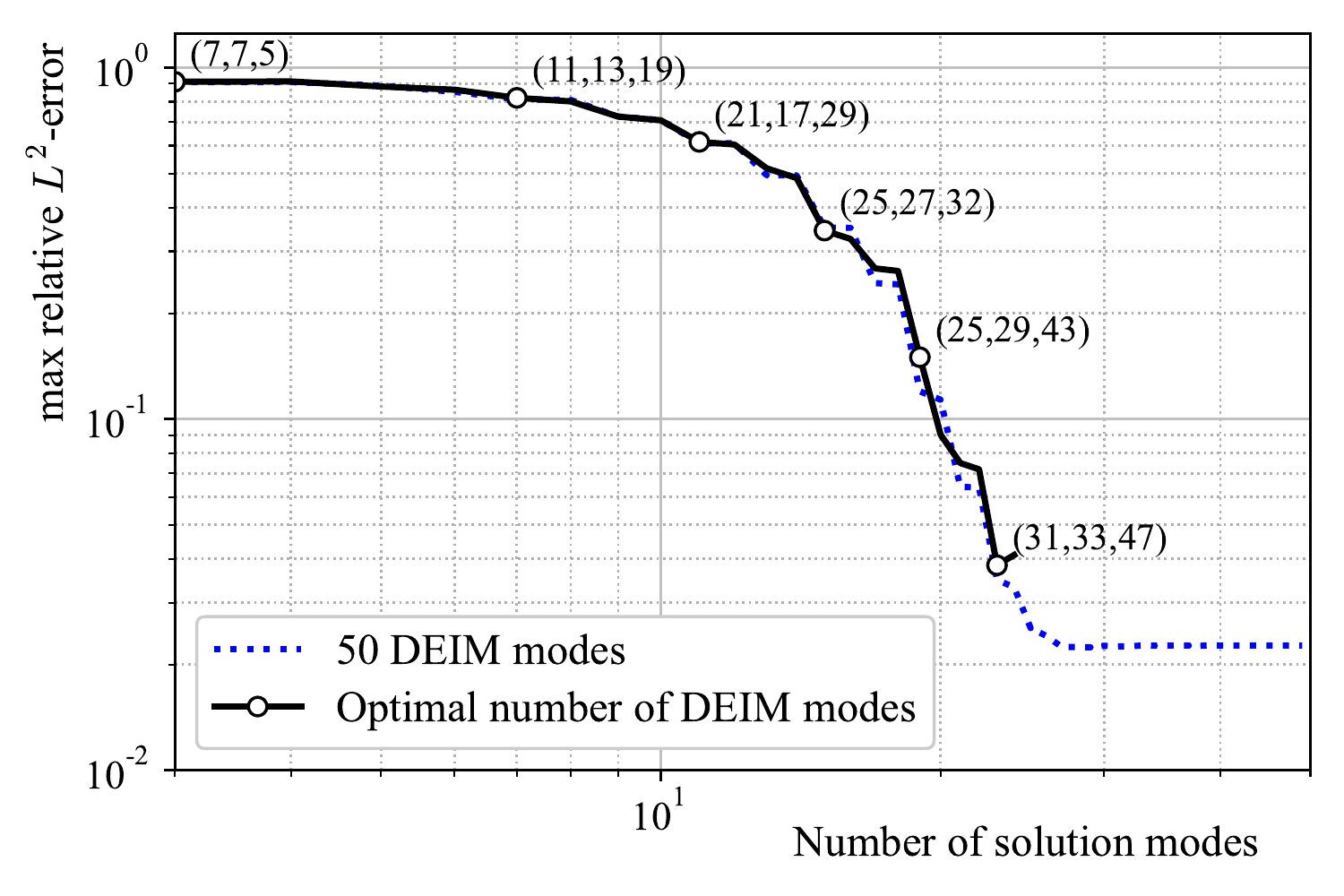} \label{fig:B3L2convb}}
    \caption{Convergence of the relative maximum $L^2$-error over the complete parameter space as a function of the number of reduced basis modes for different number of DEIM modes for benchmark 2.}
    \label{fig:B2L2conv}
\end{figure}

\begin{figure}[!t]
    \centering
    \subfloat[Naive choice of number of DEIM modes.]{\includegraphics[trim=0 0 0 0, clip,width=0.48\linewidth]{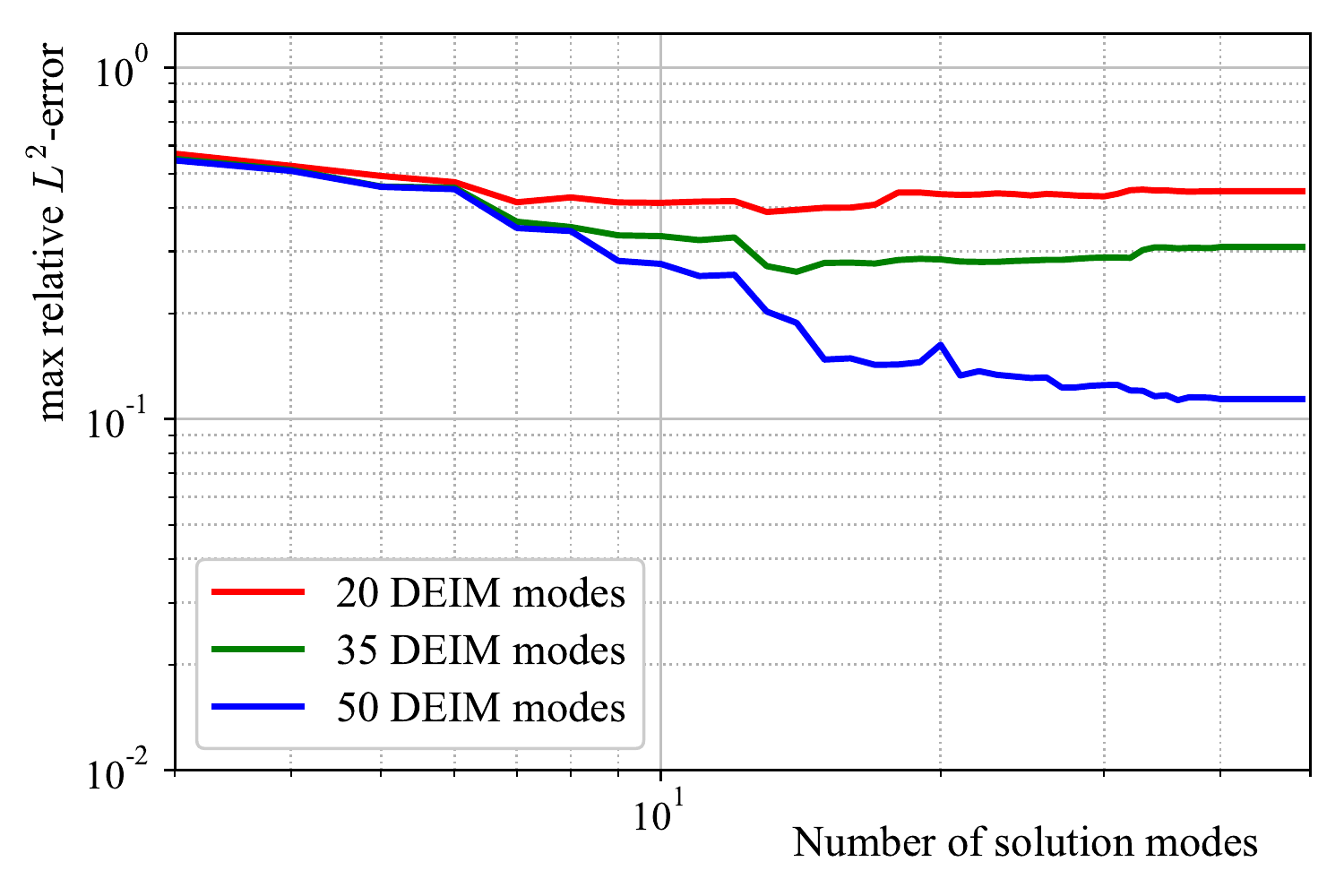} \label{fig:B2L2conva} }
    \subfloat[Guided choice of number of DEIM modes.]{\includegraphics[trim=0 0 0 0, clip,width=0.48\linewidth]{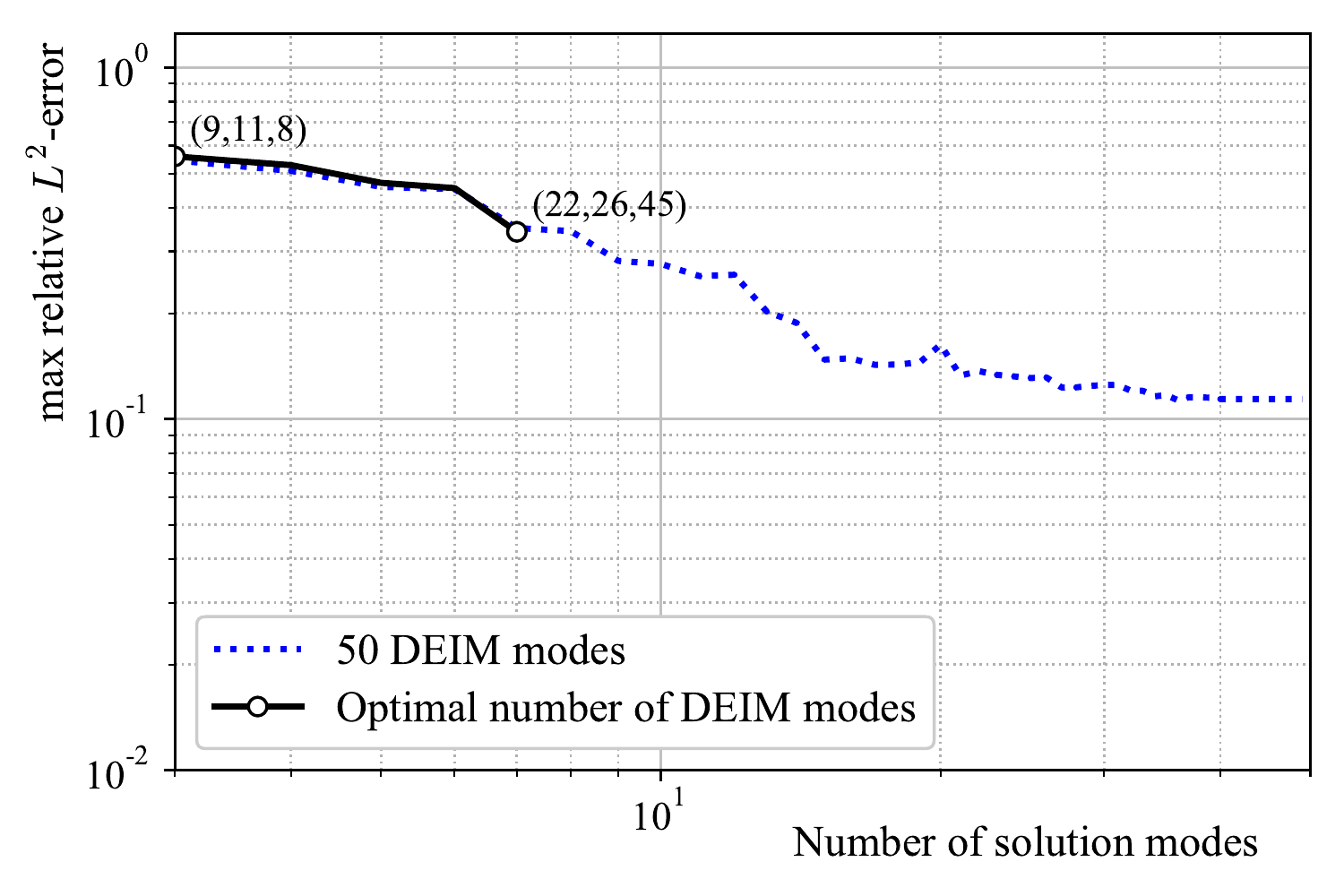} \label{fig:B2L2convb} }
    \caption{Convergence of the relative maximum $L^2$-error over the complete parameter space as a function of the number of reduced basis modes for different number of DEIM modes for benchmark 3.}
    \label{fig:B3L2conv}
\end{figure}

As a first indicator of the quality of the reduced order model, and its dependence on the DEIM approximation of the phase-field geometry representation, we compute the maximum $L^2$-error of the velocity field of the reduced order model over the complete (discrete) parameter space as a function of $r$ for various numbers of DEIM interpolation modes. The results are plotted in \cref{fig:B1L2conva,fig:B2L2conva,fig:B3L2conva} for each of the three benchmark problems. The general trend in each of the figures is the same: for small $r$ values, the maximum $L^2$-error is independent of the number of modes for the DEIM approximations. As $r$ increases above a certain threshold, the $L^2$-error starts to drop. This drop plateaus at different values for the different numbers of DEIM modes. For the first two benchmark problems, an acceptable maximum relative error of $\scriptstyle{\sim}$3$\%$ is achieved for fewer than 25 solution modes. For benchmark 3, however, the maximum error does not drop below 10\%. Apparently, 50 DEIM modes are insufficient to accurately represent the internal geometry. Indeed, the $\epsilon$-values for the phase-field quantities of benchmark 3, previously graphed in \cref{fig:B3singularvalues}, decayed significantly slower than those of benchmarks 1 and 2.


In \cref{sec:DEIM}, we show that the DEIM approximation of the three phase-field related quantities performs differently well. At the same time, the error in the reconstruction of the three fields may have different impact on the final solution error. These observations imply that a different number of modes may be required for the DEIM approximation of the three fields. A more targeted and optimized choice of the number of DEIM interpolation modes could reduce the computational expense of operating the reduced order model without adversely affecting the quality of the final solution. We propose to equate the $\epsilon$-value for the particular choice of number of solution modes (i.e., those in \cref{fig:Uepsilon}) to the average of the $\epsilon$-values for the DEIM interpolation modes (i.e., those in \cref{fig:B1singularvalues,fig:B2singularvalues,fig:B3singularvalues}). That is, we use the minimal values $n_\xi$,  $n_\zeta$ and $n_t$ for the DEIM reconstruction of the respective fields such that:
\begin{align}\label{opt_modes}
    \epsilon_{RB}(n) = \frac{1}{3}\big( \epsilon_{\xi}(n_\xi) +  \epsilon_{\zeta}(n_\zeta) + \epsilon_{t}(n_t)  \big)\,,
\end{align}
which we solve by iteratively incrementing $n_\xi$,  $n_\zeta$ or $n_t$ depending on which reduces the right-hand-side of \cref{opt_modes} most. We limit the number of DEIM modes for each field to 50 to avoid excessive computational expense.

With the choice of number of DEIM interpolation modes for $\xi$, $\zeta$ and $\B{t}$ based on \cref{opt_modes}, we recompute the $L^2$-errors throughout the entire parameter space. The results are graphed for each of the three benchmark problems in \cref{fig:B1L2convb,fig:B2L2convb,fig:B3L2convb}. The graphs markers at intermediate points in the graphs state the predicted number of DEIM interpolation modes as $(n_{\xi},  n_{\zeta}, n_{t})$-tuples. The figures also include a copy of the blue 50-modes lines from \cref{fig:B1L2conva,fig:B2L2conva,fig:B3L2conva} as references. We observe that the maximum $L^2$-errors of the optimized approach closely follows the (``overkill'') 50-modes line while substantially reducing the required number of DEIM interpolation modes. The optimized mode number lines break shortly before the 50-modes lines plateau, confirming that more than 50 DEIM modes are required for at least one of the fields to maintain a drop of error. These results indicate that \cref{opt_modes} provides an effective estimation for the required number of DEIM modes for all fields such that the DEIM approximations do not dominate the source of error of the reduced order model.

\newpage

\section{Example application: in-situ leach mining}
\label{sec:SolutionMining}

\begin{figure}[!t]
    \centering
    \subfloat[3-spot pattern]{\includegraphics[width=0.32\linewidth]{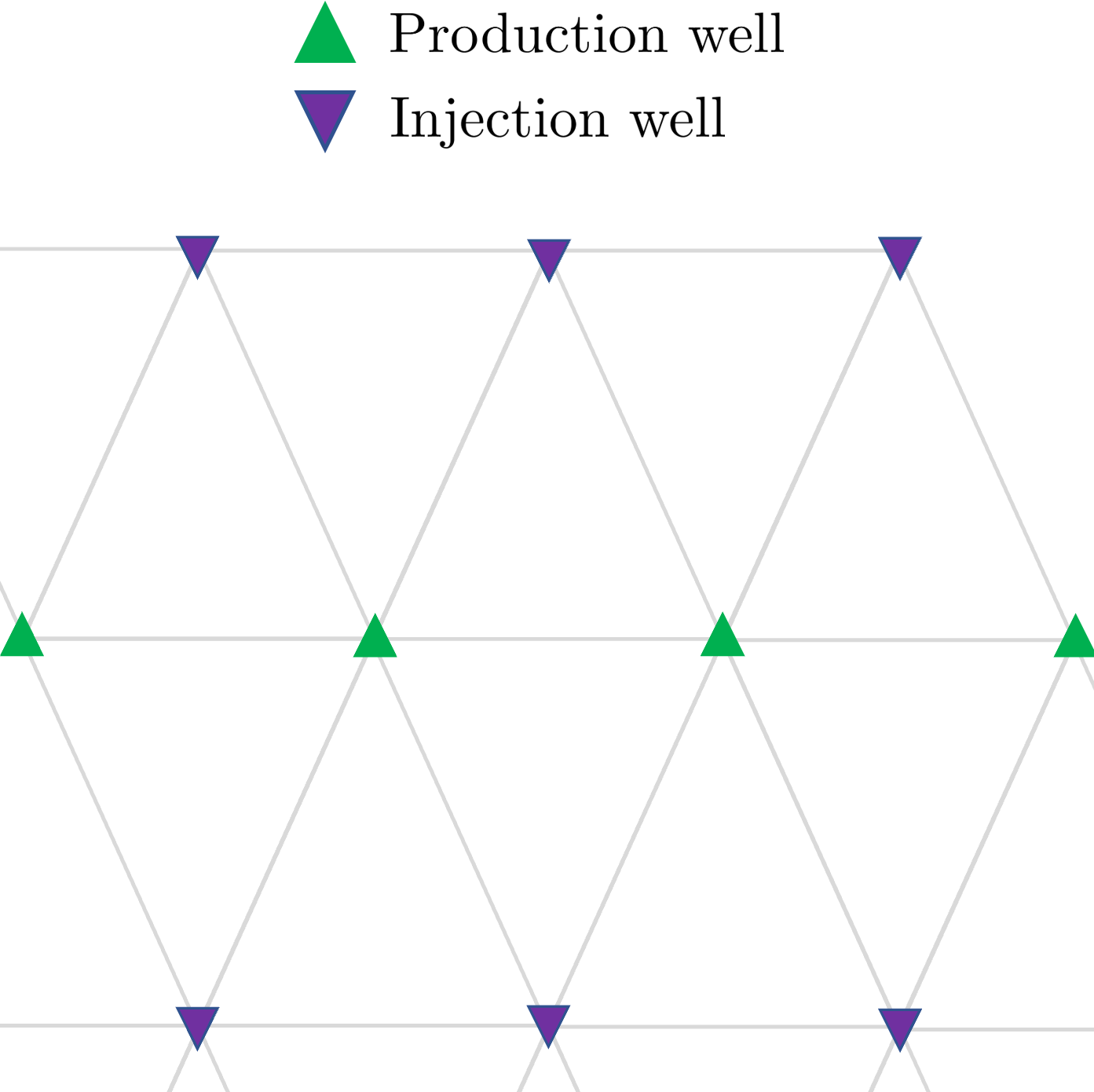} }\hfill
    \subfloat[5-spot pattern]{\includegraphics[width=0.25\linewidth]{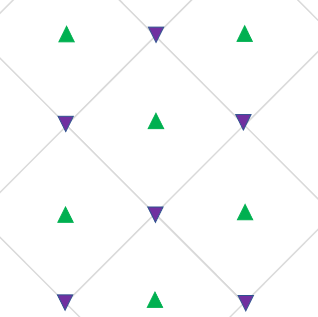} }\hfill
    \subfloat[7-spot pattern]{\includegraphics[width=0.27\linewidth]{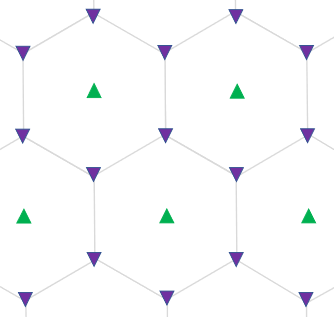} }
    \caption{Patterns of injection and production wells used in practice for ISL mining.\\[-0.25cm]}
    \label{fig:wellPatterns}
\end{figure}

In-situ leaching (ISL), also known as solution mining, is a mining process for recovering underground minerals such as copper and uranium. In 2019, ISL accounted for 57\% of the total uranium mining worldwide \cite{worldNuclear}. 
During ISL, a solution (typically a mixture of native groundwater, a complex agent and an oxidant) is pumped through an injection well down to an ore deposit. 
At this depth, the solution flows towards the production well, while dissolving minerals from the ore. The pregnant solution is then pumped to the surface at the production well, where the dissolved uranium is later extracted from the solution. 
The injection and production wells 
are distributed in regular patterns. 
The most frequently used pattern types are the 
3-spot, 5-spot, and 7-spot patterns, shown in \Cref{fig:wellPatterns}. 
The choice of pattern depends on factors such as subsurface permeability, deposit type, ore grade and installation cost. 
The distance between the injection and productions wells can be as high as 50m$-$60m for the 3-spot pattern, and as low as 15 - 30m for the 5-spot and 7-spot patterns. As a result, the 5-spot and 7-spot patterns have a higher installation cost but increased uranium recovery rate and operation flexibility. 

Challenges that are faced while constructing and operating an ISL mining site include (i) ensuring that nearby groundwater is not contaminated by only operating sections of the mine and continuously taking measurements at monitoring wells surrounding the mine, and (ii) managing the permeability of the ore deposit with hydraulic fracturing or controlled explosives. In both cases, a predictive simulation tool with fast response time could facilitate decision making processes, may serve to optimize long-term operation of the mining site and enables more targeted interventions. We make use of the reduced order model described in the previous sections to develop such a simulation tool. We illustrate its operation capability by solving the inverse problem of predicting the damage state of the ore deposit based on inflow and outflow measurements at the injection and production wells.

\subsection{The high-fidelity model}
\label{sec:hexagonModel}

\begin{figure}[!t]
    \centering
    \vspace{1cm}
    \subfloat[Proposed well field design, adapted from \cite{honeymoon}.]{\includegraphics[width=1\linewidth]{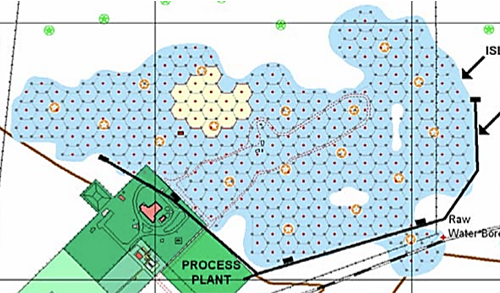} \label{fig:honeymoonSchematica}}\\
    \subfloat[Hexagon pattern.]{\includegraphics[width=0.9\linewidth]{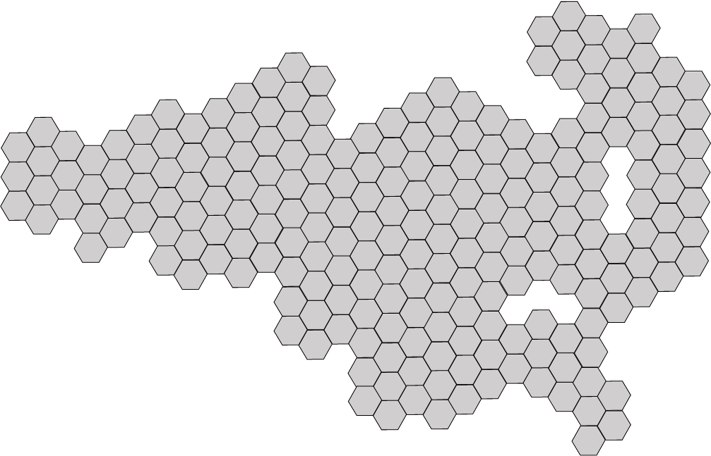} \label{fig:honeymoonSchematicb}}
    \caption{Well field of the Honeymoon mining project. It contains 225 hexagons.}
    \label{fig:honeymoonSchematic}
\end{figure}

In the following, we examine the Honeymoon mining project illustrated in \cref{fig:honeymoonSchematica}.
This mine is located in South Australia and was active from 2011 to 2013, after which operations were halted due to the low price point of uranium \cite{honeymoon}. The mine features a 7-spot pattern. A distinctive property of this pattern 
is its tiling into hexagonal segments, illustrated in \cref{fig:honeymoonSchematicb}. The complete mining site consists of 225 such hexagons. We exploit this regular pattern in our computational model for the subsurface flow throughout the entire mining site: we develop a reduced basis method that can represent the flow in a single hexagon, and reuse this model in each tile.
For each hexagon, the distance between the single production well to each of the six injection wells is 15m. The diameter of the injection and production wells is 0.5m. 
The uranium-bearing sand lies 100 - 150m underground, under sheets of gravel, clay, and sand. We treat it as coarse-grained pyritic sand, with a permeability of 100 millidarcy ($10^{-13}$m$^2$). Since the solution consists largely of groundwater, we assume a fluid viscosity of 1mPa$\cdot$s. The boundary conditions that we assume for each hexagon are a variable positive pressures at the injection wells, a zero pressure at the production well and symmetry conditions on the remaining boundaries. These conditions correspond to the ones from \cref{sec:Diffuse} (i.e., \cref{BCs}).

\begin{figure}[!t]
    \centering
    \vspace{-0.5cm}
    \subfloat[$D = 0.25$\mbox{\hspace{-1.5cm}}]{\includegraphics[trim=70 40 150 40, clip,width=0.39\linewidth]{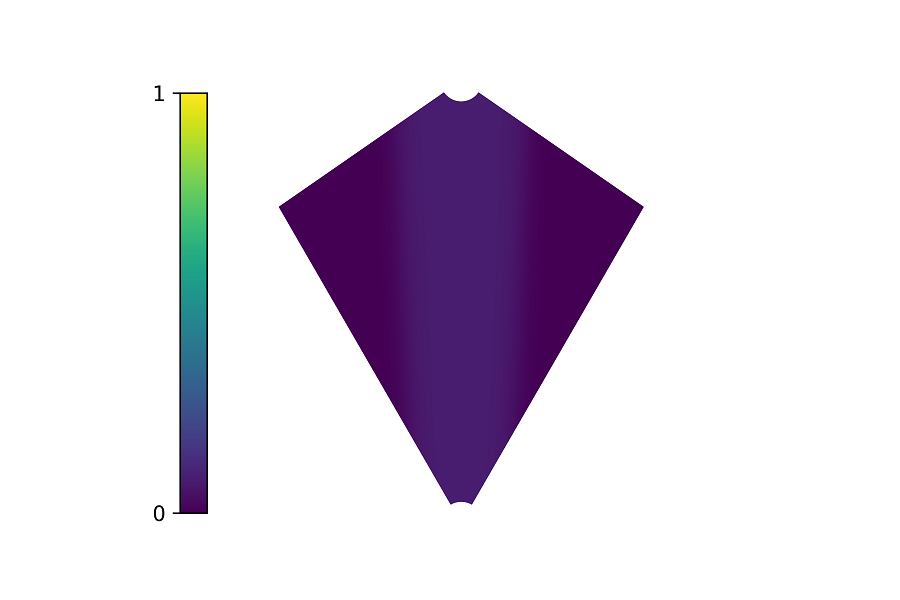}}\hfill
    \subfloat[$D = 0.5$\mbox{\hspace{0.25cm}}]{\includegraphics[trim=160 40 140 40, clip,width=0.293\linewidth]{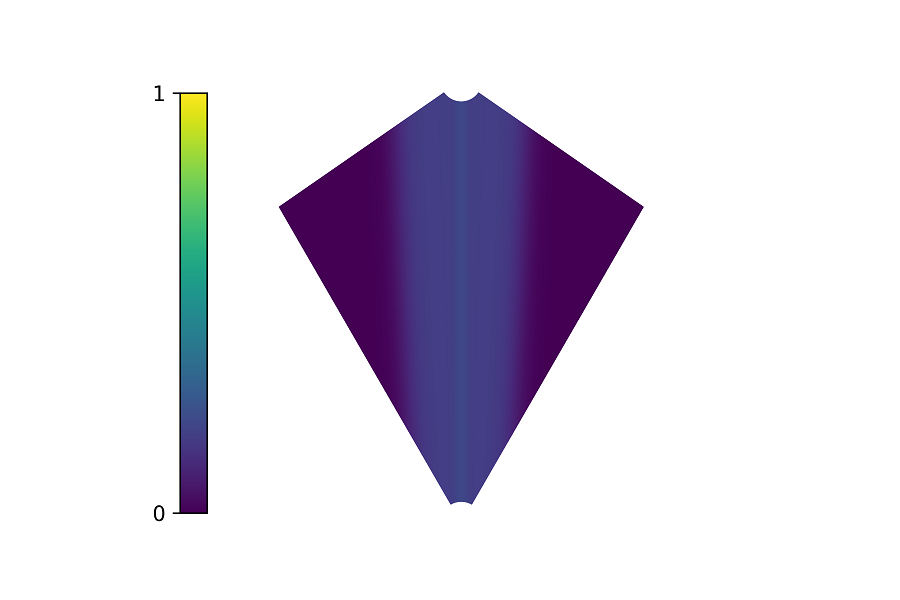}}\hfill
    \subfloat[$D = 1.0$\mbox{\hspace{0.5cm}}]{\includegraphics[trim=160 40 140 40, clip,width=0.293\linewidth]{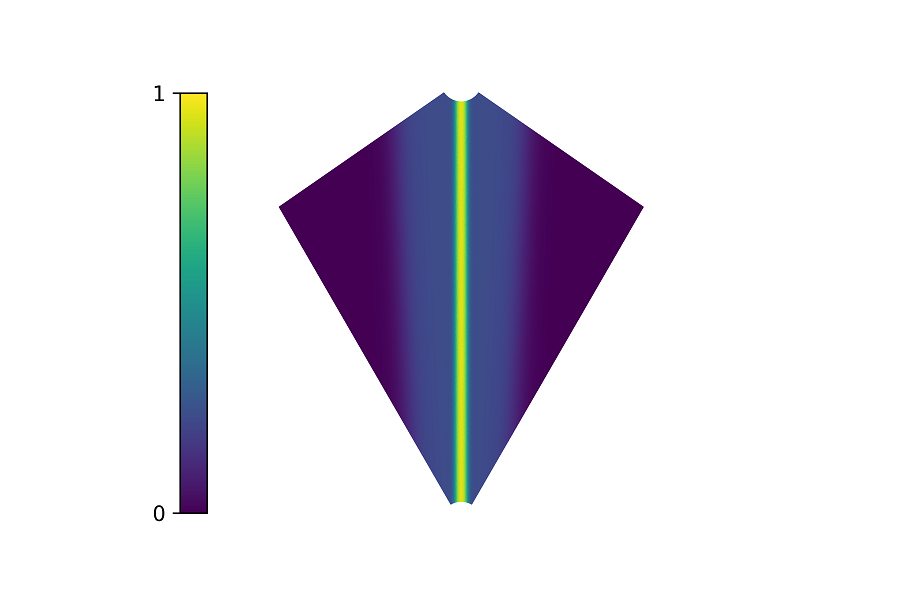}}\\
    \caption{Damage field through sextant for different damage values $D$.}
    \label{fig:channelPatterns}
\end{figure}

\begin{figure}[!t]
    \centering
    \vspace{-0.5cm}
    \subfloat[$D = 0.25$\mbox{\hspace{-1.5cm}}]{\includegraphics[trim=70 40 160 40, clip,width=0.39\linewidth]{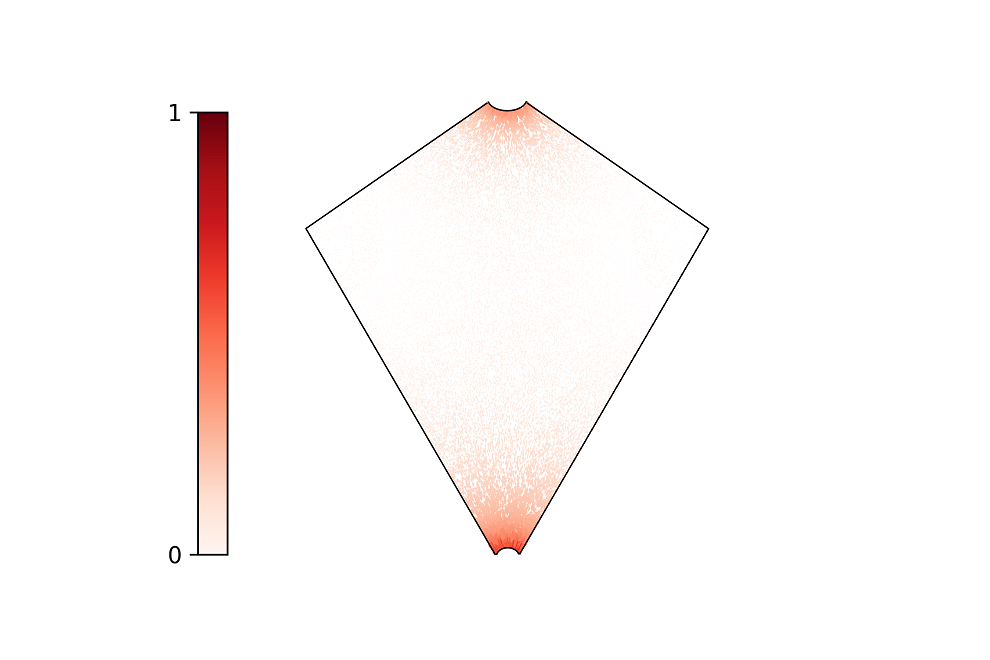}}\hfill
    \subfloat[$D = 0.5$\mbox{\hspace{0.25cm}}]{\includegraphics[trim=170 40 150 40, clip,width=0.293\linewidth]{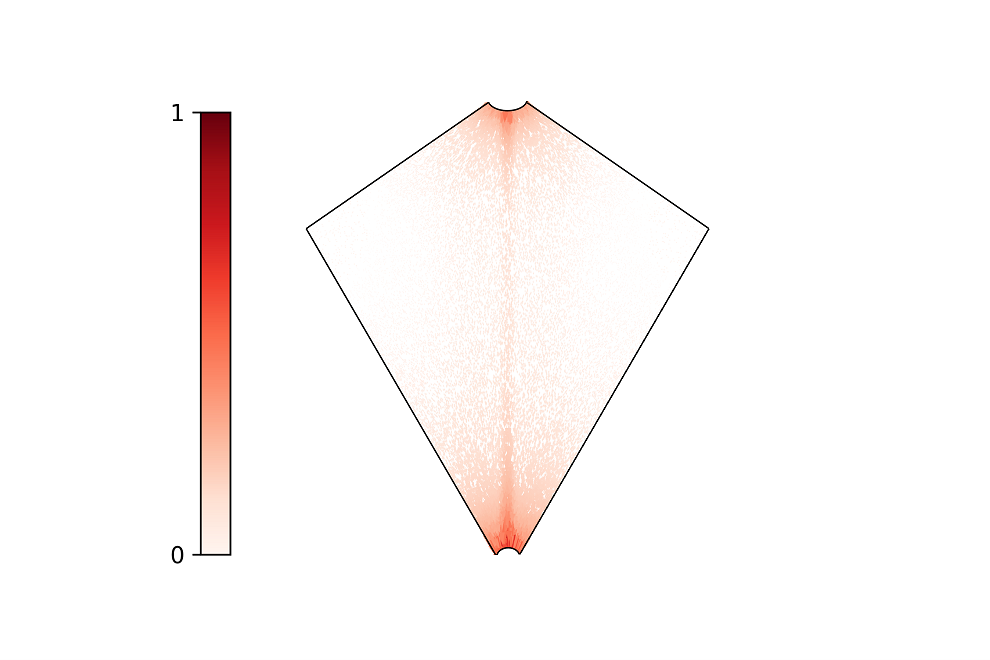}}\hfill
    \subfloat[$D = 1.0$\mbox{\hspace{0.5cm}}]{\includegraphics[trim=180 40 140 40, clip,width=0.293\linewidth]{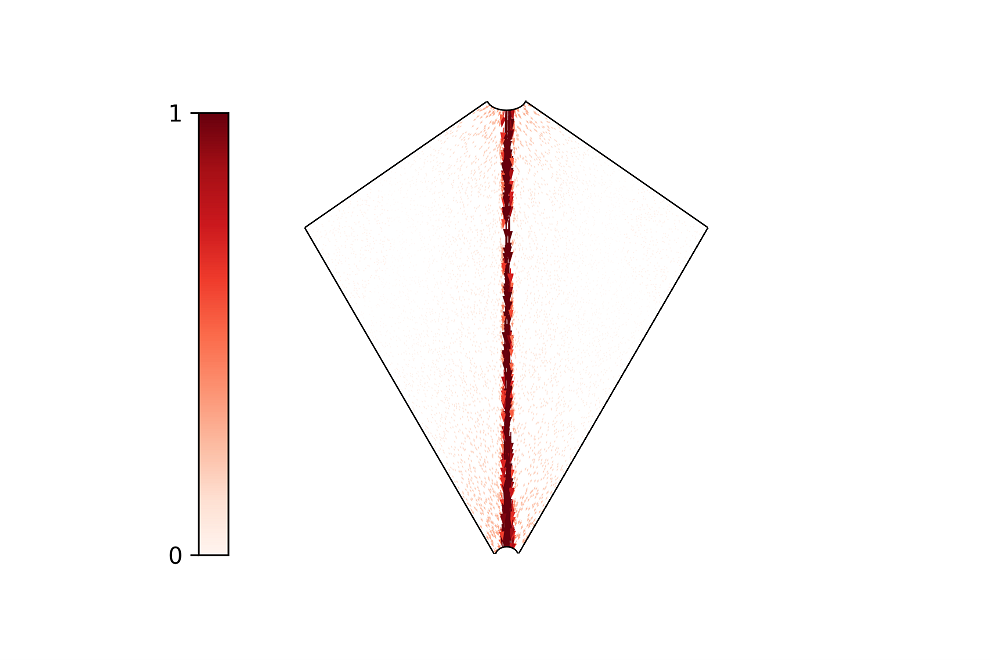}}
    \caption{Velocity profile through sextant for different damage values $D$.}
    \label{fig:damagePatterns}
\end{figure}

The rate at which uranium can be leached depends not only on the permeability of the sandstone but also 
on the crack patterns running through it.
These cracks create pathways in the ore deposit for the solution to penetrate. Once cracks become too large they reduce the effectiveness of the uranium recovery. To include the crack patterns in our model, we construct an evolving damage field that runs from each of the six injection wells to the production well. The damage field depends only on a single parameter, the damage parameter $D$, and represents the severity of the cracked state for each sextant of the hexagon. The damage parameter varies between 0 and 1, signifying an undamaged state or a fully developed crack from injection well to production well.
\Cref{fig:channelPatterns,fig:damagePatterns} show the evolving damage field and the corresponding velocity profile for $D$ between 0 and 1 in one sextant. For $D=0$ there is no damage and the resulting velocity profile follows from pure Darcy flow. While $D$ increases it creates a mixed Stokes/Darcy domain, which can be interpreted as a homogenization of smaller cracks in the sandstone. The resulting velocities in this created channel are therefore higher than in the pure Darcy domain. This difference increases with higher values of $D$. At $D=0.5$ a thinner channel on top of the existing one is developing and represents the growth of a single straight crack. At $D=1$ this thin crack is fully developed, represented by a domain of Stokes flow.

\subsection{The reduced order model in one hexagon}
\label{ssec:ROM}

Computing finite element approximations of the high-fidelity model on each of the 225 hexagons would yield an excessive amount of degrees of freedom for the complete system. Solving the resulting discrete system would be too time consuming for optimization purposes or for solving inverse problems. Moreover, for each new damage field in each different hexagon, the system would have to be reassembled, adding to the severe computational expense. Hence, we create a reduced basis method with a very limited number of degrees of freedom within each hexagon. To ensure that our reduced order model is capable of accurately representing the possible solution states up to a tolerance level, we perform the procedure discussed at the end of \cref{sec:RB} to determine the required number of DEIM and RB modes.

\begin{figure}[!t]
    \centering
    \subfloat[Phase-field quantities on one hexagon.]{\includegraphics[trim=0 0 0 0, clip,width=0.48\linewidth]{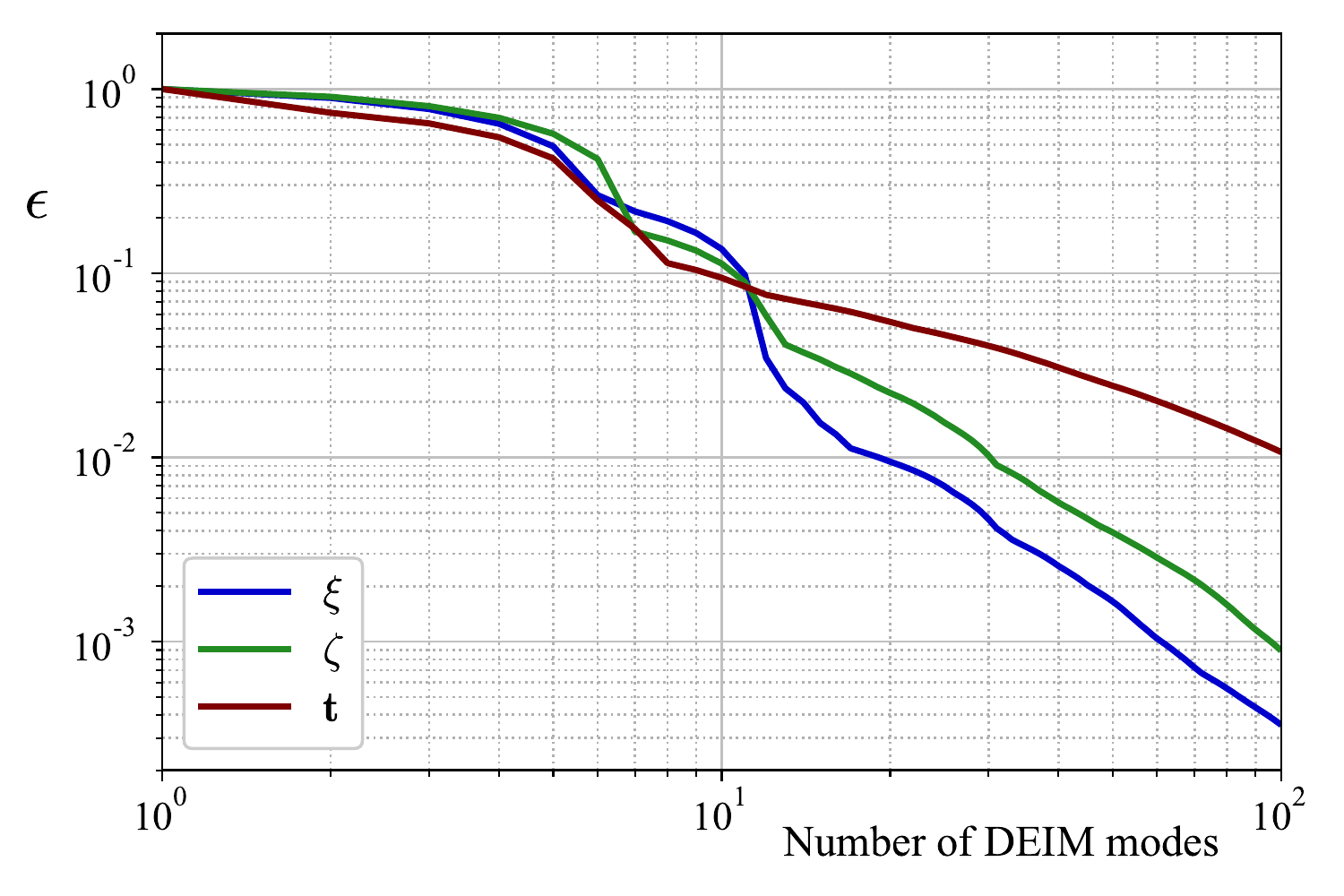} \label{fig:ISLsingularvalues_pf} }\hspace{.2cm}
    \subfloat[Solution fields on one hexagon.]{\includegraphics[trim=0 0 0 0, clip,width=0.48\linewidth]{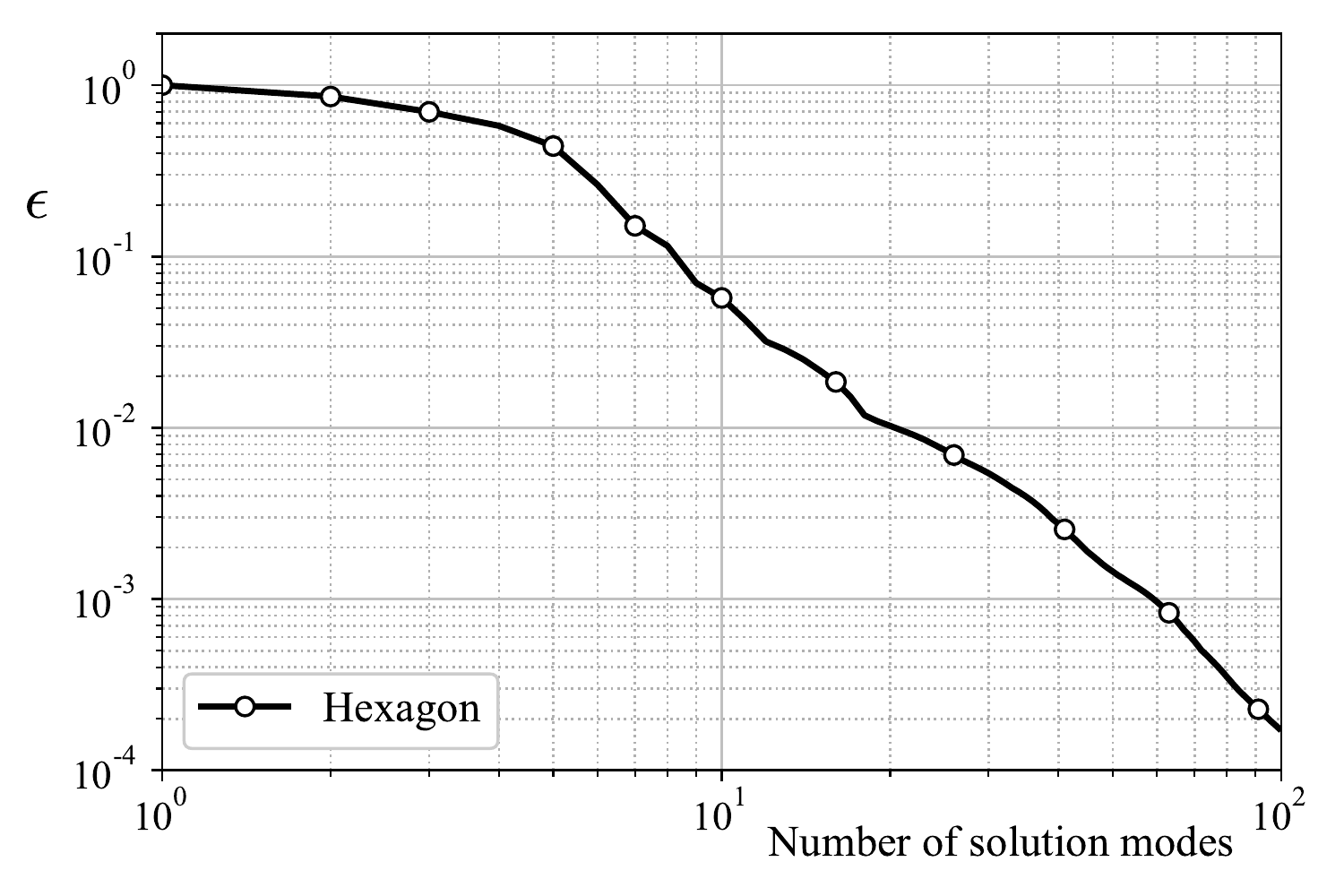} \label{fig:ISLsingularvalues_sol}}
    \caption{Convergence of $\epsilon$ as a function of the number of modes.}
    \label{fig:ISLsingularvalues}
\end{figure}

In each hexagon, there are six parameters running from 0 to 1 that determine the overall damage state. Additionally, there are six pressure boundary conditions, which we vary between 0 and $10^5$. The resulting parameter space is thus twelve-dimensional: $\mathbb{P}=[0,1]^6\times[0,10^5]^6$. To explore the solution manifold, we combine a tensor product sampling grid with a random sampling set, giving $\mathbb{P}^h = \mathbb{P}^h_{U} \bigcup \mathbb{P}^h_{R}$. We construct $\mathbb{P}^h_{U}$ as $\{ 0.25 , 0.5, 0.75, 1 \}^6\times\{ 10^5 \}^6$ such that it contains $4096$ points representing different internal geometries (damaged states). To also include the variable pressure, we construct $\mathbb{P}^h_{R}$ from another 6000 points that are sampled with a uniform random distribution along every axis in $\mathbb{P}$. 
For each parameter point in $\mathbb{P}^h$, we compute the phase-field related fields $\xi$, $\zeta$, $\B{t}$, as well as the solutions $(\B{u}^h,p^h)$ on the high-fidelity mesh. \Cref{fig:ISLsingularvalues_pf,fig:ISLsingularvalues_sol} show the $\epsilon$ measures from \cref{errormeasure} as a function of the number of modes for the phase-field quantities and the solution fields respectively. Based on these results we choose $25$ solution modes per hexagon and, using the guided choice discussed above,
$(21, 30, 48)$ for the DEIM modes respectively.

With this set-up, we obtain the following reduced order model in each individual hexagon:
\begin{subequations}
\begin{align}\label{WFHex}
& \text{Find }(\B{u}^{RB},p^{RB})\in RB^r \text{ s.t. } \forall\,\, (\B{v},p)\in RB^r: \nonumber\\
\begin{split}
    &\int\limits_{\Omega} 2 \phi \mu \nabla^s \B{u}^{RB} : \nabla^s \B{v} \dO  + \int\limits_{\Omega} (1-\phi) \mu \B{\kappa}^{-1} \B{u}^{RB}\cdot \B{v} \dO - \int\limits_{\Omega}  p^{RB} \, \nabla \cdot \B{v}  \dO + \int\limits_{\Omega}  \B{\alpha} \B{u}^{RB} \,  \cdot \B{v}  \dO  \\[-0.2cm]
    & \hspace{2cm}= \sum_{j=1}^6 \, \int\limits_{\partial\Omega_{i}} p_j \,\B{n} \cdot \B{v} \dO \,,
\end{split} \\
&\int\limits_{\Omega} q\, \nabla \cdot \B{u}^{RB}  \dO = 0 \,,
\end{align}
\end{subequations}
or, in matrix representation:
\begin{align}
\underbar{\underbar{K}}(D_1,\cdots,D_6)\, \hat{\underbar{u}}^{RB} = \sum_{j=1}^6 p_j \underbar{F}_j\,, 
\end{align}
where, of course, the parameter dependence of $\underbar{\underbar{K}}$ is affine due to the DEIM approximation of the geometry representation.

\subsection{Connecting reduced order models}
\label{ssec:Connecting}

Due to the symmetry conditions between the hexagons, there is no mass-flow through any of the boundary facets, nor is there a shear stress acting on those facets. These homogeneous conditions are satisfied by each of the solution snapshots in the snapshot matrix, and thus by each of the reduced basis functions. All coupling between the hexagonal cells then occurs through shared connectivity with the injection well at each of the corners. The available measurement data at each of these injection wells is the total amount of fluid that is pumped into the underground system from the surface. Since the distribution of the mass-flow among the neighbouring hexagons is unknown, this produces a coupled system of equations.

\begin{figure}[!t]
    \centering
    \subfloat[Adjacent hexagons]{\label{fig:coupledHexagonsa} \includegraphics[trim=0 -30 0 0,width=0.3\linewidth]{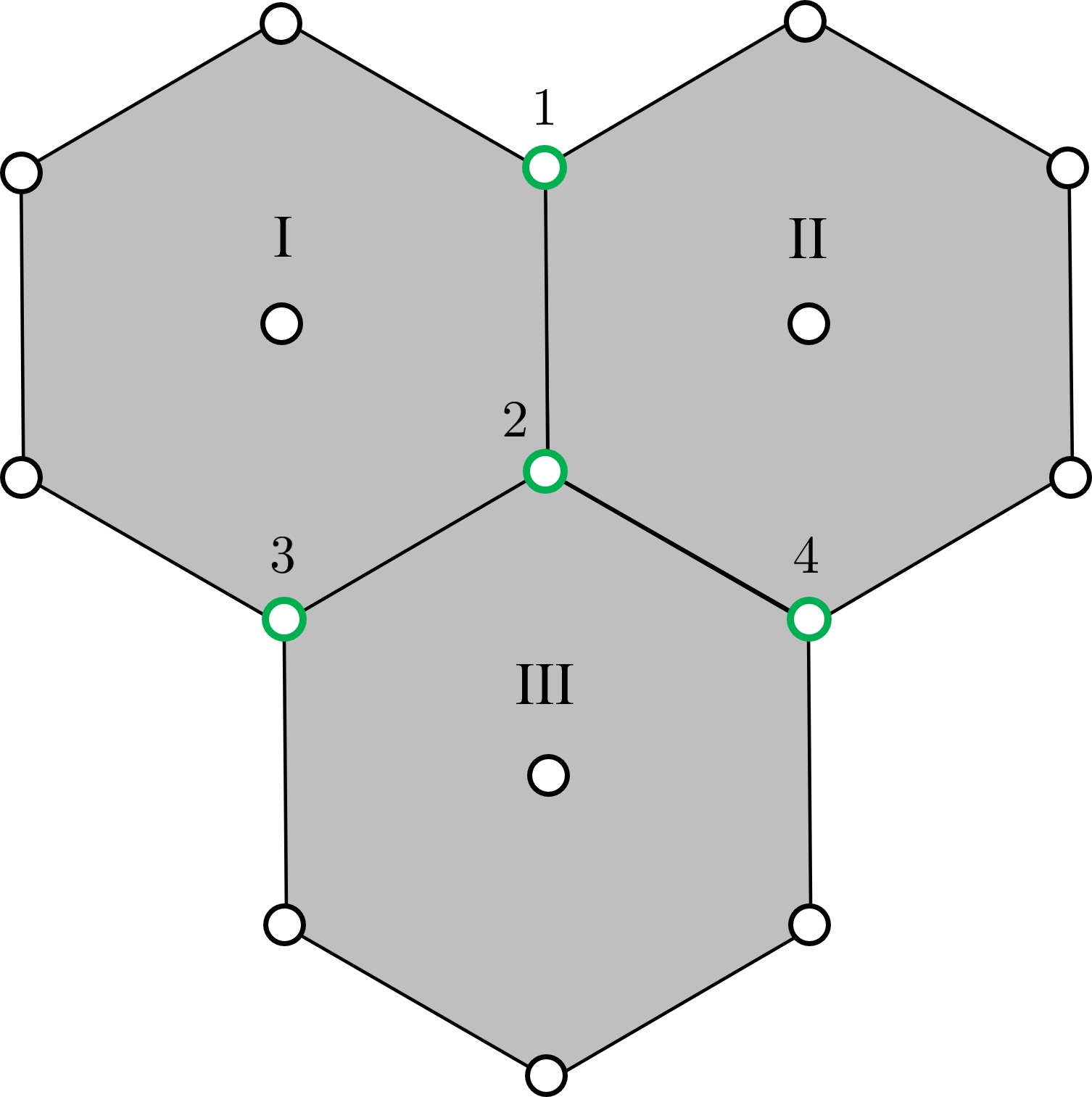} }\hspace{1cm}
    \subfloat[Coupling matrix]{\includegraphics[width=0.35\linewidth]{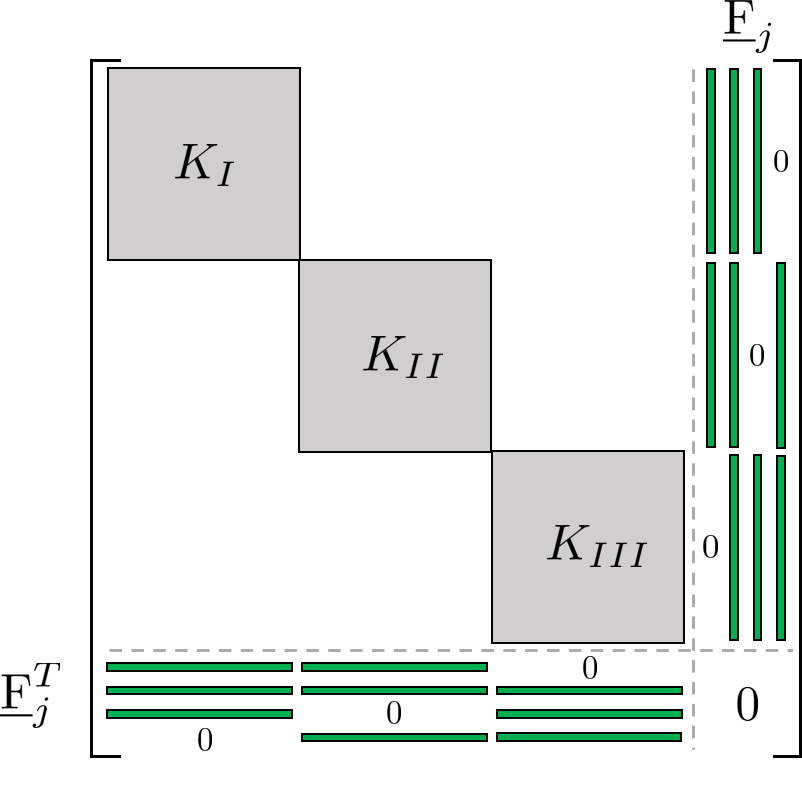}}
    \caption{Example of three adjacent hexagons with corresponding coupling matrix (green bars indicating the position of the coupling vectors $\underbar{F}_j$).}
    \label{fig:coupledHexagons}
\end{figure}

To introduce this connectivity, we add a Lagrange multiplier constraint that enforces a prescribed total inflow for each of the injection wells. For injection well 2 in the example of \cref{fig:coupledHexagonsa}, the variational statement corresponding to the inflow condition reads:
\begin{align}\label{inflowconstraint}
    q_2 u_{2,in} = \int\limits_{\partial\Omega_{\text{I},j}} q_2 \B{n} \cdot \B{u}_{\text{I}} \dS + \int\limits_{\partial\Omega_{\text{II},j}} q_2 \B{n} \cdot \B{u}_{\text{II}} \dS  +
    \int\limits_{\partial\Omega_{\text{III},j}} q_2 \B{n} \cdot \B{u}_{\text{III}} \qquad \forall\, q_2\dS 
\end{align}
where $q_2$ is the scalar test function that corresponding to the constraint. The inflow conditions of \cref{inflowconstraint} can be recognized as the transpose of the pressure boundary conditions of \cref{WFHex}. Indeed, the Lagrange multipliers corresponding to the constraint of $q_2$ is the pressure $p_2$, which now becomes a degree of freedom. 
The full system of equations is then assembled as follows. The 25 degrees of freedom for the reduced basis approximation within each hexagon form a $25 \times 25$ block on the diagonal. For each of the injection wells, another row and column are added to the system. The column and the row are populated with three $\underbar{F}_j$-vectors (respectively their transposes) corresponding to the three neighboring hexagons. This procedure is illustrated in \cref{fig:coupledHexagons}.

For the complete Honeymoon mining-site shown in \cref{fig:honeymoonSchematic}, this results in a system of equations of ``merely'' 6,077 degrees of freedom. Without the reduced basis representation of the flow in each hexagon, the full system would exceed 15 million degrees of freedom. For any given new subsurface damage state in any particular hexagon, the corresponding matrix $\underbar{\underbar{K}}(D_1,\cdots,D_6)$ has to be recomputed (which is a cheap operation due to its affine dependence on the local damage parameters), and then the block on the system matrix can simply be overwritten. As a result, we are able to recompute a solution state for a completely new damage parameter in around $0.1$s on an Intel i9-9900k @5Ghz, computing in serial. The majority of the compute time is spent on assembling the 225 diagonal blocks. Since these blocks are completely independent of one another, this computation could alternatively be performed in parallel very effectively. For our purposes, though, the $0.1$s computation time suffices. \\[-1cm]

\subsection{Solving an inverse problem: calibration of the subsurface state}

\begin{figure}[!t]
    \centering
    \includegraphics[width=0.85\linewidth]{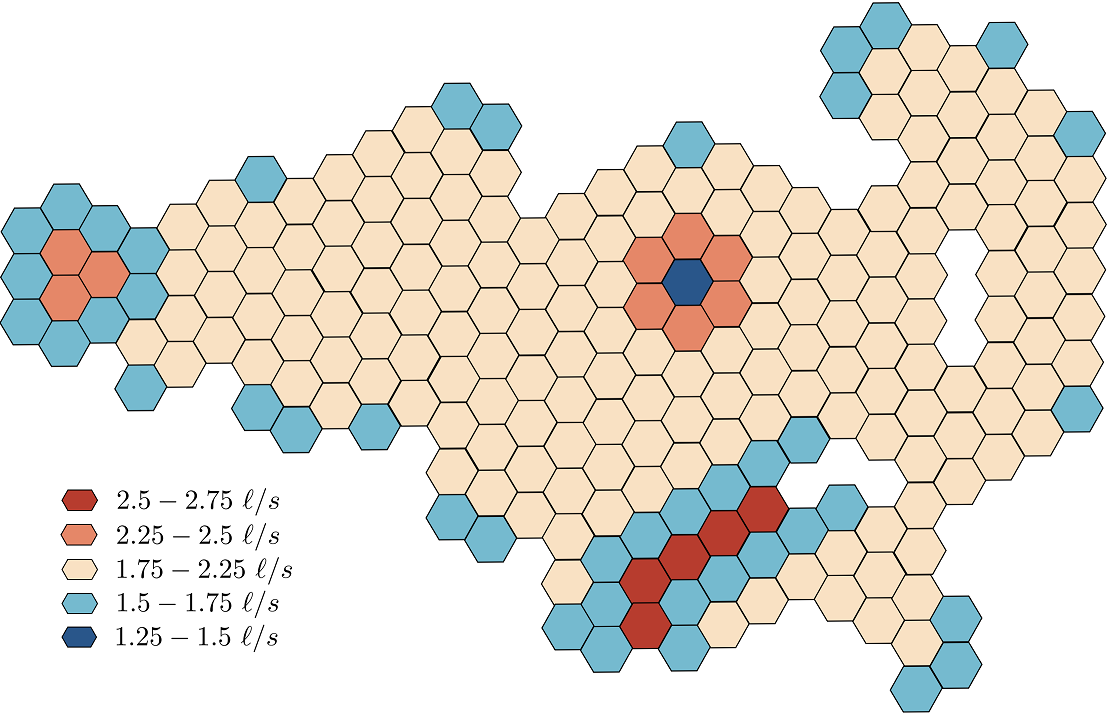}
    \caption{Assumed measurements of outflow at production wells over entire well field.\\[-0.2cm]}
    \label{fig:outflowMeasurements}
    \vspace{0.1cm}
\end{figure}

The subsurface state affects the flow profile, for which measurement data is available in terms of mass throughputs at the injection and production wells. Determining the subsurface state based on the measurement data thus involves solving an inverse problem. Based on \cref{sec:hexagonModel,ssec:ROM,ssec:Connecting}, we have a low-order computational model for the subsurface flow through the entire Honeymoon mining site of \cref{fig:honeymoonSchematic}. To produce an example problem, we artificially generate a distribution of inflow and outflow measurements. We do so by assuming smoothly varying inflow pressure values for the injection sites, ranging between $8.5\cdot 10^4$Pa and $1.0\cdot 10^5$Pa. We also assume a smoothly varying damage field. We overlay this damage pattern with random noise of about 1\% and add some larger damage values to represent local cracks. We then compute the outflow field in $\ell$/s of the described problem by multiplying the resulting throughput in each production well by a height of $5$m. This is the assumed intake height of the production wells inside the uranium deposits. Finally, we add random noise to these computed values to ensures that our model will not trivially be able to obtain the exact solution. The resulting incompatibility between inflow and outflow also occurs in actual mining sites due to leakage of the solute into different ground layers and due to extraction of groundwater. The final pattern of outflow measurements is shown in \cref{fig:outflowMeasurements}. 


\begin{figure}[!b]
    \vspace{-0.1cm}
    \centering
    \includegraphics[trim=0 0 0 0, clip,width=0.55\linewidth]{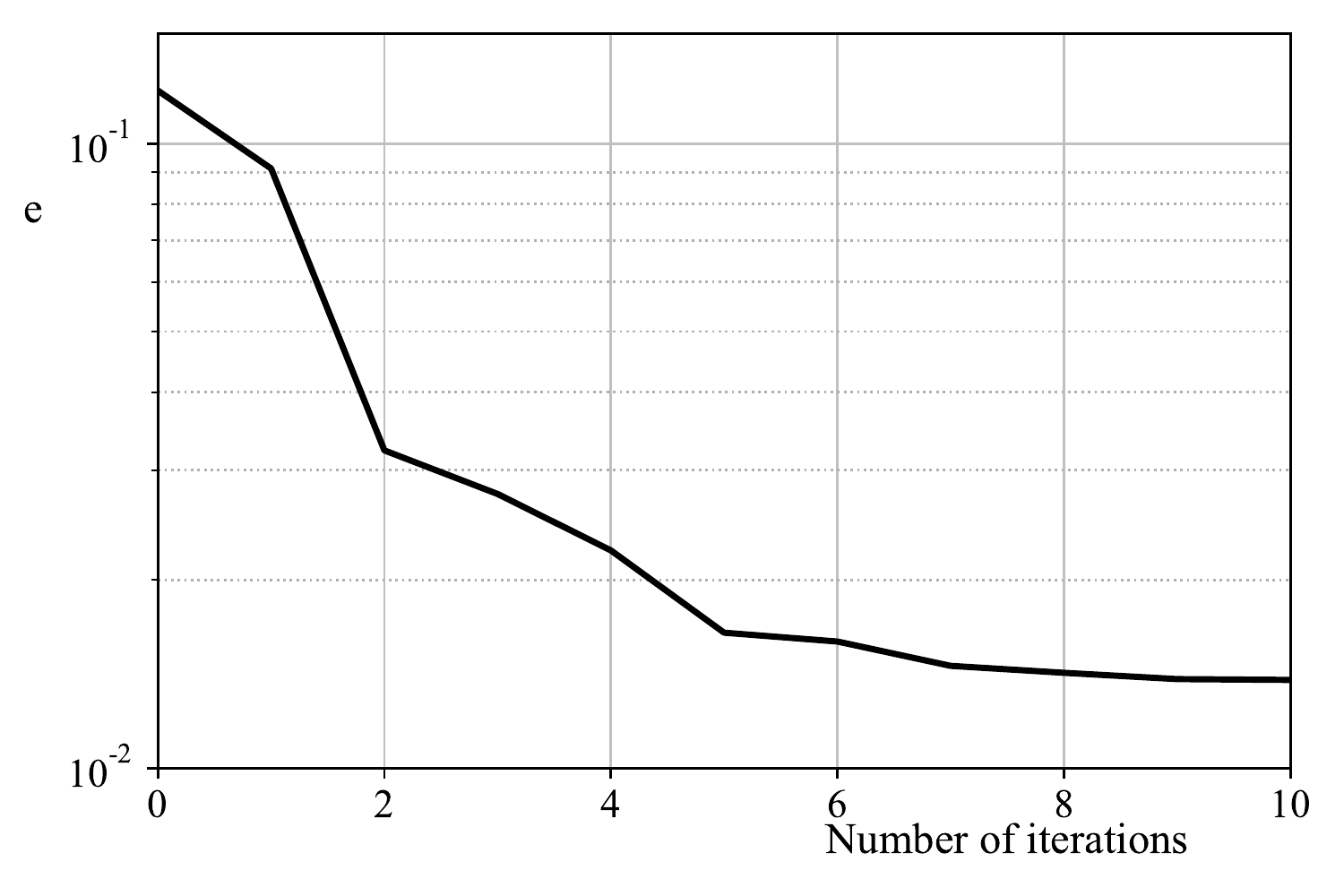}
    \caption{Error between measurements and computed outflows over the number of iterations.\\[-0.9cm]}
    \label{fig:convergence_optimization}
    \vspace{-0.9cm}
\end{figure}

Given only this measurement data, we now solve an optimization problem for the damage parameter $D$ in each sextant of each hexagon. 
The objective is to minimize the error between the given measurements of \cref{fig:outflowMeasurements} and the resulting outflow at each production well of the coupled computation. As an initial state, the damage parameter is homogeneously set to 0. We employ a gradient-descent method to iteratively adjust each damage parameter, until the model predicts the correct outflows for the given inflow. The gradient is computed numerically with a central difference finite difference scheme, which results in two function evaluations per sextant for changes in the one-dimensional damage parameter. Such an approach 
is only feasible due to the reduced-order model. 
To track the achieved accuracy of the scheme, we determine the error between the computed outflow and the measured outflow according to:
\begin{align}
    e = \sqrt{ \frac{1}{225} \sum_{i=1}^{225} \left( \frac{u_{i,out}-u^h_{i,out}}{u_{i,out}} \right)^2  }\,, \label{flowmeasure}
\end{align}
where $u_{i,out}$ is the measured outflow for production well $i$, and $u^h_{i,out}$ is the computed outflow. The convergence of this relative error with respect to the number of iterations is shown in \cref{fig:convergence_optimization}. As the graph indicates, we reach a root mean square relative error below $1.5$\% in under $10$ iterations. This residual error can be attributed to the artificially added effects of leakage and groundwater extraction.

\begin{figure}[!b]
    \centering
    \includegraphics[width=0.95\linewidth]{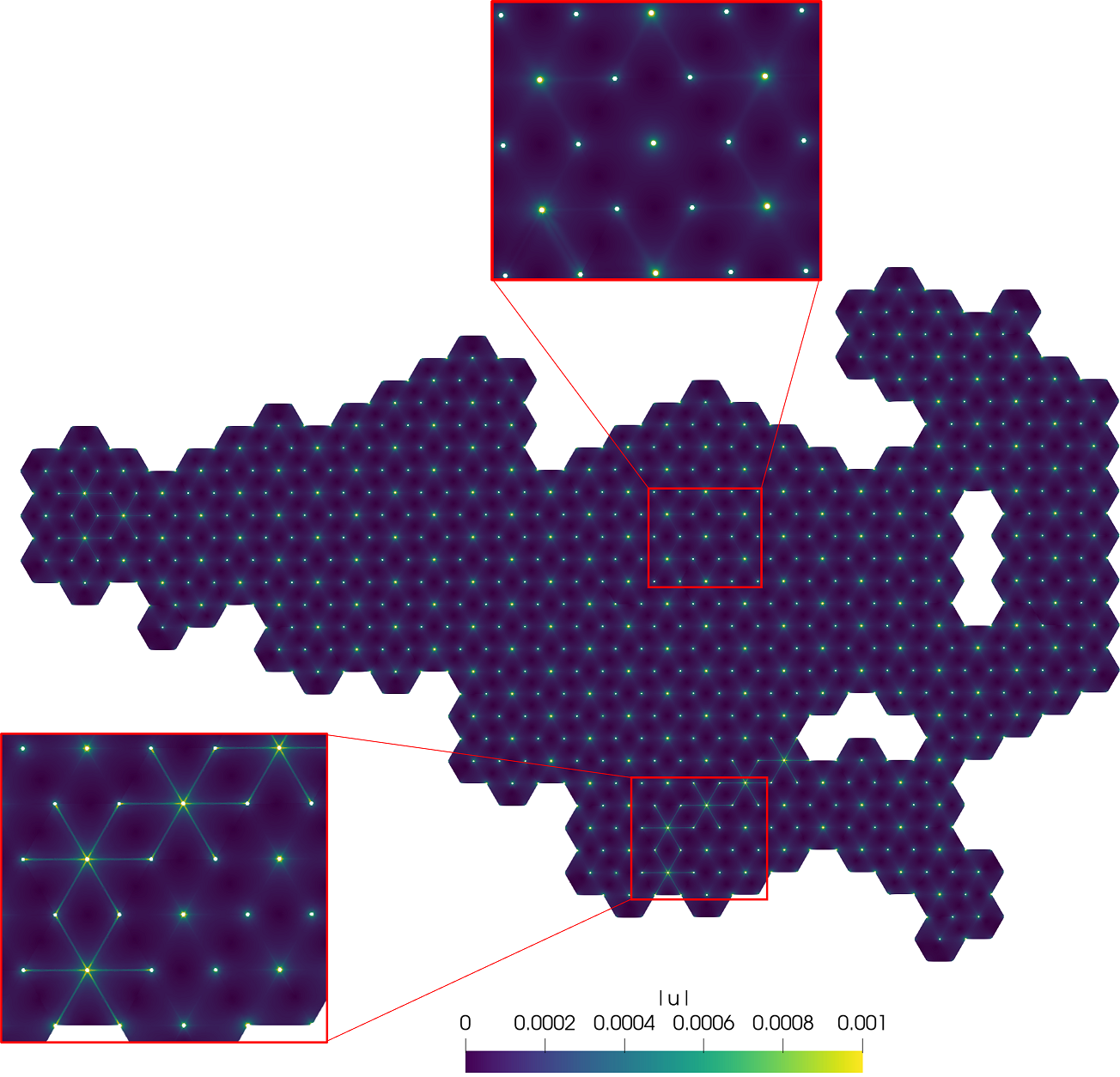}
    \caption{Velocity profile over entire well field with optimized damage parameters.\\[0.5cm]}
    \label{fig:honeymoonOptimized}
    \vspace{0.5cm}
\end{figure}

Having obtained a prediction for the damaged state, we can plot the resulting velocity field in \cref{fig:honeymoonOptimized}, illustrated as a distribution of velocity magnitudes. These velocity magnitudes confirm that the crack patterns are more developed in the areas where higher mass flows are measured compared to the areas where low mass flows are measured. In particular, in the bottom center a row of hexagons developed cracks with damage parameters around $0.95$. This area is of particular interest, because it is connected to the outside of the well field and could lead to the uranium-rich solution escaping into the groundwater. Such results should prompt the installment of additional monitoring wells in the vicinity of the area. In the top center of the mine, where low mass-flow (under $1.5 \ell$/s) for one hexagon is measured, the subsurface soil is highly impermeable with damage parameters around $0.05$. Such results may instigate hydraulic fracturing measures at these areas.


\section{Conclusions}
\label{sec:Conclusion}

In this article, we have developed a reduced basis method capable of handling the coupled Stokes/Darcy equations on variable internal geometry. We based our model on a diffuse interface representation of the Stokes and Darcy domains, which permits changes in topology of the parameter dependent internal geometry. We showed that all three of the Beavers-Joseph-Saffman coupling conditions can be treated straightforwardly in a diffuse geometry setting, and we highlighted the equivalency with a particular type of Brinkman model.

We used the discrete empirical interpolation method (DEIM) to handle the non-affinity of the phase-field quantities that arise in the coupled equations. To ensure non-negativity of the domain indicator required by physics, we introduced a non-negativity preserving version of DEIM. We showed that in total, we require three such DEIM representations of phase-field related quantities that arise from the diffuse geometry representation and the Beavers-Joseph-Saffman coupling conditions. We studied the performance of the non-negativity preserving DEIM reconstruction for three benchmark problems. The benchmarks differed in complexity, which could also be observed in the convergence graphs for $\epsilon$, a relative error measure depending on the singular values for the DEIM interpolation modes. 

The same three benchmarks were used to investigate the performance of a complete reduced basis formulation. The solution modes that act as reduced basis functions were obtained from the singular value decomposition of a snapshot matrix. We investigated the dependence of the number of solution modes and the number of DEIM interpolation modes for each of the three fields on the maximum $L^2$-error of the reduced order model throughout the parameter space, and we proposed a strategy for determining a required number of DEIM modes for each of the phase-field quantities for a given number of solution modes based on precomputed $\epsilon$-values. Making use of this strategy, we obtained the optimum maximum relative $L^2$-error for a given number of solution modes while significantly reducing the number of DEIM interpolation modes compared to a naive approach. 

We then illustrated the application of our reduced order model for a large-scale computations via the analysis of an in-situ leach mining site. A reduced order model for the subsurface flow throughout the entire mining site was constructed from reduced basis approximations on subdomains and by subsequently coupling the different subdomains with Lagrange multipliers. We showcased the capability of the resulting model for predicting subsurface crack patterns for an existing in-situ leach mining site based on measured inflow and outflow data at injection and production wells. A complete model evaluation required approximately $0.1$ seconds on a desktop machine, and involved solving a system of equation of 6,077 degrees of freedom. This constitutes a significant cost reduction compared to the 15 million degrees of freedom that are required to solve the full problem without the reduced order model at virtually the same fidelity level.

\vspace{0.3cm}

\textbf{Acknowledgements:} 
The results presented in this paper were achieved as part of the ERC Starting Grant project ``ImageToSim'' that has received funding from the European Research Council (ERC) under the European Union’s Horizon 2020 research and innovation programme (Grant agreement No. 759001). The authors gratefully acknowledge this support.

\bibliographystyle{ieeetr}
\bibliography{MyBib.bib}   

\end{document}